\documentclass[a4paper,12pt, reqno]{amsart}
\usepackage{amssymb,amsmath,latexsym}
\usepackage{color}
\usepackage{xcolor}
\usepackage{tikz}
\usepackage{comment}
\allowdisplaybreaks

\def\1{{\bf 1}}
\def\ind{{\bf 1}}
\def\nn{\nonumber}
\def\bee{\begin{equation}}
\def\eee{\end{equation}}
 \def\sB {{\mathcal B}}

\def\R {{\mathbb R}}

\newtheorem{thm}{Theorem}[section]
\newtheorem{lemma}[thm]{Lemma}

\newtheorem{prop}[thm]{Proposition}
\newtheorem{corollary}[thm]{Corollary}
\newtheorem{remark}[thm]{Remark}

\numberwithin{equation}{section}

\def\qed{{\hfill $\Box$ \bigskip}}
\def\NN{{\mathcal N}}

\def\BB{{\mathcal B}}
\def\CC{{\mathcal C}}
\def\LL{{\mathcal L}}
\def\DD{{\mathcal D}}
\def\FF{{\mathcal F}}
\def\EE{{\mathcal E}}
\def\QQ{{\mathcal Q}}

\def\R{{\mathbb R}}

\def\E{{\mathbb E}}

\def\P{{\mathbb P}}
\def\N{{\mathbb N}}
\def\eps{\varepsilon}
\def\wh{\widehat}
\def\wt{\widetilde}
\def\pf{\noindent{\bf Proof.} }

\usepackage{hyperref}
\begin{document}
\title[Green function estimates]{Sharp two-sided Green function estimates for Dirichlet forms degenerate at the boundary}
\author{ Panki Kim \quad Renming Song \quad and \quad Zoran Vondra\v{c}ek}
\thanks{Panki Kim: This work was  supported by the National Research Foundation of
Korea(NRF) grant funded by the Korea government(MSIP) 
(No. 2016R1E1A1A01941893)
}
\thanks{Renming Song: Research supported in part by a grant from
the Simons Foundation (\#429343, Renming Song)}
\thanks{Zoran Vondra\v{c}ek: Research supported in part by the Croatian Science Foundation under the project 4197.}

 \date{}

\begin{abstract}
The goal of this paper is to establish Green function estimates for
 a class of purely discontinuous symmetric Markov processes with jump kernels
degenerate at the boundary and critical killing potentials. The jump kernel and the killing potential 
depend
on several parameters. We establish sharp two-sided estimates on the Green functions of these processes for all admissible values of 
the parameters involved. Depending on the regions where the parameters belong,
the estimates on the Green functions are different. In fact, the estimates have
three different forms.
As applications, we prove that the boundary Harnack principle holds in 
certain region of the parameters
and fails in some other region of the parameters. 
Combined with the main results of \cite{KSV},
we completely  determine the  region of the parameters where the boundary Harnack principle holds.
\end{abstract}
\maketitle

\bigskip



\smallskip
\section{Introduction and main results}
In the last few decades, 
many important results have been obtained 
in the study of potential theoretic properties  for various types of jump processes in open subsets of $\R^d$. These include  isotropic $\alpha$-stable processes,  more general symmetric L\'evy and L\'evy-type processes and their censored versions. 
The main results include the boundary Harnack principle, 
see \cite{Bog97, SW99, BBC, BKK, CS07, KSV14, GK18}, 
sharp two-sided Green function estimates, 
see  \cite{Kul97, CS98, CK02, CT, KSV12b, CKS14}
and sharp two-sided Dirichlet heat kernel estimates,  
see  \cite{BGR10, CKS10a, CKS10b, BGR14, CKS14, KM18, GKK}.
In all these results, the jump kernel $J^D(x,y)$ of the process in the open set $D$ is either the restriction of the jump kernel of the original  process in $\R^d$ or comparable to such a kernel and it does
 not tend to zero as $x$ or $y$ tends to the boundary of $D$.
In this sense, 
one can say that the corresponding  integro-differential operator 
is uniformly elliptic.

Subordinate killed Brownian motions, and more generally, subordinate killed L\'evy 
processes, form another important class of Markov processes.
In case of a stable subordinator, the generator of the subordinate killed Brownian motion is 
the spectral fractional Laplacian. The spectral fractional Laplacian and, more generally, 
fractional powers of  elliptic differential operators in domains have
been studied by 
many people in the PDE community,  see \cite{ST,  BV, CS16, Gru, BFV,BSV}. 
In contrast with killed L\'evy processes and censored processes,  the jump kernel of 
a subordinate killed L\'evy process in an open
subset $D\subset \R^d$ tends to zero near the boundary of $D$, 
see \cite{SV, KSV19, KSV20}.
In this sense, the Dirichlet forms of subordinate killed L\'evy processes are degenerate near the boundary. 
Partial differential equations degenerate at the boundary have been studied 
intensively in the PDE literature, see, for instance, 
\cite{DL, Kim, FMPS, FP, WX} and the references therein.

In our recent paper \cite{KSV}, we introduced a class of symmetric Markov processes in open subsets
$D\subset \R^d$ whose Dirichlet forms are degenerate at the boundary of $D$. This
class of processes includes  subordinate killed L\'evy processes as special cases. 

This paper is the second part of our investigation of the potential theory of Markov processes with 
jump kernels degenerate at the boundary. 
In \cite{KSV} we studied Markov processes in open sets $D\subset \R^d$ defined via Dirichlet
forms with 
jump kernels $J^D(x,y)=j(|x-y|)\sB(x,y)$
(where $j(|x|)$ is the density of a pure jump isotropic L\'evy process)  and 
critical killing potentials $\kappa$. 
The function $\sB(x,y)$ 
is assumed to satisfy certain conditions, and is allowed to decay at the boundary of the state space $D$. 
This is in contrast with all the works mentioned in the first paragraph where $\sB(x,y)$
is assumed to be bounded between two positive constants, which can be viewed as a uniform ellipticity condition for non-local operators. In this sense, our paper \cite{KSV} is the first systematic attempt to study the potential theory of 
general degenerate non-local operators 
defined in terms of  Dirichlet forms.
We proved in  \cite{KSV} that the Harnack inequality and Carleson's estimate are valid for non-negative harmonic functions with respect  
to these Markov processes.

When $D=\R_+^d=\{x=(\wt{x},x_d):\, x_d>0\}$, $j(|x-y|)=|x-y|^{-\alpha-d}$, $0<\alpha<2$,  
and $\kappa(x)=c x_d^{-\alpha}$,
we showed in  \cite{KSV} that for certain values of 
the parameters involved in $\sB(x,y)$ the boundary Harnack principle holds, while for some other values of the parameters the boundary Harnack principle fails (despite the fact that Carleson's estimate holds). 
The main goal of this paper is to establish sharp two-sided estimates on 
the Green functions of the corresponding processes
for all admissible values of the parameters involved  in $\sB(x,y)$. 
These estimates imply anomalous boundary behavior for certain Green potentials, 
see  Proposition \ref{p:bound-for-integral-new}, a feature recently studied
both in the probabilistic as well as in the PDE literature, 
see  \cite{AGV, BFV, KSV20}.
As an application of these Green function estimates, we give a complete answer to the question for 
which values of the parameters the boundary Harnack principle holds true.

We first repeat the assumptions on 
$\sB$ 
that were introduced in \cite{KSV}.
 Here and below, 
 $a\wedge b:=\min \{a, b\}$ and  $a\vee b:=\max\{a, b\}$. 
\medskip

\noindent
\textbf{(A1)} $\sB(x,y)=\sB(y,x)$ for all $x,y\in  \R^d_+$.

\medskip
\noindent
\textbf{(A2)}  
If $\alpha \ge 1$, then there exist $\theta>\alpha-1$ and  
$C_1>0$ such that 
$$
|\sB(x, x)-\sB(x,y)|\le 
C_1\left(\frac{|x-y|}{x_d\wedge y_d}\right)^{\theta}\,.
$$ 

\medskip
\noindent
\textbf{(A3)}
There exist $C_2\ge 1$ and 
parameters $\beta_1, \beta_2, \beta_3,  \beta_4  \ge 0$,  
with $\beta_1>0$ if $\beta_3 >0$, and $\beta_2>0$ if $\beta_4>0$, 
such that
\begin{equation}\label{e:B7}
C_2^{-1}\wt{B}(x,y)\le \sB(x,y)\le C_2 \wt{B}(x,y)\, ,\qquad x,y\in \R^d_+\, ,
\end{equation}
where
\begin{eqnarray}\label{e:B(x,y)}
\wt{B}(x,y)&:=& \Big(\frac{x_d\wedge y_d}{|x-y|}\wedge 1\Big)^{\beta_1}\Big(\frac{x_d\vee y_d}{|x-y|}\wedge 1\Big)^{\beta_2} \left[ \log\Big(1+\frac{(x_d\vee y_d)\wedge |x-y|}{x_d\wedge y_d\wedge |x-y|}\Big)\right]^{\beta_3}\nn \\  
 & &  \times  \left[\log \Big(1+\frac{|x-y|}{(x_d\vee y_d)\wedge |x-y|}\Big)\right]^{\beta_4}. 
\end{eqnarray}

\medskip
\noindent
\textbf{(A4)} 
For all $x,y\in \R^d_+$ and  $a>0$, $\sB(ax,ay)=\sB(x,y)$. 
In case $d\ge 2$, for  all 
$x,y\in \R^d_+$ and $\wt{z}\in \R^{d-1}$, $\sB(x+(\wt{z},0), y+(\wt{z},0))=\sB(x, y)$.

\medskip
Other than the requirements $\beta_1>0$ if $\beta_3 >0$ and $\beta_2>0$ if $\beta_4>0$, the parameters 
$\beta_1, \beta_2, \beta_3$ and  $\beta_4$ are arbitrary.
They control the rate at which $\sB$ goes to 0 at the boundary. Note that the term
$$
 \Big(\frac{x_d\wedge y_d}{|x-y|}\wedge 1\Big)^{\beta_1} \left[ \log\Big(1+\frac{(x_d\vee y_d)\wedge |x-y|}{x_d\wedge y_d\wedge |x-y|}\Big)\right]^{\beta_3}$$
goes to 0 when one of $x$ and $y$ goes to the boundary, while the term
$$
\Big(\frac{x_d\vee y_d}{|x-y|}\wedge 1\Big)^{\beta_2} \left[\log \Big(1+\frac{|x-y|}{(x_d\vee y_d)\wedge |x-y|}\Big)\right]^{\beta_4}
$$
goes to 0 when both $x$ and $y$ go to the boundary.
Note that, if 
$\sB(x,y) \equiv c \wt{B}(x,y)$ for some positive constant $c$, 
then \textbf{(A1)}-\textbf{(A4)} trivially hold.

\medskip
{\it 
In the remainder of this paper, we always assume that}
\begin{align*}
& \quad \quad d > (\alpha+ \beta_1 +\beta_2) \wedge 2,  \quad  
 p\in ((\alpha-1)_+, \alpha+\beta_1)\quad \text{ and } \\
&  J(x,y)=|x-y|^{-d-\alpha}
\sB(x,y) \text{ on }  \R_+^d\times \R_+^d \text{ with  } \sB
\text{ satisfying } \textbf{(A1)}-\textbf{(A4)}.
\end{align*}

\medskip
To every parameter $p\in ((\alpha-1)_+, \alpha+\beta_1)$, we associate a constant 
$C(\alpha, p, \sB)\in (0, \infty)$
depending on $\alpha$, $p$ and $\sB$ 
defined as
\begin{equation}\label{e:explicit-C}
C(\alpha, p, \sB)=
\int_{\R^{d-1}}\frac{1}{(|\wt{u}|^2+1)^{(d+\alpha)/2}}
\int_0^1 \frac{(s^p-1)(1-s^{\alpha-p-1})}{(1-s)^{1+\alpha}}
\sB\big((1-s)\wt{u}, 1), s\mathbf{e}_d \big)\, ds
d\wt{u}\, ,
\end{equation}
where $\mathbf{e}_d=(\tilde{0}, 1)$.
In case $d=1$, $  C(\alpha, p, \sB)$ is defined as
$$
C(\alpha, p, \sB)=\int_0^1 \frac{(s^p-1)(1-s^{\alpha-p-1})}{(1-s)^{1+\alpha}}
\sB\big(1, s \big)\, ds.
$$
Note that $\lim_{p\downarrow (\alpha-1)_+}C(\alpha, p, \sB)=0$, $\lim_{p\uparrow \alpha+\beta_1}C(\alpha, p, \sB)=\infty$ and that the function $p\mapsto C(\alpha, p, \sB)$ is 
strictly increasing
(see \cite[Lemma  5.4  and Remark  5.5]{KSV}). 
Thus, the interval $((\alpha-1)_+, \alpha+\beta_1)$ is the full admissible range for the  parameter $p$. 

Let
\begin{align}
\label{e:kappap}
\kappa(x)
=C(\alpha, p,  \sB)x_d^{-\alpha}, \qquad x\in \R^d_+,
\end{align} 
be the killing potential.
Note that $\kappa$ depends on $p$, but we omit this dependence from the notation for simplicity.
We denote  by $Y$  the Hunt process with jump kernel $J$ and 
killing potential $\kappa$.

To be more precise,
let us define
$$
\EE^{\R^d_+}(u,v):=\frac12 \int_{\R^d_+}\int_{\R^d_+} (u(x)-u(y))(v(x)-v(y))J(x,y)\, dy\, dx,
$$
which is a symmetric form degenerate at the boundary due to \textbf{(A1)} and \textbf{(A3)}.
By Fatou's lemma,  $(\EE^{\R^d_+}, C_c^{\infty}({\R^d_+}))$
is closable in $L^2({\R^d_+}, dx)$. 
Let $\FF^{\R^d_+}$ be the closure of $C_c^{\infty}({\R^d_+})$ under
$\EE^{\R^d_+}_1:=\EE^{\R^d_+}+(\cdot, \cdot)_{L^2({\R^d_+},dx)}$. 
Then $(\EE^{\R^d_+}, \FF^{\R^d_+})$ 
is a regular Dirichlet form on $L^2({\R^d_+}, dx)$. 
Set
$$
\EE(u,v):=\EE^{\R^d_+}(u,v)+\int_{\R^d_+} u(x)v(x)\kappa(x)\, dx\, .
$$
Since $\kappa$ is locally bounded, the measure $\kappa(x)dx$ is a positive Radon measure  charging no set of zero capacity. Let $\FF:=\wt{\FF^{\R^d_+}}\cap L^2({\R^d_+}, \kappa(x)dx)$, where $\wt{\FF^{\R^d_+}}$ is the family of all quasi-continuous functions in $\FF^{\R^d_+}$. By \cite[Theorems 6.1.1 and 6.1.2]{FOT}, 
$(\EE, \FF)$ 
is a regular Dirichlet form on $L^2({\R^d_+}, dx)$ with $C_c^{\infty}({\R^d_+})$ as a special standard core.
Let $((Y_t)_{t\ge 0}, (\P_x)_{x\in {\R^d_+}\setminus \NN})$ be the associated Hunt 
process with lifetime $\zeta$. By \cite[Proposition 3.2]{KSV},  the exceptional set 
$\NN$ can be taken as the empty set. 
We add a cemetery point $\partial$ to the state space ${\R^d_+}$ and define $Y_t=\partial$
for $t\ge \zeta$. 

The process $Y$ enjoys the following important scaling property shown in \cite[Lemma 5.1]{KSV}: For any $r>0$ define the process $Y^{(r)}$ by $Y_t^{(r)}:=rY_{r^{-\alpha}t}$. Then under \textbf{(A1)}, the boundedness of $\sB$ and \textbf{(A4)}, $(Y^{(r)}, \P_{x/r})$ has the same law as $(Y, \P_x)$. The homogeneity property of $\sB$ from \textbf{(A4)} is crucial to establish this fact.

Recall that a  Borel function $f:\R^d_+\to [0, \infty)$ is said to be  \emph{harmonic} in an open set $V\subset \R^d_+$ 
with respect to $Y$ if for every bounded open set $U\subset\overline{U}\subset V$,
\begin{align}
\label{e:har}
f(x)= \E_x \left[ f(Y_{\tau_{U}})\right], \qquad
\hbox{for all } x\in U,
\end{align}
where $\tau_U:=\inf\{t>0:\, Y_t\notin U\}$ is 
the first exit time 
of $Y$ from $U$. 
We say $f$ is \emph{regular harmonic} in $V$ if \eqref{e:har} holds for $V$.

Let $G(x,y)$ denote the Green function of the process $Y$. 
The following theorem is our main result on Green function estimates.
For two functions $f$ and $g$, we use 
the notation $f\asymp g$ to denote that the quotient $f/g$ stays bounded between two positive constants.

\begin{thm}\label{t:Green}
Assume that 
\textbf{(A1)}-\textbf{(A4)} and \eqref{e:kappap} hold true.
Suppose that 
$d > (\alpha+ \beta_1 +\beta_2) \wedge 2$ and $p\in ((\alpha-1)_+, \alpha+\beta_1)$.
Then the process $Y$ admits a Green function $G:\R^d_+\times \R^d_+\to [0,\infty]$
 such that $G(x, \cdot)$ is 
 continuous in $\R^d_+\setminus \{x\}$ and regular harmonic with respect to $Y$ in $\R^d_+\setminus B(x, \epsilon)$ for any $\epsilon>0$. 
Moreover, $G(x,y)$ has the following estimates:

\medskip
\noindent
(1) If
$p\in ((\alpha-1)_+, \alpha+\frac12[\beta_1+(\beta_1 \wedge \beta_2)])$, 
then on $\R^d_+\times \R^d_+$,
\begin{align}\label{e:Green}
G (x,y)\asymp \frac{1}{|x-y|^{d-\alpha}}\left(\frac{x_d}{|x-y|}  \wedge 1 \right)^p\left(\frac{y_d}{|x-y|}  \wedge 1 \right)^p .
\end{align}

\medskip
\noindent
(2) If  
$p=\alpha+\frac{\beta_1+\beta_2}2$, then on $\R^d_+\times \R^d_+$,
$$
G (x,y)\asymp \frac{1}{|x-y|^{d-\alpha}}\left(\frac{x_d}{|x-y|}  \wedge 1 \right)^p\left(\frac{y_d}{|x-y|}  \wedge 1 \right)^p 
 \left( \log\left(1+\frac{|x-y|}{(x_d\vee y_d)\wedge |x-y|}\right) \right)^{ \beta_4+1}. 
$$

\medskip
\noindent
(3)
If $p\in (\alpha+\frac{\beta_1+\beta_2}2, \alpha+\beta_1)$,
then on $\R^d_+\times \R^d_+$,
\begin{align*}
&G(x, y) \\
&\asymp \frac{1}{|x-y|^{d-\alpha}}\left( \frac{ x_d \wedge y_d}{|x-y|} \wedge 1 \right)^p
\left(\frac{ x_d \vee y_d}{|x-y|}\wedge 1\right)^{2\alpha-p+\beta_1+\beta_2} 
\left( \log\left(1+\frac{|x-y|}{(x_d\vee y_d)\wedge |x-y|}\right) \right)^{\beta_4} 
\\
& =
\frac{1}{|x-y|^{d-\alpha}}\left( \frac{ x_d }{|x-y|} \wedge 1 \right)^p\left( \frac{  y_d}{|x-y|} \wedge 1 \right)^p
\left(\frac{ x_d \vee y_d}{|x-y|}\wedge 1\right)^{-2(p-\alpha-(\beta_1+ \beta_2)/2)} \\
& \quad \times  
 \left( \log\left(1+\frac{|x-y|}{(x_d\vee y_d)\wedge |x-y|}\right) \right)^{\beta_4}.  
\end{align*}
\end{thm}

\bigskip
Note that when $\beta_1\le \beta_2$, then case (1) covers all possible values of the parameter $p$, while when $\beta_2<\beta_1$ 
the regimes of $p$ in cases (1), (2) and (3) are disjoint and exhaustive.

In fact, for lower bounds of Green functions, we have more general results, see Theorems \ref{t:GB} and \ref{t:f4e}.
In these theorems, we establish
lower bounds on the Green function 
$G^{B(w, R)\cap \R^d_+}(x, y)$ for  $Y$ killed upon exiting $B(w, R)\cap \R^d_+$ (where $w \in \partial \R^d_+$) in $B(w, (1-\eps)R)\cap \R^d_+$.
The lower bounds 
on $G(x, y)$ in the theorem above 
are corollaries of these more general results.

Note that 
$$
p \mapsto 2\alpha-p+\beta_1+\beta_2 =(\alpha+\beta_2)+(\alpha+ \beta_1-p)
$$ 
is decreasing on $\alpha+\frac{\beta_1+\beta_2}2 \le p <\alpha+\beta_1$, 
which has a somewhat strange and interesting consequence. Namely, 
the power of $\frac{ x_d \wedge y_d}{|x-y|}\wedge 1$ is always $p$ 
and we can increase the exponent $p$
of $\frac{ x_d \wedge y_d}{|x-y|}\wedge 1$ all the way up to (just below) $\alpha+ \beta_1$. 
But the exponent  of $\frac{ x_d \vee y_d}{|x-y|}\wedge 1$ is $p$
only up to  $\alpha+\frac{\beta_1+[\beta_1 \wedge \beta_2]}2$ 
and one can increase the exponent only up to  $\alpha+\frac{\beta_1+[\beta_1 \wedge \beta_2]}2$.
In the case $\beta_2<\beta_1$, once $p$ reaches $\alpha+\frac{\beta_1+\beta_2}2$,
the exponent of $\frac{ x_d \vee y_d}{|x-y|}\wedge 1$ starts decreasing.

Estimates \eqref{e:Green} can be equivalently stated as
\begin{equation}\label{e:Green1}
G (x,y) \asymp \left(\frac{x_dy_d}{|x-y|^2}  \wedge 1 \right)^p 
\frac{1}{|x-y|^{d-\alpha}} \quad \mbox{ on } \R^d_+\times \R^d_+.
\end{equation}
Note that, when $d \ge 3$ and  $p=(d-\alpha)/(d-2) \in (1, d/(d-2))$, the estimates in \eqref{e:Green1} are those of a power 
of the Green function of killed Brownian motion in $\R^d_+$.
See  \cite{DDMMV}.

Moreover, we can rewrite the estimates in Theorem \ref{t:Green}  in a unified way:
 Let $a_p=2(p-\alpha-\frac{\beta_1+[\beta_1 \wedge \beta_2]}2)$. Then on $\R^d_+\times \R^d_+$,
\begin{align*}
&G(x, y) \asymp \\
 & \frac{1}{|x-y|^{d-\alpha}}\left( \frac{ x_d \wedge y_d}{|x-y|} \wedge 1 \right)^p\left(\frac{ x_d \vee y_d}{|x-y|}\wedge 1\right)^{p-{a_p}_+} 
 \log\left(2+{\bf 1}_{a_p\le 0}\frac{|x-y|}{(x_d \vee y_d)\wedge |x-y|}\right)^{
 \beta_4 +{\bf 1}_{a_p= 0}}.
\end{align*}

\bigskip

In \cite[Theorem 1.3]{KSV} we have proved that the boundary Harnack principle holds 
when either (a) $\beta_1=\beta_2$ and $\beta_3 =\beta_4  =0$, or (b) $p<\alpha$.
 In \cite[Theorem 1.4]{KSV} we have showed that when $\alpha+\beta_2<p<\alpha+\beta_1$
 the boundary Harnack principle fails. However, we were unable to determine what happens 
 with the boundary Harnack principle in the remaining regions of the admissible parameters.
 As applications of our Green function estimates, 
 we can completely resolve this issue 
 and prove the following two results. 
In the remainder of this paper, 
we will only give the statements and proofs of the results for $d\ge 2$. 
The counterparts in the $d=1$ case are similar and simpler.
 
For any $a, b>0$ and $\wt{w} \in \R^{d-1}$, we define a box
 $$D_{\wt{w}}(a,b):=\{x=(\wt{x}, x_d)\in \R^d:\, |\wt{x}-\wt{w}|<a, 0<x_d<b\}.$$ 
 
\begin{thm}\label{t:BHPnew}
Assume that 
\textbf{(A1)}-\textbf{(A4)} and \eqref{e:kappap} hold true.
Suppose that 
$d > (\alpha+ \beta_1 +\beta_2) \wedge 2$ and $p\in ((\alpha-1)_+, \alpha+(\beta_1 \wedge \beta_2))$.
Then there exists $C_3 \ge 1$
 such that for all $r>0$, $\wt{w} \in \R^{d-1}$, 
and any non-negative function $f$ in $\R^d_+$ which is harmonic in $D_{\wt{w}}(2r, 2r)$ with respect to $Y$ and vanishes continuously on 
$B(({\wt{w}},0), 2r)\cap \partial \R^d_+$, 
we have 
\begin{equation}\label{e:TAMSe1.8new}
\frac{f(x)}{x^p_d}\le C_3\frac{f(y)}{y^p_d}, 
\quad x, y\in D_{\wt{w}}(r/2, r/2) .
\end{equation}
\end{thm}

\bigskip

Theorem \ref{t:BHPnew} implies that, 
if two functions $f, g$ in $\R^d_+$ both satisfy the assumptions in 
Theorem \ref{t:BHPnew}, then
$$
\frac{f(x)}{f(y)}\,\le C_3^2\,\frac{g(x)}{g(y)}, \quad x, y\in D_{\wt{w}}(r/2, r/2) .
$$

\medskip
We say that the non-scale-invariant boundary Harnack principle holds
near the boundary of $\R^d_+$ if there is a constant $\wh{R}\in (0,1)$  such that for any 
$r \in (0, \wh{R}\, ]$, there exists  a constant $c=c(r)\ge 1$ such that for all 
$\wt w\in \R^{d-1}$ and non-negative functions $f, g$ in $\R^d_+$ which are harmonic  
in $\R^d_+ \cap B((\wt w, 0),r)$ with respect to $Y$ and vanish continuously on $ \partial \R^d_+ \cap B((\wt w, 0), r)$, we have
$$
\frac{f(x)}{f(y)}\,\le c\,\frac{g(x)}{g(y)} \qquad
\hbox{for all } x, y\in  
B((\wt w, 0), r/2)\cap \R^d_+.
$$
\medskip

\begin{thm}\label{t:counterexample}
Suppose  
$d > \alpha+ \beta_1 +\beta_2$
and $d \ge 2$.
Assume that \textbf{(A1)}-\textbf{(A4)} and \eqref{e:kappap} hold true.
If $\alpha+\beta_2 \le p<\alpha+\beta_1$, then
the non-scale-invariant boundary Harnack principle is not valid for $Y$. 
\end{thm}
Thus, when $\alpha+\beta_2 \le p<\alpha+(\beta_1+\beta_2)/2$, the boundary Harnack principle is not valid for $Y$ even though we have 
the standard form of the Green function
estimates \eqref{e:Green1}. This phenomenon has already been observed by 
the authors in \cite{KSV20} for subordinate killed L\'evy processes.

The following two results proved in \cite{KSV}  will be fundamental for this paper.
Note that, by the scaling property of $Y$
(see \cite[Lemma 5.1]{KSV}), 
we can allow $r>0$ instead of $r\in (0,1]$.
\begin{thm} 
 [{Harnack inequality, \cite[Theorem 1.1]{KSV}}]
 \label{t:uhp}
Assume that 
\textbf{(A1)}-\textbf{(A4)} and \eqref{e:kappap} hold true and $p\in ((\alpha-1)_+, \alpha+\beta_1)$.
\begin{itemize}
\item[(a)] 
There exists a constant $C_4>0$ such that for any
$r>0$, any $B(x_0, r) \subset \R^d_+$ and any non-negative function $f$ in $\R^d_+$ which is 
harmonic in $B(x_0, r)$ with respect to $Y$, we have
$$
f(x)\le C_4 f(y), \qquad \text{ for all } x, y\in  B(x_0, r/2).
$$

\item[(b)] 
There exists a constant $C_5>0$ 
such that for any  $L>0$, any  $r>0$,
any $x_1,x_2 \in \R^d_+$ with $|x_1-x_2|<Lr$ and $B(x_1,r)\cup B(x_2,r) \subset \R^d_+$ 
and any non-negative function $f$ in $\R^d_+$ which is  harmonic in $B(x_1,r)\cup B(x_2,r)$ with respect to $Y$, we have
$$
f(x_2)\le C_5
(L+1)^{\beta_1+\beta_2+d+\alpha} f(x_1)\, .
$$
\end{itemize}
\end{thm}

Since the half-space $\R^d_+$ is $\kappa$-fat with characteristics $(R, 1/2)$ for any $R>0$, we also have

\begin{thm}
[{Carleson's estimate, \cite[Theorem 1.2]{KSV}}]
\label{t:carleson}
Assume that 
\textbf{(A1)}-\textbf{(A4)} and \eqref{e:kappap} hold true and $p\in ((\alpha-1)_+, \alpha+\beta_1)$. Then there exists a constant $C_6>0$  
such that for any $w \in\partial \R^d_+$, $r>0$, 
and any non-negative  function $f$ in $\R^d_+$ that is harmonic in $\R^d_+ \cap B(w, r)$ with respect to $Y$ and vanishes continuously on $ \partial \R^d_+ \cap B(w, r)$, we have
\begin{equation}\label{e:carleson}
f(x)\le C_6 f(\wh{x}) \qquad \hbox{for all }  x\in \R^d_+\cap B(w,r/2),
\end{equation}
where $\wh{x}\in \R^d_+\cap B(w,r)$ with $\wh{x}_d\ge r/4$.
\end{thm}
 
The assumptions 
\textbf{(A1)}, \textbf{(A2)}, \textbf{(A3)} and \textbf{(A4)} in this paper are the assumptions \textbf{(B1)}, \textbf{(B4)}, \textbf{(B7)} and \textbf{(B8)} in \cite{KSV}, respectively. 
As a consequence of assumptions 
\textbf{(A1)}-\textbf{(A4)},
$\sB(x,y)$ also satisfies assumptions \textbf{(B2)}, \textbf{(B3)}, \textbf{(B5)} and \textbf{(B6)} in \cite{KSV}. 

\medskip

Now we explain the content of this paper and our strategy for proving the main results.

In Section \ref{s:existence-GF} we first show that the process  
$Y$ is transient and 
admits a symmetric Green function $G(x,y)$,  
see Proposition \ref{p:existenceGF}. 
This is  quite
standard once we establish that the occupation measure $G(x, \cdot)$ of $Y$ is absolutely continuous. 
We also show that $x\mapsto G(x,y)$ is  harmonic away from $y$. 
As a consequence of the scaling property of $Y$ and the invariance property of the half space under scaling, one 
gets the following scaling property of the Green function: For all $x,y\in \R^d_+$,
$$
G(x,y)=|x-y|^{\alpha-d}G\left(\frac{x}{|x-y|}, \frac{y}{|x-y|}\right).
$$
In this paper, we use this property several times  so that,
to prove Theorem \ref{t:Green},  we 
mainly deal with the case of $x,y\in \R^d_+$ satisfying $|x-y| \asymp 1$.

In Section \ref{s:decayofGf}, we show that the Green function $G(x, y)$ tends to $0$ when $x$ or $y$ tends to the boundary. The proof of this result depends in a fundamental way on several lemmas from \cite{KSV}. The decay of the Green function at the boundary allows us to apply Theorem \ref{t:carleson} in later sections.

Section \ref{s:interior} is devoted to proving interior estimates on the Green function $G(x, y)$. Roughly, we show that if the points $x,y\in \R^d_+$ are closer to each other than to the boundary, then $G(x,y)\asymp |x-y|^{-d+\alpha}$.  
For the lower bound given in Proposition \ref{p:green-lower-bound}, we use a capacity argument. 
The upper bound is more difficult and relies on 
the Hardy inequality in 
\cite[Corollary 3]{BDK}  
and the heat kernel estimates of symmetric
jump processes with 
large jumps of lower intensity in \cite{BKKL}.
 This is where the assumption 
$d > (\alpha+ \beta_1 +\beta_2) \wedge 2$
is needed. 
The key to obtaining the interior upper estimate is to get a uniform estimate on the $L^2$ norm of $\int_{B(z, 4)}G(x,y)dy$ on $B(z, 4)$
for all $z$
sufficiently away from the boundary, 
see Proposition \ref{p:l2norm}.

In Section \ref{s:prelgfe}, we give a lower bound for the Green function of 
the process $Y$ killed upon exiting a half-ball centered at the boundary of $\R^d_+$
and a preliminary upper bound for the Green function. 
The lower bound given in Theorem \ref{t:GB} is proved for $G^{B(w,R)\cap \R^d_+}(x,y)$, 
the Green function of the process $Y$ killed upon exiting 
$B(w,R)\cap \R^d_+$, $w\in \partial \R^d_+$, for $x,y\in B(w, (1-\epsilon)R  
)\cap \R^d_+$. 
This gives 
the sharp lower bound of Green function for $p\in ((\alpha-1)_+, \alpha+\frac12[\beta_1+(\beta_1 \wedge \beta_2)])$.
A preliminary estimate of the upper bound is given in Lemma \ref{l:prelub}. Proofs of these estimates use 
the already mentioned fundamental lemmas from \cite{KSV} and Theorem \ref{t:carleson}.

Section \ref{s:thm1.1} is central to the paper. We first prove a technical Lemma \ref{l:key} modeled after \cite[Lemma 3.3]{AGV} and its Corollary \ref{c:key}. 
They are both used throughout this section. In proving 
Theorem \ref{t:Green}, one is led to double integrals involving the Green function (or the Green 
function of the killed process) twice and the jump kernel. 
The sharp bounds of these double integrals are essential in the proof of Theorem \ref{t:Green}. 
To obtain the correct bound,  we have  to divide the region of integration into several parts and deal with them separately. 
These estimates are 
quite difficult and delicate,
see Remark \ref{r:6.8} below. 
By using the preliminary estimates of the Green function obtained in Section \ref{s:prelgfe}
and the explicit form of $\wt{B}$, those integrals are successfully estimated by means of Lemma \ref{l:key} and Corollary \ref{c:key}. As an application of the Green function estimates, we end the section with sharp two-sided estimates on some killed potentials of the process $Y$, or in analytical language, with estimates of $\int_D G^D(x,y)y_d^{\beta}dy$ where $D$ is a box of arbitrary size and 
$\beta>-p-1$
(see Proposition \ref{p:bound-for-integral-new} below), 
as well as estimates of $\int_{\R^d_+}G(x,y)y_d^{\beta}dy$. The latter estimates give precise information on the expected lifetime of the process $Y$. 

In Section \ref{s:BHP} we prove Theorems \ref{t:BHPnew} and \ref{t:counterexample}. 
Our Proposition \ref{p:bound-for-integral-new} 
is powerful enough for us  to cover  the full range of the parameters.

We end this introduction by discussing 
some examples of  explicit processes 
satisfying our assumptions, as well as a process which does not fall in the class considered here.

The first (and the motivating) example is a subordinate killed stable process $Y$ whose infinitesimal generator is
$-((-\Delta)^{\delta/2}_{\ |\R^d_+})^{\gamma/2}$, where $\delta\in (0,2]$  and $\gamma\in (0,2)$. Its jump kernel is $J(x,y)=|x-y|^{-d-\alpha}\sB(x,y)$ with $\alpha=\gamma\delta/2$ and $\sB(x,y)$ satisfying \textbf{(A3)} with parameters as follows: If $\delta=2$ , then $\beta_1=\beta_2=1$, $\beta_3=\beta_4=0$. For $\delta\in (0,2)$,  (i) when $\gamma\in (1,2)$, then $\beta_1=\delta(1-\gamma/2)$, $\beta_2=\beta_3=\beta_4=0$; (ii) when $\gamma=1$, then $\beta_1=\delta/2$, $\beta_3=0$, $\beta_2=\beta_4=0$, (iii) when $\gamma\in (0,1)$, 
then $\beta_1=\delta/2$, 
$\beta_2=(1-\gamma)\delta/2$, $\beta_3=\beta_4=0$.
For more details see \cite[(1.1), (1.2) and Section 2]{KSV}. In all cases it holds that $p=\delta/2$ which can be deduced by comparing Green function estimates in Theorem \ref{t:Green} and \cite[Theorem 6.4]{KSV20}.

An example of a process with $\beta_4>0$ has been recently discovered in \cite{CKSV}. Let $\delta\in (0, 2)$, and let $X$ be the 
reflected symmetric $\delta$-stable process in 
$\overline \R^d_+=\{x=(\wt{x},x_d):\, x_d \ge0\}$ killed leaving upon $\R^d_+$, 
whose infinitesimal generator is the regional fractional Laplacian
$$
\LL f(x)=c(d, \delta) \lim_{\varepsilon \to 0} \int_{\R^d_+, |y-x|> \varepsilon}(f(y)-f(x))|x-y|^{-d-\delta}\, dy,
$$
see \cite[pp. 232--234]{CKSV} for details.
Let $q\in [\delta-1, \delta)\cap(0, \delta)$ and $Z$ be the process corresponding to the the Feynman-Kac semigroup via the
multiplicative functional
$$
\exp\left(-C(d, \delta, q)\int^t_0(X^d_s)^{-\delta}ds\right), 
$$
where $C(d, \delta, q)$ is the positive constant (involving parameter $q$)
defined in \cite[p. 233]{CKSV0},  see also \cite[(3.5)]{CKSV0}.
Let    $S$ be an independent 
$\gamma/2$-stable
subordinator with $\gamma\in (0, 2)$ and set $\alpha=\delta\gamma/2$.
 Define a process $Y$ by $Y_t=Z_{S_t}$ whose infinitesimal generator is
$-(-\LL+ C(d, \delta, q) (x_d)^{-\delta} )^{\gamma/2}$. 
The jump kernel of $Y$ is of the form $J(x,y)=|x-y|^{-d-\alpha}\sB(x,y)$, 
with $\sB(x,y)$ satisfying \textbf{(A1)}--\textbf{(A4)}. 
Moreover, the parameter $\beta_4$ in \textbf{(A3)} is equal to $1$ for 
certain value of $q$.
 For details, see \cite[Example 7.3]{CKSV} and the paragraph above \cite[Lemma 2.2]{KSV}.
 
The jump kernels of this paper are degenerate since $\sB$ approaches 0 at the boundary. There exist processes in $\R^d_+$ whose jump kernels are of the form
$|x-y|^{-d-\alpha}\sB(x, y)$ with $\sB$ blowing up at the boundary. Here is an example.  Let $X$ be an isotopic $\alpha$-stable process, and let $Y$ be the process obtained from $X$ by deleting the parts of the path outside $\R^d_+$. More precisely, let
$$
A_t=\int^t_01_{\{X_s\in \R^d_+\}}ds, 
$$
be the occupation time in $\R^d_+$ up to time $t$ and let $\gamma_t=\inf\{s>0: A_s>t\}$. The process $Y$ defined by $Y_t=X_{\gamma_t}$ is the trace of $X$
in $\R^d_+$.  It is called the path-censored $\alpha$-stable process 
in \cite{KPW}.  Using \cite[Theorem 6.1]{BGPR21}, one can show that the jump kernel of $Y$ is of the form $|x-y|^{-d-\alpha}\sB(x, y)$ with $\sB$ blowing up at the boundary. 

\textbf{Notation:} Throughout this paper,
the positive constants $\beta_1$, $\beta_2$, $\beta_3$, $\beta_4$,  
$\theta$ will remain the same.
We will use the following convention:
Capital letters $C, C_i, i=1,2,  \dots$ will denote constants
in the statements of results and assumptions. The labeling of these constants will remain the same. Lower case letters 
$c, c_i, i=1,2,  \dots$ are used to denote constants in the proofs
and the labeling of these constants starts anew in each proof.
The notation $c_i=c_i(a,b,c,\ldots)$, $i=0,1,2,  \dots$ indicates  constants depending on $a, b, c, \ldots$.
We will use ``$:=$" to denote a
definition, which is read as ``is defined to be".
For any $x\in \R^d$ and $r>0$, we use $B(x, r)$ to denote the open ball of radius $r$ centered at $x$.
For a Borel subset $V$ in $\R^d$, $|V|$ denotes  the Lebesgue measure of $V$ in $\R^d$,
$\delta_{V}:=\mathrm{dist}(V, \partial D)$.
We use the superscript instead of the subscript for the coordinate of processes as $Y=(Y^1, \dots, Y^d)$.


\section{Existence of the Green function}\label{s:existence-GF}

Recall that  $\zeta$ is the  lifetime of $Y$.
Let $f:\R^d_+\to [0,\infty)$ be a Borel function and $\lambda \ge 0$. 
The $\lambda$-potential of $f$ is defined by
$$
G_{\lambda} f(x):=\E_x\int_0^\zeta e^{-\lambda t}f(Y_t)\, dt\, , \quad x\in \R^d_+.
$$
When $\lambda=0$, we write $Gf$ instead of $G_0f$ and call $Gf$ the Green potential of $f$. If $g:\R^d_+\to [0,\infty)$  is another Borel function, then by the symmetry of $Y$ we have that 
\begin{equation}\label{e:symmetry-of-G}
\int_{\R^d_+}G_{\lambda} f(x)g(x)\, dx = \int_{\R^d_+}f(x) G_{\lambda} g(x)\, dx \, .
\end{equation}
For $A\in \BB( \R^d_+)$, we let $G_{\lambda}(x,A):=G_{\lambda} \ind_A (x)$ be the $\lambda$-occupation measure of $A$. 
In this section we show the existence of the Green function of the process $Y$, that is, the density of the $0$-occupation measure. 
We start by recalling some of the results of \cite[Subsection 3.1]{KSV}. 

Let $U$ be a relatively compact $C^{1,1}$ open subset of $\R^d_+$. 
For $\gamma>0$ small enough, 
define a kernel $J_{\gamma}(x,y)$ on $\R^d\times \R^d$ by $J_{\gamma}(x,y)=J(x,y)$ for $x,y\in U$, and 
$J_{\gamma}(x,y)= \gamma|x-y|^{-d-\alpha}$
otherwise. 
Then there exist $c_1>0$ and $c_2>0$  such that 
(see  the first display below \cite[(3.3)]{KSV})
$$
c_1|x-y|^{-d-\alpha}\le J_{\gamma}(x,y)\le c_2|x-y|^{-d-\alpha}\, ,\quad x,y\in \R^d\, .
$$
For $u\in L^2(\R^d, dx)$, define
$$
\CC(u,u):=\frac12 \int_{\R^d}\int_{\R^d} (u(x)-u(y))^2 J_{\gamma}(x,y)\, dx\, dy \text{ and }
\DD(\CC):=\{u\in L^2(\R^d):\, \CC(u,u)<\infty\}\, .
$$
Then there exists  a conservative Feller and strongly Feller process $Z$ associated with 
$(\CC, \DD(\CC))$ 
which has a continuous transition density 
(with respect to the Lebesgue measure), 
see  \cite{CK08}.
Let $Z^U$ be the process $Z$ killed upon exiting $U$ and let $A_t:=\int_0^t\wt{\kappa}(Z_s^U)\, ds$ where $\wt{\kappa}$ is a certain non-negative function defined in \cite[Subsection 3.1]{KSV} 
($\wt{\kappa}$ is non-negative when $\gamma>0$ is small enough). 
Let $Y^U$ be the process $Y$ killed upon exiting $U$, and let $(Q_t^U)_{t\ge 0}$ denote its semigroup: For $f:U\to [0, \infty)$,
$$
Q_t^U f(x)=\E_x [f(Y_t^U)]=\E_x[f(Y_t), t<\tau_U],
$$
where $\tau_U=\inf\{t>0:\, Y_t\notin U\}$ is the first exit time from $U$.  It is shown in \cite[Subsection 3.1]{KSV} that 
$$
Q_t^U f(x)=\E_x[\exp(-A_t)f(Z_t^U)]\, , \quad t>0, \ x\in U.
$$
Moreover,  $Q_t^U$ has a  transition density $q^U(t,x,y)$ (with respect to the Lebesgue measure) which is symmetric in $x$ and $y$, and such that for all $y\in U$, $(t,x)\mapsto q^U(t,x,y)$ is continuous. 

Let $G_{\lambda}^U f(x):=\int_0^{\infty}e^{-\lambda t}Q_t^U f(x)\, dt=\E_x \int_0^{\tau_U}e^{-\lambda t} f(Y_t)\, dt$ denote the $\lambda$-potential of $Y^U$ and  $G_{\lambda}^U(x,y):=\int_0^{\infty}e^{-\lambda t}q^U(t,x,y)\, dt$ the $\lambda$-potential density of $Y^U$. 
We will write $G^U$ for $G_{0}^U$ for simplicity.
Then 
$G_{\lambda}^U(x,\cdot)$ is the density of the $\lambda$-occupation measure. In particular this shows that $G_{\lambda}^U(x,\cdot)$ is absolutely continuous with respect to the Lebesgue measure.
Moreover, since $x\mapsto q^U(t,x,y)$ is continuous, we see that $x\mapsto G^U_{\lambda}(x,y)$ is lower semi-continuous. By Fatou's lemma this implies that $G_{\lambda}^U f$ is also lower semi-continuous.

Let $(U_n)_{n\ge 1}$ be a sequence of bounded $C^{1,1}$ open sets such that $U_n\subset \overline{U_n}\subset U_{n+1}$ and $\cup_{n\ge 1}U_n=\R^d_+$. 
For any Borel $f: \R^d_+\to [0,\infty)$, it holds that

\begin{equation}\label{e:Gf-increasing-limit}
G_{\lambda}f(x)=\E_x \int_0^{\zeta}e^{-\lambda t}f(Y_t)\, dt =\, \uparrow\!\!\! \lim_{n\to \infty}\E_x \int_0^{\tau_{U_n}}e^{-\lambda t}f(Y_t)\, dt =\, \uparrow\!\!\!\lim_{n\to \infty}G_{\lambda}^{U_n}f(x)\, ,
\end{equation}
where $\uparrow\!\!\!\ \lim $ 
denotes an increasing limit. 
In particular, if 
$A\in \BB( \R^d_+)$ is of Lebesgue measure zero, then for every $x\in \R^d_+$,
$$
G_{\lambda}(x,A)=\lim_{n\to \infty}G_{\lambda}^{U_n}(x,A)=\lim_{n\to \infty}G_{\lambda}^{U_n}(x,A\cap U_n)=0\, .
$$
Thus,  $G_{\lambda}(x,\cdot)$ is absolutely continuous with respect to the Lebesgue measure 
for each $\lambda\ge 0$ and $x\in \R^d_+$. Together with \eqref{e:symmetry-of-G} 
this shows that the conditions of \cite[VI Theorem (1.4)]{BG} are 
satisfied, which  implies 
that the resolvent $(G_{\lambda})_{\lambda>0}$ is self dual. In particular, 
see \cite[pp.256--257]{BG}, 
there exists a symmetric function $G(x,y)$ excessive in both variables such that
$$
Gf(x)=\int_{\R^d_+}G(x,y)f(y)\, dy\, ,\quad x\in \R^d_+.
$$
We recall, see \cite[II, Definition (2.1)]{BG}, that a measurable function $f:\R^d_+\to [0,\infty)$ is $\lambda$-excessive, $\lambda \ge 0$, with respect to the process $Y$ if for every $t\ge 0$ it holds that 
$\E_x[e^{-\lambda t}Y_t]\le f(x)$ and $\lim_{t\to 0}\E_x[e^{-\lambda t}Y_t] =f(x)$, for every $x\in \R^d_+$.
$0$-excessive functions are simply called excessive functions.

We now show that $Y$ is transient. 

\begin{lemma}\label{l:transience}
The process $Y$ is transient in the sense that there exists $f:\R^d_+\to (0, \infty)$ such that $Gf<\infty$. More precisely, $G\kappa\le 1$.
\end{lemma}
\pf Let $(Q_t)_{t\ge 0}$ denote the semigroup of $Y$. For any $A\in \BB( \R^d_+)$, we use \cite[(4.5.6)]{FOT} with $h=\ind_A$, $f=1$, and let $t\to \infty$ to obtain
$$
\E_{\ind_A dx}(\zeta <\infty) \ge 
\E_{\ind_A dx}(Y_{\zeta-} \in \R_+^d,  \zeta <\infty)
=\int_0^{\infty}\int_{\R^d_+} \kappa(x) Q_s \ind_A(x)\, dx \, dt.
$$
This can be rewritten as
$$
\int_A
 \P_x(\zeta<\infty)
\, dx \ge \int_{\R^d_+} \kappa(x) G\ind_A (x)\, dx =\int_A G\kappa(x)\, dx.
$$
Since this inequality holds for every $A\in \BB( \R^d_+)$, we conclude that 
$
\P_x( \zeta<\infty)
 \ge G\kappa(x)$ for a.e.~$x\in \R^d_+$. 
Both functions $x\mapsto \P_x(\zeta<\infty)$ and $G\kappa$ are excessive. Since $G(x, \cdot)$ is absolutely continuous with respect to the Lebesgue measure 
(i.e.,~Hypotesis (L) holds, see \cite[p.112]{CW}), 
by  \cite[Proposition 9, p.113]{CW}, 
we conclude that  $G\kappa(x) \le \P_x(\zeta <\infty)\le 1$ for all $x\in \R^d_+$. 
\qed

As a consequence of Lemma \ref{l:transience},
we have that  $G(x,y)<\infty$ for a.e.~$y\in \R^d_+$.
Another consequence is that, for every compact $K\subset \R^d_+$,  $G\1_K$ is bounded. Indeed, by the definition of $\kappa$, we see that $\inf_K \kappa(x)=:c_K>0$. Thus 
\begin{equation}\label{e:Green-potential-bounded}
G\1_K\le c_K^{-1}G\kappa \le c_K^{-1}.
\end{equation}

Note that it follows from \eqref{e:Gf-increasing-limit} that, for every non-negative Borel $f$, $G_{\lambda}f$ is lower semi-continuous, as an increasing limit of lower semi-continuous functions. Since every $\lambda$-excessive function is an increasing limit
of $\lambda$-potentials,   
see  \cite[II Proposition (2.6)]{BG}, 
we conclude that all $\lambda$-excessive functions of $Y$ are lower semi-continuous. In particular, for every $y\in \R^d_+$, $G_{\lambda}(\cdot, y)$ is lower semi-continuous. Since $G(\cdot, y)$ is the increasing limit of $G_{\lambda}(\cdot, y)$ as $\lambda\to 0$, we see that $G(\cdot, y)$ is also lower semi-continuous.

Fix an open set $B$ in $\R^d_+$ and $x \in \R^d_+$ and 
let $f$ be a non-negative Borel function on $\R^d_+$. By Hunt's switching identity, \cite[VI, Theorem (1.16)]{BG}, 
$$
\E_x [G f(Y_{\tau_B})]=\int_{\R^d_+} \E_x [G(Y_{\tau_B},y)] f(y)\, dy 
=\int_{\R^d_+} \E_y [G(x,Y_{\tau_B})] f(y)\, dy. 
$$
Suppose, further, that $f=0$ on $B$. Then by the strong Markov property, 
\cite[I, Definition (8.1)]{BG},
$$\int_{\R^d_+} G(x,y) f(y)\, dy
= \E_x \int_{\tau_B}^{\infty} f(Y_t)\, dt 
= \E_x[G f(Y_{\tau_B})]=\int_{\R^d_+ \setminus B} \E_y [G(x,Y_{\tau_B})] f(y)\, dy\, ,$$
and hence $G(x,y)=\E_y [G(x,Y_{\tau_B})]$ for a.e. $y\in \R^d_+ \setminus B$. Since both sides are excessive (and thus excessive for the killed process $Y^{\R^d_+\setminus B}$), equality holds for every $y\in \R^d_+ \setminus B$. By using Hunt's switching identity one more time, we arrive at
$$
G(x,y)=\E_x [G (Y_{\tau_B},y)]\, ,\quad \text{for all } x\in \R^d_+,\ y\in
 \R^d_+\setminus B\, .
$$
In particular, if $y\in \R^d_+\setminus B$ is fixed, then the above equality says that $x\mapsto G(x,y)$ is regular harmonic in $B $ with respect to $Y$. By symmetry, $y\mapsto G(x,y)$ is regular harmonic in $B $ as well. By the Harnack inequality, Theorem \ref{t:uhp}, we conclude that $G(x,y)<\infty$ for all $y\in \R^d\setminus \{x\}$. 
This proves  the following 
proposition on the existence of the Green function.

\begin{prop}\label{p:existenceGF}
There exists a symmetric function $G:\R^d_+\times \R^d_+\to [0,\infty]$ which is lower semi-continuous in each variable and finite outside the diagonal such that for every non-negative Borel $f$,
$$
Gf(x)=\int_{\R^d_+}G(x,y)f(y)\, dy\, .
$$
Moreover, $G(x, \cdot)$ is harmonic with respect to $Y$ in $\R^d_+\setminus \{x\}$
and regular harmonic with respect to $Y$ in $\R^d_+\setminus B(x, \epsilon)$ for any $\epsilon>0$
\end{prop}

\begin{remark}
{\rm
We note in passing that all the results established above are valid, 
with the same proofs,  for the process 
 in any open set $D$ (not necessarily $\R^d_+$), under conditions  (1.3)-(1.6) and  \textbf{(B1)}-\textbf{(B3)} from \cite{KSV}. In particular, in the setup of \cite{KSV}, the process  in any open set $D$
 studied there has a symmetric Green function.
}
\end{remark}

For further use, we recall now the formula for the Green function of the process $Y$ killed upon exiting an open set $B\subset \R^d_+$. Let $f:\R^d\to [0,\infty]$ be a measurable function vanishing on $\R^d_+\setminus B$. By the strong Markov property, for $x\in B$,
\begin{eqnarray*}
\int_{\R^d}G(x,y)f(y)dy &=&\E_x \int_0^{\infty}f(Y_s)ds =\E_x\int_0^{\tau_B}f(Y_s) ds+\E_x\left(\E_{Y_{\tau_B}}\int_0^{\infty} f(Y_s)ds\right)\\
&=&\E_x \int_0^{\infty}f(Y^B_s )ds +\E_x Gf(Y_{\tau_B})\\
&=&\E_x \int_0^{\infty}f(Y^B_s )ds+ \int_{\R^d_+}\E_x [G(Y_{\tau_B},y)]f(y)dy.
\end{eqnarray*}
By rearranging, we see that
\begin{equation}\label{e:gf-formula}
G^B(x,y):=G(x,y)-\E_x [G(Y_{\tau_B},y)]
\end{equation}
is the Green function of $Y^B$.

We end this section with the scaling property of the Green 
function, which will be used several times later in this paper.

\begin{prop}\label{p:scaling-of-G}
For all $x,y\in \R^d_+$, $x\neq y$,  it holds that
\begin{equation}\label{e:scaling-of-G}
G(x,y)=G\left(\frac{x}{|x-y|}, \frac{y}{|x-y|}\right)|x-y|^{\alpha-d}\, .
\end{equation}
\end{prop}
\pf  Let $r>0$ and $Y^{(r)}_t:=rY_{r^{-\alpha} t}$. 
Let $(\EE^{(r)}, \DD(\EE^{(r)}))$ be the Dirichlet form of $Y^{(r)}$. 
It was shown in the proof of \cite[Lemma 5.1]{KSV} that,
for $f, g\in C^\infty_c(\R^d_+)$,
it holds that 
$\EE^{(r)}(f,g)=\EE(f,g)$. Since $\EE(Gf,g)=\int_{\R^d_+}f(x)g(x)\, dx$, 
we see that $Gf$ is the $0$-potential operator of $Y^{(r)}$. In particular, 
$G^{(r)}(x,y) := G(x,y)$ 
is the Green function of $Y^{(r)}$.

Let $(Q_t)$ be the semigroup of $Y$ and $(Q^{(r)}_t)$ the semigroup of $Y^{(r)}$. 
For $f:\R^d_+\to [0,\infty)$ define $f^{(r)}(x)=f(rx)$. 
Then $Q^{(r)}_t f(x)=Q_{r^{-\alpha} t}f^{(r)}(x/r)$, 
implying that
$$
G^{(r)}f(x)=\int_0^{\infty} Q^{(r)}_t  f(x)\, dt=\int_0^{\infty}
Q_{r^{-\alpha}t}f^{(r)}(x/r)\, dt=r^{\alpha}\int_0^{\infty}Q_s f^{(r)}(x/r)\, ds=
r^{\alpha}Gf^{(r)}(x/r)\, .
$$
Then
\begin{align*}
&\int_{\R^d_+} G(x,y)f(y)\,dy =Gf(x)= r^{\alpha}Gf^{(r)}(x/r)=r^{\alpha}\int_{\R^d_+} G(x/r, y)f^{(r)}(y)\, dy\\
&=r^{\alpha-d}\int_{\R^d_+}G(x/r, z/r)f^{(r)}(z/r)\, dz
= r^{\alpha-d}\int_{\R^d_+}G(x/r, y/r)f(y)\, dy\, .
\end{align*}
This implies that for every $x\in \R^d_+$, $G(x,y)=r^{\alpha-d} G(x/r, y/r)$ for a.e.~$y$. 

Note that since $(Y_t, \P_x)\stackrel{d}{=}(Y^{(r)}, \P_{x/r})$, the processes $Y$ and $Y^{(r)}$ have same excessive functions. Thus, if $f$ is excessive for $Y$, it is also excessive for $Y^{(r)}$ 
and therefore $Q_{r^{-\alpha}t}f^{(r)}f(x/r)=Q^{(r)}_t f(x)\uparrow f(x)$ as $t\to 0$. 
Thus we also have $Q_{t}f^{(r)}f(y)\uparrow f(ry)=f^{(r)}(y)$ as $t\to 0$, 
proving that $f^{(r)}$ is also excessive for $Y$.  In particular, for every $x\in \R^d$, $y\mapsto r^{\alpha-d}G(x/r, y/r)$ is excessive for $Y$. Since this function is for a.e.~$y$ equal to the excessive function $y\mapsto G(x,y)$, it follows that they are equal everywhere. Thus for all $x,y\in \R^d_+$,
$$
G(x,y)=r^{\alpha-d}G(x/r,y/r).
$$
By taking $r=|x-y|$ we obtain \eqref{e:scaling-of-G}. \qed


\section{Decay of the Green function}\label{s:decayofGf}
The goal of this section is to show that the Green function $G(x,y)$ vanishes at the boundary of $\R^d_+$.
Recall that for $a,b>0$ and $\wt{w} \in \R^{d-1}$, 
$$
D_{\wt{w}}(a,b)=\{x=(\wt{x}, x_d)\in \R^d:\, |\wt{x}-\wt{w}|<a, 0<x_d<b\}.
$$
Due to \textbf{(A4)}, without loss of generality, 
we mainly deal with the case $\wt w=\wt{0}$. 
We will write $D(a,b)$ for $D_{\wt{0}}(a,b)$ and,
for $r>0$, $U(r)=D_{\wt{0}}(\frac{r}2, \frac{r}2).$
Further  we write  $U$ for $U(1)$.

In several places below
we will need the following upper bound for $\sB(x,y)$ proved in \cite[Lemma 5.2(a)]{KSV}: There exists a constant $C>0$ such that for all $x,y\in \R^d_+$ satisfying $|x-y|\ge x_d$, it holds  that
\begin{equation}\label{e:key-estimate-for-B}
\sB(x,y)\le C x_d^{\beta_1}(|\log x_d|^{\beta_3}\vee 1)\big(1+{\bf 1}_{|y|\ge1}(\log|y|)^{\beta_3}\big)|x-y|^{-\beta_1}.
\end{equation}

We now recall three key lemmas from \cite{KSV}. 
Recall that  $Y_t=(Y_t^1, \dots, Y_t^d)$.

\begin{lemma}
[{\cite[Lemma 5.7]{KSV}}] \label{l:upper-bound-for-integral}
For all $x\in U$, 
$$
\E_x \int_0^{\tau_{ U}} (Y_t^d)^{\beta_1}  |\log Y_t^d|^{\beta_3}\, dt \le  
x_d^p\, .
$$
\end{lemma}

In the next two lemmas, we have used the scaling property of $Y$.

\begin{lemma}
[{\cite[Lemma 5.10]{KSV}}]  \label{l:exit-probability-estimate}  
There exists  $C_7\in (0, 1)$ such that 
for all $r>0$ 
and all $x=(\wt{0}, x_d)\in D(r/8,r/8)$,
$$
\P_x(Y_{\tau_{D(r/4,r/4)}}\in D(r/4,r)\setminus D(r/4,3r/4)) \ge 
C_7\left(\frac{x_d}{r}\right)^p\, .
$$
\end{lemma}

\begin{lemma}
 [{\cite[Lemma 6.2]{KSV}}]\label{l:POTAl7.4} 
There exists $C_8>0$  such that 
for all $r>0$ and all $x \in D(2^{-5}r, 2^{-5}r)$,
$$
\P_x\left(Y_{\tau_{U(r)}}\in D(r, r)\right)\le C_8
 \left(\frac{x_d}{r}\right)^p.
$$
\end{lemma}

The L\'evy system formula (see \cite[Theorem 5.3.1]{FOT} and the arguments in \cite[p.40]{CK03}) states 
that for any 
non-negative Borel function $F$ on $\R^d_+\times \R^d_+$ vanishing on the diagonal and any stopping time $T$, it holds that
\begin{align}
\label{e:levys}
\E_x\sum_{s\le T} F (Y_{s-}, Y_s)=
\E_x\left(\int^T_0\int_{\R^d_+} F (Y_s, y)J(Y_s, y) dyds\right), \quad x\in \R^d_+.
\end{align} 
Here $Y_{s-}=\lim_{t\uparrow s}Y_t$ denotes the left limit of the process $Y$ at time $s>0$. We will use \eqref{e:levys} in the following form: Let $f:\R^d_+\to [0,\infty)$ be a Borel function, and  let 
 $V,W$ be two  Borel subsets of $\R^d_+$
with disjoint closures. 
If $F(x,y):=\1_V(x)\1_W(y)f(y)$, and $T=\tau_V$, then \eqref{e:levys} reads
\begin{align}
& \E_x\left[f(Y_{\tau_V}), Y_{\tau_V}\in W\right]=\E_x \sum_{s\le \tau_V} \1_V(Y_{s-})\1_W(Y_s)f(Y_s)\nn \\
&=\E_x \int_0^{\tau_V}\int_{\R^d_+}\1_V(Y_s)\1_W(y)f(y)J(Y_s,y)dy\, ds 
= \E_x \int_0^{\tau_V} \int_W f(y)J(Y_s,y)dy\, ds \label{e:levys2}\\
&= \int_V G^V(x,z) \int_W f(y)J(z,y)dy\, dz. \label{e:levys3}
\end{align}
The last line follows from the formula for the Green potential already described in Section \ref{s:existence-GF}.

\medskip
The following lemma is an improvement of Lemma \ref{l:POTAl7.4}, since $\R^d_+$ is a larger set than any $D(r,r)$.

\begin{lemma}\label{l:new-lemma}
There exists $C_9>0$ 
such that for all $r>0$ and  $x\in D(2^{-5}r, 2^{-5}r)$ we have that
\begin{equation}\label{e:new-lemma}
\P_x(Y_{\tau_{U(r)}}\in \R^d_+)\le C_9\left(\frac{x_d}{r}\right)^p.
\end{equation}
\end{lemma}
\pf By scaling, it suffices to prove \eqref{e:new-lemma} for $r=1$. Let $U=U(1)$ and $D=D(1,1)$. By Lemma \ref{l:POTAl7.4} we only need to show that $\P_x(Y_{\tau_U}\in \R^d_+\setminus D)\le c_1x_d^p$ for some $c_1>0$. 
By using \eqref{e:levys2} (with $f\equiv 1$) in the first line and
\eqref{e:key-estimate-for-B} in the second, 
\begin{align*}
&\P_x(Y_{\tau_U}\in \R^d_+\setminus D) = \E_x \int_0^{\tau_U} 
\int_{\R^d_+\setminus D}  J  (w, Y_t)\, dw\, dt\\
&\le  c_2 \E_x \int_0^{\tau_U} (Y_t^{ d})^{\beta_1}|\log Y_t^d|^{\beta_3}\, dt \int_{\R^d_+\setminus D}\frac{1+\1_{|w|>1}(\log|w|)^{\beta_3}}{|w|^{d+\alpha+\beta_1}}
\, dw.
\end{align*}
Since
$$
\int_{\R^d_+\setminus D}\frac{1+\1_{|w|>1}(\log|w|)^{\beta_3}}{|w|^{d+\alpha+\beta_1}}\, dw<\infty,
$$
it follows from Lemma \ref{l:upper-bound-for-integral} 
that $\P_x(Y_{\tau_U}\in \R^d_+\setminus D)\le c_3x_d^p$. \qed

The next result allows us to apply Theorem \ref{t:carleson} to get the Proposition \ref{p:gfcnub}, which is a key for us to get sharp two-sided Green function estimates.

\begin{thm}\label{t:green-function-decay}
For each $y\in \R^d_+$, it holds that $\lim_{x_d\to 0}G(x,y)=0$.
\end{thm}
\pf By translation invariance it suffices to show that $\lim_{|x|\to 0}G(x,y)=0$.  We fix $y\in \R^d_+$ and consider $x\in \R^d_+$ with $|x|<2^{-10}y_d$. Let $B_1=B(y, y_d/2)$ and $B_2=B(y,y_d/4)$. For $z\in B_1$ we have 
$z_d\ge y_d/2$ so that $|z-y|\le y_d/2\le z_d$. 
Moreover, $|z-x|\ge y_d/2-x_d\ge (7/16)y_d$. Thus, by the regular harmonicity of $G(\cdot, y)$ 
(see Proposition \ref{p:existenceGF}), 
\begin{equation}\label{e:green-function-decay}
G(x,y)=\E_x[G(Y_{T_{B_1}}, y), Y_{T_{B_1}}\in B_1\setminus B_2]+\E_x[G(Y_{T_{B_1}}, y), Y_{T_{B_1}}\in B_2]=:
I_1+I_2,
\end{equation}
where, for any $V\subset \R^d_+$, $T_V:=\inf\{t>0: Y_t\in V\}$.
By the Harnack inequality and Lemma \ref{l:transience}, 
\begin{align*}
\sup_{z\in B_1\setminus B_2}G(z,y) \le &
\frac{c_1}{|B_1\setminus B_2|} \int_{B_1\setminus B_2} G(z,y)dz \le 
c_2\frac{ y_d^\alpha }{y_d^d}\int_{B_1\setminus B_2} G(y,z) \kappa(z) dz\\
 \le &c_2 y_d^{\alpha-d} G\kappa (y) \le c_2 y_d^{\alpha-d}.
\end{align*}
In the second inequality we used the definition of $\kappa$ in \eqref{e:kappap}, that $z_d\asymp y_d$ in $B_1\setminus B_2$, and the fact that $|B_1\setminus B_2|\asymp y_d^d$.
Now we have
$$
I_1\le \sup_{z\in B_1\setminus B_2}G(z,y)\P_x(Y_{T_{B_1}}\in B_1\setminus B_2)\le 
\frac{c_2}{y^{d-\alpha}}\P_x(Y_{T_{B_1}}\in B_1\setminus B_2).
$$
Further, it is easy to check that 
$J(w,z)\asymp J(w,y)$ for all $w\in \R^d_+\setminus B_1$ and $z\in B_2$. 
Moreover, by 
Lemma \ref{l:transience},
$$
\int_{B_2}G(y,z)\, dz
\le c_3 y_d^\alpha \int_{B_2}G(y,z)\kappa(z) dz \le  c_3 y_d^\alpha G\kappa (y) \le  c_3 y_d^\alpha.
$$
Therefore, 
by \eqref{e:levys2} 
(with $f=G(\cdot, y)$) in the first line,
\begin{eqnarray*}
I_2 &= & \E_x \int_0^{T_{B_1}} \int_{B_2} J (Y_t,z)G(z,y)\, dz\, dt \\
&\le &c_4  \E_x \int_0^{T_{B_1}} J(Y_t, y)y_d^{\alpha} \, dt 
\le  c_5 y_d^{\alpha}\, \E_x \int_0^{T_{B_1}} \left(\frac{1}{|B_2|}\int_{B_2}
 J (Y_t,z)\, dz \right)dt\\
&=& \frac{c_6}{y_d^{d-\alpha}}\P_x(Y_{T_{B_1}}\in B_2).
\end{eqnarray*}
Inserting the estimates for $I_1$ and $I_2$ into \eqref{e:green-function-decay} and using Lemma \ref{l:new-lemma} we get that
$$
G(x,y)\le \frac{c_7}{y_d^{d-\alpha}}\P_x(Y_{T_{B_1}}\in \R^d_+)\le \frac{c_7}{y_d^{d-\alpha}}\P_x(Y_{\tau_{U(y_d/4)}}\in \R^d_+)\le \frac{c_8}{y_d^{d-\alpha-p}}x_d^p,
$$
which implies the claim. \qed


\section{Interior estimate of Green functions}\label{s:interior}

\subsection{Lower bound }

We first use a capacity argument to show
that there exists $c>0$ such that $G(x,y)\ge c$ 
for all $x,y\in \R^d_+$ satisfying $|x-y|=1$ 
and $x_d\wedge y_d\ge 10$.
 For such $x$ and $y$, let $U=B(x,5)$, $V=B(x,3)$ and $W_y=B(y,1/2)$. 
Recall that, for any $W\subset \R^d_+$, $T_W=\inf\{t>0: Y_t\in W\}$.
By the 
Krylov-Safonov type estimate \cite[Lemma 3.12]{KSV} , there exists a constant $c_1>0$ such that
\begin{align}
\label{e:TW}
\P_x(T_{W_y}<\tau_U)\ge c_1\frac{|W_y|}{|U|}=c_2 >0\, .
\end{align}
Recall that $Y^U$ is the process $Y$ killed upon exiting $U$ and  
$G^U(\cdot, \cdot)$ is  the Green function of $Y^U$. 
The Dirichlet form of $Y^U$  is $(\EE, \FF_U)$, where
$$
\EE(u,v)=\frac{1}{2}\int_U \int_U (u(x)-u(y))(v(x)-v(y))J(x,y)\, dy \, dx+\int_U u(x)^2\kappa_U(x)\, dx,
$$
\begin{equation}\label{e:kappaU}
\kappa_U(x)=\int_{{\R^d_+}\setminus U} J(x,y)\, dy +\kappa(x)\, , \quad x\in U\, ,
\end{equation}
and $\FF_U=\{u\in \FF:\, u=0 \textrm{ q.e.~on  } \R^d_+\setminus U\}$.
Here q.e.~means that the equality holds quasi-everywhere, that is, except on a set of capacity zero with respect to $Y$. 
Let $\mu$ be the capacitary measure  of $W_y$ with respect to $Y^U$ (i.e., with respect to the corresponding Dirichlet form). Then $\mu$ is concentrated on $\overline{W_y}$, $\mu(U)=\mathrm{Cap}^{Y^U}(W_y)$ and $\P_x(T_{W_y}<\tau_U)=G^U\mu(x)$.  By
\eqref{e:TW} and applying
Theorem  \ref{t:uhp} (Harnack inequality) to the function $G(x,\cdot)$,  we get
\begin{align}
&c_2\le \P_x(T_{W_y}<\tau_U)=G^U\mu(x)=\int_U G^U(x,z)\mu (dz)\le \int_U G(x,z)\mu (dz)
\nonumber \\ 
&\le c_3  G(x,y) \mu(U)=
c_3 G(x,y)\mathrm{Cap}^{Y^U}(W_y)\, .\label{e:1}
\end{align}

Let $X$ be the isotropic $\alpha$-stable process in $\R^d$ with the jump kernel $j(x,y)=|x-y|^{-d-\alpha}$. For $u,v:\R^d\to \R$, let
\begin{eqnarray*}
\QQ(u,v)&:=&\frac{1}{2}\int_{\R^d} \int_{\R^d} (u(x)-u(y))(v(x)-v(y))j(|x-y|)\, dy \, dx\, ,\\
\DD(\QQ)&:=& \{u\in L^2(\R^d, dx): \QQ(u,u)<\infty\}.
\end{eqnarray*}
Then $(\QQ, \DD(\QQ))$ is the regular Dirichlet form corresponding to $X$. Let $X^U$ 
denote the part of the process $X$ in $U$. 
The Dirichlet form of $X^U$ is $(\QQ, \DD_U(\QQ))$, where
$$
\QQ^{U}(u,v)=\frac{1}{2}\int_U \int_U (u(x)-u(y))(v(x)-v(y))j(|x-y|)\, dy \, dx+\int_U u(x)^2\kappa^X_U(x)\, dx, 
$$
$$
\kappa^{X}_U(x)=\int_{\R^d\setminus U}j(|x-y|)\, dy\, ,\quad x\in U\, ,
$$
and $\DD_U(\QQ)=\{u\in \DD(\QQ) :\, u=0 \textrm{ q.e.~on  } \R^d\setminus U\}$. 
Using calculations similar to those in \cite[p.13]{KSV}, 
one can show that  $\kappa_U(x)\asymp \kappa_U^X(x)$ for $x\in U$. 
Thus, there exists $c_4>0$ such that $\EE(u,u)\le c_4 \QQ^{U}(u,u)$ for all 
$u\in C_c^{\infty}(U)$ which is a core for both 
$(\QQ, \DD_U(\QQ))$ and $(\EE, \FF_U)$. 
This implies that 
$$
\mathrm{Cap}^{Y^U}(W_y)\le c_4 \mathrm{Cap}^{X^U}(W_y)\le  c_4 \mathrm{Cap}^{X^U}(V)\, .
$$
The last term, $\mathrm{Cap}^{X^U}(V)$, the capacity of $V$ with respect to $X^U$, is just a number, say $c_5$,  depending only on the radii of $V$ and $U$. Hence, $\mathrm{Cap}^{Y^U}(W_y)\le c_4 c_5$. 
Inserting in \eqref{e:1}, we get that
$$
G(x,y)\ge c_2c_3^{-1}c_4^{-1}c_5^{-1}.
$$
Combining this with 
the Harnack inequality (Theorem  \ref{t:uhp})
and \eqref{e:scaling-of-G}, we immediately get the following

\begin{prop}\label{p:green-lower-bound}
For any $C_{10}>0$, there exists a constant $C_{11}>0$ 
such that for all $x,y\in \R^d_+$ satisfying 
$|x-y|\le C_{10}(x_d\wedge y_d)$, it holds that
$$
G(x,y)\ge C_{11}|x-y|^{-d+\alpha}.
$$
\end{prop}
\pf 
We have shown above that there is $c_1>0$ such that 
$G(z,w)\ge c_1$ for all $z,w\in \R^d_+$ with $|z-w|=1$ and $z_d\wedge w_d\ge 10$.
By 
the Harnack inequality (Theorem  \ref{t:uhp}), 
there exists $c_2>0$ such that 
$G(z,w)\ge c_2$ for all $z,w\in \R^d_+$ with $|z-w|=1$ and $z_d\wedge w_d>C^{-1}_{10}$.

Now let $x,y\in \R^d_+$ satisfy $|x-y|\le  C_{10}(x_d\wedge y_d)$ and set
$$
x^{(0)}=\frac{x}{|x-y|}, \quad y^{(0)}=\frac{y}{|x-y|}.
$$
Then $|x^{(0)}-y^{(0)}|=1$ and 
$x^{(0)}_d\wedge y^{(0)}_d> C^{-1}_{10}$ 
so that $G(x^{(0)},y^{(0)})\ge c_2$. 
By scaling (Proposition \ref{p:scaling-of-G}),
$$
G(x,y)=G(x^{(0)},y^{(0)})|x-y|^{\alpha-d}\ge \frac{c_2}{|x-y|^{d-\alpha}}. 
$$
\qed

\medskip
As a corollary of the lower bound above we get that for every $x\in \R^d_+$,
$$
\lim_{y\to x} G(x,y)=+\infty .
$$

\subsection{Upper bound}

The purpose of this subsection is to establish the interior upper bound on 
the Green function $G$, Proposition 
\ref{p:green-upper-bound}. By 
\eqref{e:scaling-of-G} and 
the Harnack inequality (Theorem  \ref{t:uhp}),
 it suffices to deal with 
$x, y\in \R^d_+$ with $|x-y|=1$ and 
$x_d=y_d>10$.

We fix now two points $x^{(0)}$ and $y^{(0)}$ in $\R^d_+$ such that $|x^{(0)}-y^{(0)}|=1$,
$x^{(0)}_d=y^{(0)}_d>10$ and $\wt{x^{(0)}}=\wt{0}$. Let $E=B(x^{(0)},1/4)$,
$F=B(y^{(0)},1/4)$ and $D=B(x^{(0)},4)$. Let 
$f=G\1_E$ and $u=G\1_D$. 
Since $z \mapsto G(y^{(0)}, z )$ is harmonic in $B(x^{(0)},1/2)$ with respect to $Y$
and $f$ is harmonic in  $B(y^{(0)},1/2)$ with respect to $Y$, 
by applying 
the Harnack inequality (Theorem  \ref{t:uhp}) 
to $f$ and $z \mapsto G(y^{(0)}, z )$, we get 
$$
f(y^{(0)}) =\int_E G(y^{(0)}, z ) dz \ge c |E|G( y^{(0)}, x^{(0)})
\quad \text{and}\quad 
  \int_F  f(y)^2dy \ge c  |F| f(y^{(0)})^2.
$$
Thus, using the symmetry of $G$, we obtain
\begin{equation}\label{e:G-upper-bound}
G(x^{(0)}, y^{(0)})\le \frac{c}{|E|} f(y^{(0)})\le \frac{c}{|E|}\left(\frac{c}{|F|} \int_F  f(y)^2dy\right)^{1/2}  \le \frac{c^{3/2}}{|E|^{3/2}}\|u\|_{L^2(D)},
\end{equation}
for some constant $c>0$.
The key is to get uniform estimate on the $L^2$ norm of $u=G\1_D$,  see Proposition \ref{p:l2norm}.
To get the desired uniform estimate, we will use the Hardy inequality in 
\cite[Corollary 3]{BDK} 
and the heat kernel estimates of symmetric
jump processes with 
large jumps of lower intensity in \cite{BKKL}.

By \textbf{(A3)}, we have
\begin{equation}\label{e:B-lower-bound}
\sB(x,y)\ge 
c_1\left\{\begin{array}{lc}
{|x-y|^{-\beta_1-\beta_2}}
 & \textrm{if }|x-y|\ge 1 \textrm{ and }x_d\wedge y_d \ge 1,\\
1 &  \textrm{if }|x-y|< 1 \textrm{ and }x_d\wedge y_d \ge 1. 
\end{array}\right. 
\end{equation}
 Define
$$
\phi(r):=r^{\alpha}\1_{\{r<1\}}+r^{\alpha+\beta_1+\beta_2}\1_{\{r\ge 1\}} \quad \text{and }\quad \Phi(r):=\frac{r^2}{\int_0^r \frac{s}{\phi(s)}ds}.
$$
Let $\overline{\beta}:=(\alpha+\beta_1+\beta_2)\wedge 2$. Then 
$$\Phi(r)
\asymp
\begin{cases}
r^{\alpha}& \text{if }r\le 1,\\
r^{\overline{\beta}}& \text{if }r>1\text{  and  }\alpha+\beta_1+\beta_2\neq 2,\\
r^2/\log(1+r)& \text{if }r>1\text{ and }\alpha+\beta_1+\beta_2= 2,
\end{cases}
$$
which implies that
\begin{align}\label{e:Phnew}
c_2\left( \frac{R}{r} \right)^{\alpha}  \le \frac{\Phi(R)}{\Phi(r)} \le c_3
\left( \frac{R}{r} \right)^{\overline{\beta}}, \quad 0<r\le R<\infty.
\end{align}

For $a>0$, 
let $\R^d_{a+} :=\{x\in \R^d_+:\, x_d \ge a\}$. 
Define
\begin{equation}
K(r):=\left\{\begin{array}{lc}
r^{-d-\alpha}, & \textrm{if } r\le 1, \\
r^{-d-\alpha-\beta_1-\beta_2}, & \textrm{if } r> 1,
\end{array}\right.
\end{equation}
and 
\begin{equation}
Q(u,u):=
\int_{\R^d_{1+} }\int_{ \R^d_{1+} } 
(u(x)-u(y))^2 K(|x-y|)\, dx\, dy.
\end{equation}
Note that, by  \eqref{e:B-lower-bound},
\begin{equation}\label{e:KJ}
K(|x-y|)\le c_4J(x, y) \le c_5 j(|x-y|), \quad (x,y) 
 \in \R^d_{1+}\times \R^d_{1+}
\end{equation}
for some positive constants  $c_4$ and $c_5$.
Consider the Dirichlet form $(Q, \mathcal D(Q))$ on 
$\R^d_{1+}$,
where
\begin{equation}\label{e:DQ}
\mathcal D(Q)=\{u\in 
L^2(\R^d_{1+}): 
Q(u, u)<\infty\}.
\end{equation}
Let
$$
\wt Q(u,u):=
\int_{\R^d_{1+} }\int_{ \R^d_{1+} } 
\frac{(u(x)-u(y))^2}{|x-y|^{d+\alpha}}\, dx\, dy
$$
and 
$$
\mathcal D(\wt Q)=\{u\in L^2(\R^d_{1+}): 
\wt Q(u, u)<\infty\}.
$$
It follows from \cite[Remark 2.1.(1)]{BBC} (more precisely the first sentence on \cite[p. 98]{BBC}) that $(\wt Q, \mathcal D(\wt Q))$ is a regular Dirichlet form.
Moreover, we have
\begin{align*}
&
\wt Q(u, u)=\int_{\R^d_{1+} \times \R^d_{1+} } 
\frac{(u(x)-u(y))^2}{ |x-y|^{d+\alpha}}\, dx\, 
dy\\
& = 
\int_{\R^d_{1+} \times \R^d_{1+} } 
{\bf 1}_{|x-y| \le 1} \frac{(u(x)-u(y))^2}{ |x-y|^{d+\alpha}}\, dx\, 
dy+
\int_{\R^d_{1+} \times \R^d_{1+} } 
{\bf 1}_{|x-y| > 1} \frac{(u(x)-u(y))^2}{ |x-y|^{d+\alpha}}\, dx\, 
dy\\
&\le 
Q(u, u)+
4\|u\|^2_{L^2(\R^d_{1+})}\sup_{y\in  \R^d_{1+}}\int_{ \R^d_{1+} }
{\bf 1}_{|x-y| > 1} |x-y|^{-d-\alpha}\, dx
\\
&\le 
Q(u, u)+
4\|u\|^2_{L^2(\R^d_{1+})} \int_{ \R^d }
{\bf 1}_{|z| > 1} |z|^{-d-\alpha}\, dz= 
Q(u, u)+
c_6\|u\|^2_{L^2(\R^d_{1+})}. 
\end{align*}
This implies that the Dirichlet form $(Q, \mathcal D(Q))$ is also regular on 
$L^2(\R^d_{1+}, dx)$.

 Let $X^{(1)}=(X^{(1)}_t)_{t\ge 0}$ be 
the symmetric Hunt process associated with $(Q, \mathcal D(Q))$
and denote by  $p^{(1)}(t,x,y)$   the transition density of  $X^{(1)}$. 
\cite[Theorem 4.6]{BKKL} says that there exists $c_7>0$ such that
\begin{align}
p^{(1)}(t,x,y) \le & c_7
\left(\frac{1}{\Phi^{-1}(t)^d}\wedge \frac{t}{|x-y|^d \Phi(|x-y|)}\right), 
\quad t>0, \ x,y\in  
\R^d_{1+}.
\label{e:ptxy-upper}
\end{align}
\cite[Theorem 2.19 (i)]{BKKL} says that there 
exists $c_8>0$ such that
\begin{align}
p^{(1)}(t,x,y) \ge &\frac{c_8}{\Phi^{-1}(t)^d}, 
 \quad t>0,  \ x,y\in 
 \R^d_{1+}
 \text{ with }|x-y|\le \Phi^{-1}(t).  \label{e:ptxy-lower}
\end{align}
Recall that we have assumed 
$d > \overline{\beta}$.
By using 
\eqref{e:Phnew}, \eqref{e:ptxy-upper} and \eqref{e:ptxy-lower},
 we can compute (see  \cite[p.241]{BDK}) that for every $\gamma\in (0, (d/\overline{\beta}-1)\wedge2)$,
$$
h(x,y):=\int_0^{\infty}t^{\gamma}  p^{(1)}(t,x,y)\, dt \asymp \frac{\Phi(|x-y|)^{\gamma+1}}{|x-y|^d}, \quad x,y\in 
\R^d_{1+},
$$
and
$$
\overline{h}(x,y):=\int_0^{\infty}t^{\gamma-1}  p^{(1)}(t,x,y)\, dt \asymp \frac{\Phi(|x-y|)^{\gamma}}{|x-y|^d},\quad x,y\in 
\R^d_{1+}.
$$
 This is the only place where the assumption 
$d > \overline{\beta}$
is used. 
Set $x^{\ast}=(\wt{0}, 1)$ and let
$$
q(x):=\frac{\overline{h}(x,x^{\ast})}{h(x,x^{\ast})}\asymp \frac{1}{\Phi(|x-x^{\ast}|)}.
$$
It follows from the Hardy inequality in
\cite[Theorem 2 and Corollary 3]{BDK} that there exists $c_9>0$ such that
\begin{equation}\label{e:estimate-Q}
Q(u,u)\ge c_9 
\int_{\R^d_{1+}}
u(x)^2 \frac{dx}{\Phi(|x-x^{\ast}|)} 
\quad \textrm{for all }u\in 
L^2(\R^d_{1+}).
\end{equation}
This estimate can be improved to obtain the following result.
\begin{prop}\label{p:estimate-Q}
There exists a constant $C_{12}>0$ such that for all 
$u\in \mathcal D(Q)$ and all $z_a=(\wt{0}, a)$ with $a\ge 0$, it holds that
$$
Q(u,u)\ge 
C_{12} 
\int_{\R^d_{1+}}
u(x+z_a)^2 \frac{dx}{\Phi(|x-x^{\ast}|)}.
$$
\end{prop}
\pf Let $z_a=(\wt{0}, a)$, $a\ge 0$. Then
\begin{eqnarray*}
&&
\int_{\R^d_{1+}}\int_{\R^d_{1+}} 
(u(x+z_a)-u(y+z_a))^2 K(|x-y|)\, dx \, dy\\
&=& 
\int_{\R^d_{(1+a)+}}\int_{\R^d_{(1+a)+}}
(u(x )-u(y))^2 K(|x-y|)\, dx \, dy\le  Q(u,u)<\infty.
\end{eqnarray*}
Thus, $u(\cdot+z_a)\in \mathcal D(Q)$ by \eqref{e:DQ} and 
\begin{align*}
Q(u(\cdot +z_a), u(\cdot +z_a))=
\int_{\R^d_{1+}}\int_{\R^d_{1+}} 
(u(x+z_a)-u(y+z_a))^2 K(|x-y|)\, dx \, dy\le  Q(u,u).
\end{align*}
Since clearly 
$u(\cdot+z_a)\in L^2(\R^d_{1+})$, 
the claim follows from \eqref{e:estimate-Q}. \qed

We have shown in Lemma \ref{l:transience} that  $(\EE, \FF)$  is transient. Let
$(\EE, \FF_e)$
be its extended Dirichlet space.

\begin{lemma}\label{t:extended-DF}
There exists $C_{13}>0$ such that for any 
$h\in  \FF_e$
and any $z_a=(\wt{0},a)$ with $a\ge 0$, it holds that
$$
\int_{\R^d_{1+}} 
\frac{|h(x+z_a)|^2}{\Phi(|x-x^{\ast}|)}\, dx \le C_{13}
\EE(h,h).
$$
\end{lemma}

\pf Let 
 $h\in \FF_e$.
There exists an approximating sequence $(g_n)_{n\ge 1}$ in 
$\FF$  such that 
$\EE(h,h)=\lim_{n\to \infty}\EE(g_n, g_n)$ and $h=\lim_{n\to \infty}g_n$ a.e. 
Since $g_n\in L^2(\R^d_+, dx)$, we have that 
$g_n\1_{\R^d_{1+}}\in L^2(\R^d_{1+}, dx)$. 
Further, by \eqref{e:KJ},
$$
Q(g_n \1_{\R^d_{1+}}, g_n \1_{\R^d_{1+}}) 
\le c_1 \EE(g_n,g_n) <\infty,$$ 
so that
$g_n \1_{\R^d_{1+}} \in  \mathcal D(Q)$ 
by \eqref{e:DQ}.

Now, using  Proposition \ref{p:estimate-Q} and the above inequality, we have that
$$
\EE(g_n,g_n)\ge c^{-1}_1
 Q(g_n \1_{\R^d_{1+}}, g_n \1_{\R^d_{1+}})
\ge c_2
\int_{\R^d_{1+}} 
g_n(x+z_a)^2\frac{dx}{\Phi(|x-x^{\ast}|)},
$$
for some constant $c_2>0$.
By Fatou's lemma,
\begin{align*}
\EE(h,h)&=\lim_{n\to \infty}
\EE(g_n, g_n)\ge 
c_2
\int_{\R^d_{1+}} 
\liminf_{n\to \infty}g_n(x+z_a)^2\frac{dx}{\Phi(|x-x^{\ast}|)}\\
&=c_2
\int_{\R^d_{1+}} 
h(x+z_a)^2\frac{dx}{\Phi(|x-x^{\ast}|)}. 
\end{align*}
\qed

By \cite[Theorem 1.5.4]{FOT}, for any non-negative Borel 
function $f$ satisfying $\int_{\R^d_+}f(x) Gf(x)\, dx <\infty$, we have that $Gf\in \FF_e$ and 
$\EE(Gf,Gf)=\int_{\R^d_+}f(x) Gf(x)\, dx$. Thus by Lemma \ref{t:extended-DF} we have
\begin{corollary}\label{c:extended-DF}
There exists $C_{14}>0$ such that for every non-negative Borel  
function $f$
satisfying $\int_{\R^d_+}f(x) Gf(x)\, dx <\infty$ and every $z_a=(\wt{0},a)$ 
with $a\ge 0$, it holds that
$$
\int_{\R^d_{1+}} 
\frac{|Gf(x+z_a)|^2}{\Phi(|x-x^{\ast}|)} \, dx  \le C_{14} \int_{\R^d_+}f(x) Gf(x)\, dx.
$$
\end{corollary}

\begin{prop}\label{p:l2norm} 
There exists $C_{15}>0$ such that for every $x^{(0)}\in \R^d_+$ with $x^{(0)}_d >6$,
$$
\int_{B(x^{(0)},4)}  (G1_{B(x^{(0)},4)}(x))^2\, dx \le C_{15}.
$$
\end{prop}
\pf Without loss of generality we assume that 
$x^{(0)}=(\wt{0}, x_d^{(0)})$. 
Set $B=B(x^{(0)},4)$ and let $u=G \1_B$. 
We first note that,  by \eqref{e:Green-potential-bounded} we have that $G 1_B \le c_{\overline{B}}^{-1}$, and therefore $\| u\|_{L^2(B)}<\infty$.

Let $z=(\wt{0}, x^{(0)}_d-6)$ and $\wt{B}=B((\wt{0},6),4)\subset 
\R^d_{2+}$.
 By using the change of variables $w=x-z$ and the fact that $\Phi(|w-x^{\ast}|)\asymp 1$ for $w\in\wt{B}$ in the first line, and then Corollary \ref{c:extended-DF} and the Cauchy inequality in the third line below, we have
\begin{eqnarray*}
 \| u\|_{L^2(B)}^2 &=& \int_{\wt{B}}|u(w+z)|^2\, dw \le c_1\int_{\wt{B}} |u(w+z)|^2 \frac{dw}{\Phi(|w-x^{\ast}|)}\\
&\le &c_1 
\int_{\R^d_{1+}} 
|u(w+z)|^2 \frac{dw}{\Phi(|w-x^{\ast}|)}=c_1 
\int_{\R^d_{1+}} 
|G\1_B(w+z)|^2 \frac{dw}{\Phi(|w-x^{\ast}|)}\\
&\le & c_2 \int_{\R^d_+}\1_B(x) G\1_B(x)\, dx \le c_2 |B|^{
1/2} \| u\|_{L^2(B)}.
\end{eqnarray*}
Since $\| u\|_{L^2(B)}<\infty$, we have that 
$ \| u\|_{L^2(B)} \le c_2 |B|^{
1/2}$. This completes the proof. \qed

Coming back to \eqref{e:G-upper-bound}, by Proposition \ref{p:l2norm}, we see that the right-hand side is bounded 
 above  by a constant, and therefore $G(x^{(0)}, y^{(0)})\le c$. 

\begin{prop}\label{p:green-upper-bound}
There exists a constant $C_{16}>0$ such that for all $x,y\in \R^d_+$ satisfying 
$|x-y|\le 8(x_d\wedge y_d)$, 
 it holds that
$$
G(x,y) \le C_{16} |x-y|^{-d+\alpha}.
$$
\end{prop}
\pf 
This is analogous to the proof of Proposition \ref{p:green-lower-bound}. We omit the details.
\qed

Using Theorem \ref{t:green-function-decay}, 
we can combine Proposition \ref{p:green-upper-bound}
with Theorem \ref{t:carleson} to get the following result, which is key for us to get sharp
two-sided Green functions estimates.

\begin{prop}\label{p:gfcnub}
There exists a constant $C_{17}>0$ such that for all $x, y\in \R^d_+$,
\begin{align}
\label{e:upper}
G(x,y) \le C_{17} |x-y|^{-d+\alpha}.
\end{align}
\end{prop}

\pf
By  Proposition \ref{p:green-upper-bound}, there exists 
$c_1>0$ such that $G(x, y)\le c_1$ 
for all $x, y\in \R^d_+$ with $|x-y|=1$ and $x_d\wedge y_d \ge1/8$.

Suppose that  $x, y\in \R^d_+$ with $|x-y|=1$ and $x_d \le y_d$ and $x_d<1/8<y_d$.
Since $z \to G(z, y)$ is harmonic 
in $B((\wt x, 0), 1/4)$ with respect to $Y$ and vanishes on the boundary of $\R^d_+$ by Theorem \ref{t:green-function-decay}, we can use 
Theorem \ref{t:carleson} and see that there exists 
$c_2>0$ such that
\begin{align}
\label{e:cc1}
G(x,y)  \le c_2G(x+(\wt 0, 1/8),y)  \le c_2 c_1.
\end{align}
Suppose that  $x, y\in \R^d_+$ with $|x-y|=1$ and $x_d \le y_d$ and $y_d \le 1/8$. Then, 
since $z \to G(z, y)$ is harmonic 
in $B((\wt x, 0), 1/4)$ with respect to $Y$ and vanishes on the boundary of $\R^d_+$, 
by \eqref{e:cc1} and Theorem \ref{t:carleson}, we see that
$G(x,y)  \le c_2G(x+(\wt 0, 1/8),y)  \le c_2^2c_1$.
Thus  for all $x, y\in \R^d_+$ with $|x-y|=1$, we have $G(x,y) \le C$.
Therefore, by \eqref{e:scaling-of-G}, we have
$$
G(x,y) \le C |x-y|^{-d+\alpha}, \quad x, y\in \R^d_+.
$$
\qed


\section{Preliminary Green Functions Estimates}\label{s:prelgfe}

The results of this section are valid for all $p\in ((\alpha-1)_+, \alpha+\beta_1)$.

\subsection{Lower bound}\label{ss:lb}
The goal of this subsection is to prove the following 
result, which is used later to prove the sharp lower bound of Green function for $p\in ((\alpha-1)_+, \alpha+\frac12[\beta_1+(\beta_1 \wedge \beta_2)])$.

\begin{thm}\label{t:GB}
Suppose $p\in ((\alpha-1)_+, \alpha+\beta_1)$.
For any $\eps \in (0, 1/4)$, 
there exists a constant $C_{18}>0$ such that for all 
$w \in \partial \R^d_+$,  $R>0$ and $x,y \in 
B(w, (1-\eps)R)\cap \R^d_+$, it holds that
$$
G^{B(w, R)\cap \R^d_+}(x, y)\ge C_{18}
\left(\frac{x_d}{|x-y|}  \wedge 1 \right)^p\left(\frac{y_d}{|x-y|}  \wedge 1 \right)^p \frac{1}{|x-y|^{d-\alpha}}.
$$
\end{thm}

The theorem will be proved through three lemmas.
For any $a>0$, let $B^+_a:=B(0, a)\cap \R^d_+$. 
Recall that 
$\R^d_{a+}=\{x\in \R^d_+: x_d \ge a\}$.

\begin{lemma}\label{l:GB_1}
For any $\eps \in (0, 1)$ and $M>1$, there exists a constant 
$C_{19}>0$ such that for all 
$y,z \in B^+_{1-\eps}$ with $|y-z|  \le M(y_d\wedge z_d)$, 
$$
G^{B^+_1}(y, z)\ge C_{19} |y-z|^{-d+\alpha}.
$$
\end{lemma}

\pf 
By using \eqref{e:gf-formula} in the first equality below, it 
follows from Propositions \ref{p:gfcnub} and \ref{p:green-lower-bound}
that there exists $c_1>1$ such that for all 
$y,z \in B^+_{1-\eps}$ with $|y-z|  \le M(y_d\wedge z_d)$,
$$
G^{B^+_1}(y, z)=G(y, z)-\E_y[G(Y_{\tau_{B^+_1}}, z)] 
\ge c_1^{-1}|y-z|^{-d+\alpha}-c_1 \eps^{-d+\alpha}.
$$
Now, we choose $\delta=(2c_1^2)^{-\frac1{d-\alpha}}$.  Then for all 
$y,z \in B^+_{1-\eps}$ with $|y-z|  \le (\delta \eps) \wedge M(y_d\wedge z_d)$, 
\begin{align}
\label{e:GB22}
&G^{B^+_1}(y, z) \ge c_1^{-1}|y-z|^{-d+\alpha}-c_1 (\delta^{-1}|y-z|)^{-d+\alpha}\nn \\
&\ge (c_1^{-1}-c_1 \delta^{d-\alpha})|y-z|^{-d+\alpha} 
=(2c_1)^{-1}|y-z|^{-d+\alpha}.
\end{align}

Assume that $y,z \in B^+_{1-\eps}$ with $|y-z|  \le M(y_d\wedge z_d)$ are such that also $|y-z|  \le \delta \eps$. Then clearly $|y-z|  \le (\delta \eps) \wedge M(y_d\wedge z_d)$, and \eqref{e:GB22} proves the lemma.

\begin{figure}[ht]
\begin{center}
\begin{tikzpicture}
\draw[thick,->] (0,0)--(14,0);
\draw (2,0) arc (180:0:5);
\draw[dashed] (0,0.8)--(14,0.8);
\draw (13.5,1.1) node {$\R^d_{\frac{\delta \epsilon}{M}+}$};
\draw (12,-0.5) node {$1-\epsilon$};
\draw (7,-0.5) node {$0$};
\draw (7, 0) node {$\bullet$};
\draw (5,3.5) circle (0.8);
\draw (5,3.5) node {$\cdot$};
\draw (5.14,3.11) node {$\cdot$};
\draw (4.85,3.5) node {\scriptsize $y$};
\draw (5.15,2.95) node {\scriptsize $w$};
\draw (5.97,3.60) node {\scriptsize $\frac{\delta\epsilon}{M}$};
\draw (5, 3.5)--(5.8, 3.6);
\draw (5, 3.5)--(5.14, 3.11);
\draw (5,3.5) circle (0.4); 
\draw (8,2.5) node {$\cdot$};
\draw (8.13,2.5) node {\scriptsize $z$};
\end{tikzpicture}
\end{center}
\caption{}
\label{fig:1}
\end{figure}
Now, we assume that $y,z \in B^+_{1-\eps}$ with $ |y-z| \le M(y_d\wedge z_d)$, but  $|y-z| > \delta \eps$,  see Figure \ref{fig:1}.
Since $ y_d\wedge z_d > \delta \eps/M$, 
we have
\begin{align}
\label{e:GB23}
y,z \in B^+_{1-\eps}\cap 
\R^d_{(\delta \eps/M)+}.
\end{align}
Therefore we 
can choose
a point $w\in  B(y, \delta \eps/M)$ such that $|y-w|=\delta \eps/(2M)$ 
and 
$w \in B^+_{1-\eps}\cap \R^d_{(\delta \eps/M)+}$. 
Since $M(y_d\wedge w_d)>\delta \epsilon >|y-w|$, we can use \eqref{e:GB22} for points $y$ and $w$ to conclude that 
$$
G^{B^+_1}(y, w)\ge (2c_1)^{-1}|y-w|^{-d+\alpha}=(2c_1)^{-1}(\delta\eps/(2M))^{-d+\alpha}=:c_2.
$$
Since  $G^{B^+_1}(y, \cdot)$ is harmonic in 
$B(w, \delta \eps/(4M)) \cup B(z, \delta \eps/(4M))$ by \eqref{e:GB23}, we can use Theorem \ref{t:uhp} (b) and the fact that $|y-z|<\delta \epsilon$ to get
$$
G^{B^+_1}(y, z) \ge c_3 G^{B^+_1}(y, w)  \ge c_4 \ge c_5 |y-z|^{-d+\alpha}.
$$
\qed

\begin{lemma}\label{p:GB_1}
Suppose $p\in ((\alpha-1)_+, \alpha+\beta_1)$.
For every $\eps \in (0, 1/4)$ and $M, N>1$, there exists a constant 
$C_{20}>0$ such that for 
all $x,z \in B^+_{1-\eps}$ with $x_d \le z_d$ 
satisfying $x_d/N \le |x-z|\le M z_d $, it holds that
$$
G^{B^+_1}(x, z)\ge C_{20}x^p_d|x-z|^{-d+\alpha-p}.
$$
\end{lemma}

\pf 
 Without loss of generality, we assume $M>4/\eps$. 
 If $|x-z|\le Mz_d$ and $|x-z|\ge 20M x_d$,
let $r=\frac{|x-z|}{10 M} \le \frac{1}{5 M} \le \frac{\eps}{20}$. 
Since $x\mapsto G^{B^+_1}(x, z)$ is regular harmonic in $D_{\wt x}(r, r)$, 
and $D_{\wt x}(r, 4r) \setminus D_{\wt x}(r, 3r) \subset B^+_{1-\eps/4}$, 
by Lemmas \ref{l:GB_1} and \ref{l:exit-probability-estimate}, 
we have
\begin{align*}
&G^{B^+_1}(x,z) \ge \E_x[G^{B^+_1}
(Y_{\tau_{D_{\wt x}(r,r)}}, z): Y_{\tau_{D_{\wt x}(r,r)}} \in D_{\wt x}(r, 4r) \setminus D_{\wt x}(r, 3r)]\\
&\ge c_1|x-z|^{-d+\alpha} \P_x(Y_{\tau_{D_{\wt x}(r,r)}} \in D_{\wt x}(r, 4r) \setminus D_{\wt x}(r, 3r))
\ge c_2 x^p_d|x-z|^{-d+\alpha-p}, 
\end{align*}
since, for $y\in D_{\wt x}(r, 4r) \setminus D_{\wt x}(r, 3r)$, $|y-z|\le |x-z|+|x-y|
\le 5(2M+1)r \le 2(2M+1)(y_d\wedge z_d)$.

 If $|x-z|\le Mz_d$ and $x_d/N<|x-z|< 20M x_d$, we simply use Lemma \ref{l:GB_1} (since $|x-z|< 12M (x_d \wedge z_d)$)
 and get 
 $$
G^{B^+_1}(x,z) \ge c_3|x-z|^{-d+\alpha}  \ge c_3 N^{-p}  x^p_d|x-z|^{-d+\alpha-p}.
 $$
\qed

\begin{lemma}\label{l:GB_4}
Suppose $p\in ((\alpha-1)_+, \alpha+\beta_1)$.
For every $\eps \in (0, 1/4)$ and $M \ge 40/\eps$, there exists a constant 
$C_{21}>0$ such that 
for all $x,z \in B^+_{1-\eps}$ with $x_d \le z_d$ 
satisfying $|x-z|\ge M z_d $, it holds that
$$
G^{B^+_1}(x, z)\ge C_{21}x^p_dz^p_d|x-z|^{-d+\alpha-2p}.
$$
\end{lemma}
\pf
Let $r=\frac{2|x-z|}{ M} \le \frac{4}{M} \le \frac{\eps}{10}$. 
Since $x\mapsto G^{B^+_1}(x, z)$ is regular harmonic in $D_{\wt x}(r, r)$,  
 and $D_{\wt x}(r, 4r) \setminus D_{\wt x}(r, 3r) \subset B^+_{1-\eps/4}$, 
by Lemmas  \ref{p:GB_1} and \ref{l:exit-probability-estimate}, 
we have
\begin{align*}
&G^{B^+_1}(x,z) \ge \E_x[G^{B^+_1}
(Y_{\tau_{D_{\wt x}(r,r)}}, z): Y_{\tau_{D_{\wt x}(r,r)}} \in D_{\wt x}(r, 4r) \setminus D_{\wt x}(r, 3r)]\\
&\ge c_1z_d^p|x-z|^{-d+\alpha-p} \P_x(Y_{\tau_{D_{\wt x}(r,r)}} \in D_{\wt x}(r, 4r) \setminus D_{\wt x}(r, 3r))
\ge c_2 x^p_d z_d^p|x-z|^{-d+\alpha-2p}
\end{align*}
since, for $y\in D_{\wt x}(r, 4r) \setminus D_{\wt x}(r, 3r)$, 
$|y-z|\le |x-z|+|x-y|\le (M/2+5)r\le (M/2+5)y_d$ and $|y-z|\ge |x-z|-|x-y|\ge 75r\ge 150z_d$.
 \qed

Combining the above result with scaling, we get the result of Theorem \ref{t:GB}.

\subsection{Upper bound}
The goal of this subsection is to get the following preliminary upper bound
on the Green function.

\begin{lemma}\label{l:prelub} 
Suppose $p\in ((\alpha-1)_+, \alpha+\beta_1)$.
There exists $C_{22}>0$ such that 
\begin{align}\label{e:Gu}
G(x, y) \le C_{22}
\left(\frac{x_d\wedge y_d}{|x-y|} \wedge 1\right)^p
\frac{1}{|x-y|^{d-\alpha}}, \quad x,y \in \R^d_+.
\end{align}
\end{lemma}

\pf 
Suppose  $x, y\in \R^d_+$ satisfy  
$\wt x=\wt 0$, 
$x_d \le 2^{-9}$ and $|x-y|=1$.
Let $r=2^{-8}$.
For $z\in U(r)$ and $w\in \R^d_+\setminus D(r, r)$, we have $|w-z|\asymp |w|$. Thus, by using 
\eqref{e:key-estimate-for-B} 
and  Proposition  \ref{p:gfcnub},
\begin{align}
&\int_{\R^d_+\setminus D(r, r)}G(w,y) \sB(z,w)|z-w|^{-d-\alpha} dw\nonumber\\
&\le c_1z_d^{\beta_1}(|\log z_d|^{\beta_3}\vee 1)\int_{\R^d_+\setminus D(r, r)}\frac{G(w,y)}{|w|^{d+\alpha+\beta_1}}
\big(1+{\bf 1}_{|w|\ge1}(\log|w|)^{\beta_3}\big)
 dw\label{e:GG}\\
 &\le c_2 z_d^{\beta_1} |\log z_d|^{\beta_3}\int_{\R^d_+\setminus D(r, r)}\frac{ \big(1+{\bf 1}_{|w|\ge1}(\log|w|)^{\beta_3}\big)}{|w-y|^{d-\alpha}|w|^{d+\alpha+\beta_1}}
 dw.\nonumber
\end{align}
Hence, by using \eqref{e:levys3} and \eqref{e:key-estimate-for-B} in the second line,  and Lemma \ref{l:upper-bound-for-integral}  in the third,
\begin{align*}
&\E_x\left[G(Y_{\tau_{U(r)}}, y); Y_{\tau_{U(r)}}\notin D(r, r)\right]\\
&\le c_3 \E_x\int^{\tau_{U(r)}}_0(Y^d_t)^{\beta_1}  |\log (Y^d_t)|^{\beta_3}dt\int_{\R^d_+\setminus D(r, r)}\frac{ \big(1+{\bf 1}_{|w|\ge1}(\log|w|)^{\beta_3}\big)}{|w-y|^{d-\alpha}|w|^{d+\alpha+\beta_1}}
dw\\
&\le c_4 x_d^p \int_{\R^d_+\setminus D(r, r)}\frac{ \big(1+{\bf 1}_{|w|\ge1}(\log|w|)^{\beta_3}\big)}{|w-y|^{d-\alpha}|w|^{d+\alpha+\beta_1}}
dw.
\end{align*}
Let 
\begin{align}\label{e:TAMSe}
 \int_{\R^d_+\setminus D(r, r)}\frac{ \big(1+{\bf 1}_{|w|\ge1}(\log|w|)^{\beta_3}\big)}{|w-y|^{d-\alpha}|w|^{d+\alpha+\beta_1}}
 dw =
  \int_{\R^d_+ \cap B(y, r)} 
  +
    \int_{\R^d_+\setminus (D(r, r) \cup B(y, r))}
    =: I+II.
\end{align}
It is easy to see
\begin{align}\label{e:TAMSe2}
  II \le  r^{-d+\alpha} \int_{\R^d_+\setminus (D(r, r) \cup B(y, r))}\frac{ \big(1+{\bf 1}_{|w|\ge1}(\log|w|)^{\beta_3}\big)}{|w|^{d+\alpha+\beta_1}} dw <\infty
\end{align}
and
\begin{align}\label{e:TAMS}
  I \le   c_5  \int_{\R^d_+ \cap B(y, r)}\frac{1}{|w-y|^{d-\alpha}}dw  <\infty.
\end{align}
Thus,
\begin{align}\label{e:TAMSe6.41}
\E_x\left[G(Y_{\tau_{U(r)}}, y); Y_{\tau_{U(r)}}\notin D(r, r)\right] \le 
c_6 x_d^p.
\end{align}

Let $x_0:=(\wt 0, r)$.
By Theorem \ref{t:carleson}, Proposition  \ref{p:gfcnub} and Lemma \ref{l:POTAl7.4}, we have
\begin{align}\label{e:TAMSe6.42}
&\E_x\left[G(Y_{\tau_{U(r)}}, y); Y_{\tau_{U(r)}}\in D(r, r)\right]\le 
c_7 G(x_0, y)\P_x(Y_{\tau_{U(r)}}\in D(r, r))\le c_8 x^p_d.
\end{align}
 Combining \eqref{e:TAMSe6.41} and \eqref{e:TAMSe6.42}, we get that 
 for $x, y\in \R^d_+$ satisfying $x_d \le 2^{-9}$ and $|x-y|=1$,
\begin{align*}
G(x, y)=\E_x\left[G(Y_{\tau_{U(r)}}, y); Y_{\tau_{U(r)}}\notin D(r, r)\right]+ \E_x\left[G(Y_{\tau_{U(r)}}, y); Y_{\tau_{U(r)}}\in D(r, r)\right]\le 
c_9x^p_d.
\end{align*}
Combining this with Proposition  \ref{p:gfcnub}, \eqref{e:scaling-of-G} and symmetry, we immediately get
the desired conclusion. \qed

\section{Proof of Theorem \ref{t:Green}}\label{s:thm1.1}

We begin this section  by introducing an auxiliary function that will be needed later. 
For $\gamma\in \R$ and
 $\beta\ge 0$, 
we define a function on $(0, 1]$ by
$$
F(x; \gamma, \beta)=\int^1_xh^\gamma\left(\log \frac2{h}\right)^\beta dh.
$$
Note that $F(\cdot, \gamma, \beta)$ is a decreasing function on $(0, 1]$
and that, when $\gamma>-1$, $F(0+, \gamma, \beta)$ is finite. It is obvious
that
$$
F(x; \gamma, 0)=\begin{cases}\frac1{\gamma+1}(1-x^{\gamma+1}), &\gamma\neq-1,\\
-\log x, & \gamma=-1
\end{cases}
$$
and
\begin{align}
\label{e:F1}
F(x; -1, \beta)=\frac1{1+\beta}\left(\left(\log \frac 2{x}\right)^{1+\beta}-\left(\log 2\right)^{1+\beta}\right).
\end{align}
Note also that for any $b\in (0, 1)$,  on $(0, b]$, when $\gamma>-1$,
\begin{align}
\label{e:F2}
F(0; \gamma, \beta)-F(x; \gamma, \beta)\asymp x^{\gamma+1}\left(\log \frac2{x}\right)^\beta
\end{align}
and when $\gamma<-1$,
\begin{align}
\label{e:F3}
F(x; \gamma, \beta)\asymp x^{\gamma+1}\left(\log \frac2{x}\right)^\beta,
\end{align}
with comparison constants depending on 
$\beta \ge 0$ and $\gamma<-1$. 

We first present
a technical lemma 
 inspired by \cite[Lemma 3.3]{AGV}. This lemma
will be used several times in this section. 
For $x=(\wt{0}, x_d)\in \R^d_+$ and $ \gamma, q, \delta\in \R$, 
$R>0$, $\beta\ge 0$  and $y\in \R^d_+$ with $y_d\in (0, R)$,
we define
$$
f(y; \gamma, \beta, q, \delta, x)
:=y_d^{\gamma}|x-y|^{-d+\alpha-q}  
\left(\log\left(1+\frac{2R}{y_d}\right)\right)^{\beta} 
\left( \log\left(1+\frac{|x-y|}{(x_d\vee y_d)\wedge |x-y|}\right)\right)^{\delta}
$$
and
$$
g(y; \beta , q, \delta, x)
:=\left(\frac{x_d}{|x-y|} \wedge 1 \right)^q |x-y|^{-d+\alpha}  
\left(\log\left(1+\frac{2R}{y_d}\right)\right)^{\beta}  
\left( \log\left(1+\frac{|x-y|}{(x_d\vee y_d)\wedge |x-y|}\right)\right)^{\delta}.
$$
Note that for $0<y_d<R$ we have that $\log(1+2 R/y_d)\asymp \log(2 R/y_d)$. 
In almost all our applications of Lemma \ref{l:key} and Corollary \ref{c:key}  below, the parameter $\delta$ will be 0. The only exception is Proposition \ref{p:bound-for-integral-new} where we will have $\delta$ equal to $0$, $\beta_4$ or $\beta_4+1$. 
\begin{lemma}\label{l:key}
Let $R \in (0, \infty)$ 
and $x=(\wt{0}, x_d)$ with $x_d \le 
2R/3$. 
Fix $0<a_1\le  x_d/2$ and $3x_d/2\le a_3\le a_2\le R$. 
We have the following comparison relations, with comparison constants independent of $R, a_1, a_2, a_3$ and 
$x_d\in (0, 2R/3)$:
\begin{itemize}
\item[(i)] 
If $\gamma > -1$ and $q>\alpha-1$, then
$$
I_1:=\int_{D(R, a_1)} f(y; \gamma,  \beta,  q, \delta, x) \, dy 
\asymp x_d^{\alpha-q-1}a_1^{\gamma+1}  
\left(\log \frac{2R}{a_1}\right)^{\beta}. 
$$
\item[(ii)] If $q>\alpha-1$, then
\begin{align*}
I_2&:=\int_{D(R, a_2)\setminus D(R, a_3)} f(y; \gamma,  \beta, q, \delta, x) \, dy\\
&\asymp 
R^{\gamma+\alpha-q}\left(F\left(\frac{a_3}{R};\gamma+\alpha-q-1, \beta\right)-F\left(\frac{a_2}{R};\gamma+\alpha-q-1, \beta\right)\right).
\end{align*}
\item[(iii)] 
 If $q>\alpha-1$, then 
$$
I_3:=\int_{D(R,3 x_d/2)\setminus 
D(R, x_d/2)} g(y; \beta, q, \delta, x) \, dy \asymp x_d^{\alpha} 
\left(\log \frac{2R}{x_d}\right)^{\beta}. 
$$
\end{itemize}
\end{lemma}

\pf (i) 
In $D(R, a_1)$, $y_d<x_d$. 
Without loss of generality, we replace $\log(1+2 R/y_d)$ with $\log (2 R/y_d)$. 
Thus, using the change of 
variables $y_d=x_d h$ and $r=x_d s$ in the second line below, we get
\begin{align*}
&I_1
\asymp  \int_0^R r^{d-2} \int_0^{a_1} \frac{y_d^{\gamma}}{((x_d-y_d)+r)^{d-\alpha+q}}  
\left(\log\frac{2R}{y_d}\right)^{\beta} 
\left(\log\left(1+\frac{(x_d-y_d)+r}{x_d}\right)\right)^{\delta}dy_d\, dr\\
&= x_d^{\alpha-q+\gamma} \int_0^{R/x_d} s^{d-2} \int_0^{a_1/x_d} \frac{h^{\gamma}}{[(1-h)+s]^{d-\alpha+q}} 
\left(\log\frac{2R/x_d}{h}\right)^{\beta} 
(\log(2-h+s))^{\delta} dh\, ds, 
\end{align*}
which, using $1-h\asymp 1$ (because $0<a_1\le  x_d/2$), is  comparable to 
$$ 
x_d^{\alpha-q+\gamma} \int_0^{R/x_d} \frac{s^{d-2}(\log(2+s))^\delta}{(1+s)^{d-\alpha+q}}ds \left(\int_0^{a_1/x_d} h^{\gamma} 
\left(\log\frac{2R/x_d}{h}\right)^{\beta} dh\right). 
$$
Note that, since  $q>\alpha-1$, 
$$
\int_1^{3/2}  \frac{(\log(2+s))^\delta}{s^{2-\alpha+q}}\, ds 
\le \int_1^{R/x_d}  \frac{(\log(2+s))^\delta}{s^{2-\alpha+q}}\, ds \le \int_1^{\infty}  \frac{(\log(2+s))^\delta}{s^{2-\alpha+q}}\, ds <\infty.
$$
Therefore, using 
this inequality
and  \eqref{e:F2}, we get after a change of variables 
\begin{align*}
I_1&\asymp  x_d^{\alpha-q+\gamma}
\left( \int^1_0\frac{s^{d-2}(\log(2+s))^\delta}{(1+s)^{d-\alpha+q}}ds+ \int_1^{R/x_d}  \frac{(\log(2+s))^\delta}{s^{2-\alpha+q}}\, ds\right) \times  \\
&\quad  \times\left(\frac{R}{x_d}\right)^{\gamma+1}\left(F(0; \gamma, \beta)-F\left(\frac{a_1}{R}; \gamma, \beta\right)\right)  \\
&\asymp  x_d^{\alpha-q-1} a_1^{\gamma+1} 
\left(\log \frac{2R}{a_1}\right)^{\beta} .
\end{align*}

\noindent (ii) 
In $D(R, a_2)\setminus D(R, a_3)$, $y_d>x_d$. 
Thus, using the change of variables $y_d=x_d h$ and $r=x_d s$ in the second line below, 
we get
\begin{align*}
I_2&\asymp  \int_0^R r^{d-2} \int_{a_3}^{a_2} \frac{y_d^{\gamma}}{((y_d-x_d)+r)^{d-\alpha+q}} 
\left(\log\frac{2R}{y_d}\right)^{\beta} 
\left(\log\left(1+\frac{(y_d-x_d)+r}{y_d}\right)\right)^{\delta}dy_d\, dr \nn\\
&= x_d^{\alpha-q+\gamma}   \int_{a_3/x_d}^{a_2/x_d} \int_0^{R/x_d} \frac{s^{d-2}h^{\gamma}}{[(h-1)+s]^{d-\alpha+q}} 
\left(\log \frac{2R/x_d}{h}\right)^{\beta}  
\left(\log\left(1+\frac{h-1+s}{h}\right)\right)^{\delta} 
ds\, dh,
\end{align*}
which is, by the change of variables $s=(h-1)t$, equal to 
\begin{align}\label{e:lnew1}
&x_d^{\alpha-q+\gamma} \int_{a_3/x_d}^{a_2/x_d} \int_0^{\frac{R}{(h-1)x_d}}\frac{h^{\gamma}t^{d-2}}{(h-1)^{1-\alpha+q}
(1+t)^{d-\alpha+q}} 
\left(\log \frac{2 R/x_d}{h}\right)^{\beta} \nonumber\\
&\quad \times
  \left(\log\left(1+\frac{(h-1)(1+t)}{h}\right)\right)^{\delta}
dt\, dh.
\end{align}

Note that, since 
$3x_d/2  \le a_3\le h x_d\le a_2 \le R$
 we have
$$
\frac{R}{(h-1)x_d} \ge \frac{R}{a_2-x_d} \ge 1, \quad
 a_3/x_d\le h \le a_2/x_d.
$$ 
Thus, using $q>\alpha-1$, we have that for $ a_3/x_d\le h \le a_2/x_d$,
$$
\int_{1/2}^{1}\frac{(\log(2+t))^{\delta}}{(1+t)^{2-\alpha+q}} dt  \le
\int_{1/2}^{\frac{R}{(h-1)x_d}}\frac{(\log(2+t))^{\delta}}{(1+t)^{2-\alpha+q}} dt \le
\int_{1/2}^{\infty}\frac{(\log(2+t))^{\delta}}{(1+t)^{2-\alpha+q}} dt <\infty.
$$ 
Therefore, using $(h-1)/h\asymp 1$ and the display  
above,  \eqref{e:lnew1} is comparable to 
\begin{align}
&x_d^{\alpha-q+\gamma} \int_{a_3/x_d}^{a_2/x_d} 
h^{\gamma+\alpha-q-1}
\left(\log \frac{2R/x_d}{h}\right)^{\beta} 
\int_0^{\frac{R}{(h-1)x_d}}\frac{t^{d-2}}{
(1+t)^{d-\alpha+q}}\left(\log (2+t)\right)^{\delta}
dt\, dh \nn \\
\asymp& x_d^{\alpha-q+\gamma} \int_{a_3/x_d}^{a_2/x_d} h^{\gamma+\alpha-q-1}
\left(\log \frac{2R/x_d}{h}\right)^{\beta}
\left(\int_0^{1/2} t^{d-2} dt
+\int_{1/2}^{\frac{R}{(h-1)x_d}}\frac{(\log(2+t))^{\delta}}{(1+t)^{2-\alpha+q}} dt
\right)dh \nn \\
\asymp&
 x_d^{\alpha-q+\gamma} \int_{a_3/x_d}^{a_2/x_d} h^{\gamma+\alpha-q-1}
\left(\log \frac{2 R/x_d}{h}\right)^{\beta} \, dh 
\nn\\
\asymp & 
R^{\gamma+\alpha-q}\left(F\left(\frac{a_3}{R};\gamma+\alpha-q-1, \beta\right)-F\left(\frac{a_2}{R};\gamma+\alpha-q-1, \beta\right)\right). 
\nn 
\end{align}

\noindent (iii) 
Let $B(x)=\{(\wt y, y_d): |\wt y|<x_d/2, |y_d -x_d|<x_d/2\}$.
Note that
$$
I_3=\int_{
 B(x) 
}g(y; \beta, q, 
\delta, x)dy+\int_{(D(R,3 x_d/2)\setminus D(R, x_d/2))\setminus 
 B(x)}g(y; \beta,  q, 
\delta, x)dy=:I_{31}+I_{32}.
$$
Note that in both $I_{31}$ and $I_{32}$ we have that 
$\log 2 R/y_d \asymp \log 2 R/x_d$ (since $y_d\asymp x_d$), 
and therefore this term comes out of the integral.
When $y\in  B(x) $, $x_d\asymp y_d\ge |x-y|$ so that 
$\left( \log\left(1+\frac{|x-y|}{(x_d\vee y_d)\wedge |x-y|}\right)\right)^{\delta}\asymp 1$. 
Therefore
$$
I_{31}\asymp 
\left(\log \frac{2 R}{x_d}\right)^{\beta} 
\int_{
 B(x)}|x-y|^{-d+\alpha}dy \asymp x_d^{\alpha} 
\left(\log \frac{2 R}{x_d}\right)^{\beta}.  
$$

In $(D(R,3 x_d/2)\setminus D(R, x_d/2))\setminus B(x) $, we have $y_d\asymp x_d$ and $x_d\le 2 |x-y|$. Thus, using the change of 
variables $y_d=rt+x_d$ in the third line below, we get
\begin{align*}
I_{32}&\asymp x_d^q   
\left(\log \frac{2 R}{x_d}\right)^{\beta} 
\int_{(D(R,3 x_d/2)\setminus D(R, x_d/2))\setminus 
 B(x) } |x-y|^{-d+\alpha-q}  \left(\log\left(1+ \frac{|x-y|}{x_d}\right)\right)^{\delta}dy\\
&\asymp x_d^q  
\left(\log \frac{2 R}{x_d}\right)^{\beta}
\int_{x_d/2}^R r^{d-2}\int_{x_d/2}^{3x_d/2}   {(|x_d-y_d|+r)^{-d+\alpha-q}}\left(\log\left(1+\frac{|x_d-y_d|+r}{x_d}\right)\right)^{\delta} dy_d\, dr\\
&=x_d^q  
\left(\log \frac{2 R}{x_d}\right)^{\beta} 
\int_{x_d/2}^R r^{\alpha-q-1}\int_{-\frac{x_d}{2r}}^{\frac{x_d}{2r}}(|t|+1)^{-d+\alpha-q}\left(\log\left(1+\frac{r(|t|+1)}{x_d}\right)\right)^{\delta}dt\, dr,
\end{align*}
which is, by the change of variables $r=x_d s$, comparable to
\begin{align}\label{e:lnew2}
x_d^{\alpha}  
\left(\log \frac{2 R}{x_d}\right)^{\beta} 
\int_{1/2}^{R/x_d}  s^{\alpha-q-1}\int_0^{1/s} 
\frac{\left(\log\big(1+s(t+1)\big)\right)^{\delta}}
{(t+1)^{d-\alpha+q}}dt\, ds.
\end{align}
Note that, since $q>\alpha-1$,
$$
\int_0^{1/s} 
\frac{\left(\log\big(1+s(t+1)\big)\right)^{\delta}}
{(t+1)^{d-\alpha+q}}dt\ \asymp (\log(1+s))^{\delta} \int_0^{1/s} \frac{dt}{(t+1)^{d-\alpha+q}} \asymp \frac{(\log(1+s))^{\delta}}{s}, 
\quad s >1/2
$$
and 
$$
 \int_{1/2}^{3/2}  
 \frac{(\log(1+s))^{\delta}}{s^{q+2-\alpha}}ds  \le
 \int_{1/2}^{R/x_d}   \frac{(\log(1+s))^{\delta}}{s^{q+2-\alpha}}ds  \le
 \int_{1/2}^{\infty}   \frac{(\log(1+s))^{\delta}}{s^{q+2-\alpha}}ds  <\infty.
 $$
Therefore, using the above inequalities,  \eqref{e:lnew2} is comparable to 
\begin{align*} 
& x_d^{\alpha}   
\left(\log \frac{2 R}{x_d}\right)^{\beta} 
\int_{1/2}^{R/x_d}  \frac{(\log(1+s))^{\delta}}{s^{q+2-\alpha}} ds \asymp x_d^{\alpha}  
\left(\log \frac{2 R}{x_d}\right)^{\beta}. 
\end{align*}
\qed

\begin{remark}\label{r:I1}
Note that it follows from the proof of Lemma \ref{l:key} (i) that 
$I_1=\infty$ for $\gamma\le -1$.
\end{remark}

\begin{corollary}\label{c:key}
Let $R>0$, 
$q>\alpha-1$, $\delta \in \R$,  $\gamma >-1$, $\beta\ge 0$,  and $x=(\wt{0}, x_d)$. 

 \noindent
(i) We have the following comparison result, with the comparison constant independent of $R$ and $x_d\in (0, R/2)$:
$$
\int_{D(R,R)} \left(\frac{x_d}{|x-y|}\wedge 1\right)^q
f(y;  \gamma, \beta, 0, \delta, x)\, dy
\asymp  \left\{ \begin{array}{ll}
R^{\alpha+\gamma-q} x_d^q, & 
 \text{ if } \alpha-1<q <\alpha+ \gamma;
\\
x_d^q 
\left(\log
\frac{2 R}{x_d}\right)^{\beta+1},  &
 \text{ if } q =\alpha+ \gamma;
\\
x_d^{\alpha+\gamma} 
\left(\log \frac{2 R}{x_d}\right)^{\beta} ,
& 
 \text{ if } q>\alpha+\gamma.
\end{array}\right.
$$

(ii) 
Let $a\in (0, R]$ and $\alpha-1 <q < \alpha+\gamma$.
Then there is a constant $C_{23}$ independent of 
$R$, $a$ and  $x_d\in (0, R/2)$ such that 
\begin{equation}\label{e:lemma-5-7}
\int_{D(R,a)}\left(\frac{x_d}{|x-y|}\wedge 1\right)^q 
f(y;  \gamma, \beta,  0, \delta, x)\, dy
\le C_{23}x_d^q a^{\alpha+\gamma -q}  
(\log 2 R/a)^{\beta}.  
\end{equation}
\end{corollary}
\pf (i) 
Set $a_1=x_d/2$, $a_2=R$ and $a_3=3x_d/2$ in Lemma \ref{l:key}. In $D(R, x_d/2)$
and $D(R, R)\setminus D(R, 3x_d/2)$, we have  $x_d \le c|x-y|$. 
Therefore,
\begin{align*}
&\int_{D(R, x_d/2)} 
\left(\frac{x_d}{|x-y|}\wedge 1\right)^q 
f(y; \gamma,\beta, 0, \delta, x)\, dy
\\
&\asymp x_d^q 
\int_{D(R, x_d/2)} 
f(y;  \gamma,\beta, q, \delta, x)\, dy
\asymp x_d^{\alpha+\gamma} 
\left(\log\frac{2R}{x_d}\right)^{\beta} .
\end{align*}
Using $3x_d/2 <3R/4$ (so that $3x_d/2R \le 3/4$), \eqref{e:F1} and \eqref{e:F3}, we get  
\begin{align*}
& \int_{D(R, R)\setminus D(R, 3x_d/2)} 
\left(\frac{x_d}{|x-y|}\wedge 1\right)^q 
f(y;  \gamma,\beta, 0, \delta, x)\, dy\\
&\asymp x_d^q 
 \int_{D(R, R)\setminus D(R, 3x_d/2)} 
f(y;  \gamma, \beta, q, \delta, x)\, dy\\
&\asymp x_d^q R^{\gamma+\alpha-q}F\left(\frac{3x_d}{2R}; \gamma+\alpha-q-1, \beta\right)  \\
& \asymp \left\{ \begin{array}{lll}
x_d^q R^{\alpha+\gamma-q}, &
\text{ if } \alpha-1<q <\alpha+ \gamma;\\
x_d^q \left(\log\frac{2 R}{x_d}\right)^{\beta+1}, 
 & \text{ if } q =\alpha+ \gamma;\\
x_d^{\alpha+\gamma} \left(\log\frac{2 R}{x_d}\right)^{\beta},  & 
\text{ if } q >\alpha+ \gamma.\\
\end{array}\right.\\
\end{align*}
In $D(R,3 x_d/2)\setminus D(R, x_d/2)$ we have that $y_d\asymp x_d$, so 
\begin{align*}
&\int_{D(R,3 x_d/2)\setminus D(R, x_d/2)}\left(\frac{x_d}{|x-y|}\wedge 1\right)^q
f(y;  \gamma, \beta, 0, \delta, x)\, dy\\
& \asymp x_d^{\gamma} 
\left(\log\frac{2 R}{x_d}\right)^{\beta} 
\int_{D(R,3 x_d/2)\setminus D(R, x_d/2)} g(y; q,\delta, x)\, dy 
\asymp x_d^{\alpha+\gamma} 
\left(\log\frac{2 R}{x_d}\right)^{\beta}. 
\end{align*}
By adding up these three displays we get the claim. 

\noindent (ii) If $a \le x_d/2$, then by 
Lemma \ref{l:key} (i) (with $a_1=a$)
and the assumption $\alpha-q-1<0$, we get that the integral in \eqref{e:lemma-5-7} is less than $c x_d^q (x_d^{\alpha-q-1} a^{\gamma+1} 
(\log 2 R/a)^{\beta})  \le x_d^q a^{\alpha+\gamma -q} 
 (\log2 R/a)^{\beta} $. If $x_d/2\le a \le 3x_d/2$, we split the integral into two parts -- over $D(R, x_d/2)$ and $D(R,a )\setminus D(R, x_d/2)$. The first one is by Lemma  \ref{l:key} (i) comparable with $x_d^q x_d^{\alpha-q+\gamma}  
(\log 4 R/x_d)^{\beta}  \asymp x_d^q a^{\alpha+\gamma -q} 
(\log 2 R/a)^{\beta}$,  while the second one is by Lemma \ref{l:key} (iii)  smaller than $x_d^{\gamma} x_d^{\alpha} 
 (\log 2 R/x_d)^{\beta} = x_d^q x_d^{\alpha+\gamma-q}  
(\log 2 R/x_d)^{\beta} \asymp x_d^q a^{\alpha+\gamma -q}  
(\log 2 R/a)^{\beta}$.   Finally, if $a\in (3x_d/2, R]$, then by using 
Lemma \ref{l:key} (ii) (with $a_2=a, a_3=3x_d/2$) and the assumption 
$q < \alpha+\gamma$
we get that the integral over $D(R,a)\setminus D(R, 3x_d/2)$ is bounded by above by $cx_d^q a^{\alpha+\gamma-q}  
(\log 2 R/a)^{\beta}$. 
\qed

\subsection{
Green function upper bound
for $p\in ((\alpha-1)_+, \alpha+\frac12[\beta_1+(\beta_1 \wedge \beta_2)])$}\label{ss:thm1.1(1)}
In this subsection we deal with the case
\begin{align}\label{e:fa1}
p\in ((\alpha-1)_+, \alpha+2^{-1}[\beta_1+ (\beta_1 \wedge \beta_2)]).
\end{align}
If $\beta_2>0$, then there exists $0<\wt \beta_2<\beta_2$ such that
\begin{align}\label{e:fa1-new}
p\in ((\alpha-1)_+, \alpha+2^{-1}[\beta_1+ (\beta_1 \wedge \wt \beta_2)]).
\end{align}
Further, 
if $\beta_4>0$,  there is $c>0$ such that for all $s\in (0,1)$
\begin{equation}\label{e:beta4-log}
s^{\beta_2}\log\left(1+\frac{8}{s}\right)^{\beta_4}\le c s^{\wt \beta_2}.
\end{equation} 

Let 
$$
\eps_0=
\begin{cases}
0& \text{ if }\beta_3=0 ;\\
2^{-1}( \alpha +\beta_1 -p) 
 & \text{ if }\beta_3>0.
\end{cases}
$$
Note that 
\begin{align}\label{e:flog}
[\log(1+s)]^{\beta_3} \le c s^{\eps_0}, \quad s\ge 1.
\end{align}

Recall
$$
D_{\wt{w}}(a,b)=\{x=(\wt{x}, x_d)\in \R^d:\, |\wt{x}-\wt{w}|<a, 0<x_d<b\}.
$$

\begin{lemma}\label{l:f3}
Suppose  that \eqref{e:fa1} holds.
There exists $C_{24}>0$ such that  for all $x, y\in \R^d_+$ with
$|\wt x- \wt y|>3$ and $0<x_d, y_d<1/4$, 
\begin{align}\label{e:f30}
&\int_{D_{\wt{x}}(1,1)}
 \int_{D_{\wt{y}}(1,1)}
\left(  \frac{x_d}{|w-x|} \wedge 1\right)^p
\left(  \frac{y_d}{|z-y|} \wedge 1\right)^p
\frac{(w_d \wedge z_d)^{ \beta_1}(w_d \vee z_d)^{ \beta_2}}{|x-w|^{d-\alpha}|y-z|^{d-\alpha}} 
\nn\\
& \qquad \times
\left( \log\left(1+\frac{w_d \vee z_d}{w_d \wedge z_d}\right)\right)^{ \beta_3}  \left(\log\left(1+\frac{8}{w_d\vee z_d}\right)\right)^{\beta_4} 
dz dw \le C_{24} x_d^p y_d^p.
\end{align}
\end{lemma}
\pf
Define 
$\wh \beta_1= \beta_1-\eps_0$, $\wh \beta_2=   \wt \beta_2  +\eps_0$. Note that by the definition of $\eps_0$, $p < \alpha+\wh \beta_1$.  Note first that by \eqref{e:beta4-log} we can estimate $(w_d\vee z_d)^{\beta_2}
\left(\log(1+8/(w_d\vee z_d))\right)^{\beta_4}$
 by a constant times $(w_d\vee z_d)^{\wt \beta_2}$.  
By \eqref{e:flog} and Tonelli's theorem, the left hand side of \eqref{e:f30} is less than or equal to 
\begin{align*}
&c_1 \int_{D_{\wt{x}}(1,1)}
 \int_{D_{\wt{y}}(1,1)}
\left(  \frac{x_d}{|w-x|} \wedge 1\right)^p
\left(  \frac{y_d}{|z-y|} \wedge 1\right)^p
\frac{(w_d \wedge z_d)^{\wh \beta_1}(w_d \vee z_d)^{\wh \beta_2}}{|x-w|^{d-\alpha}|y-z|^{d-\alpha}} dzdw\\
&=c_1 \left(\int_{\{(z, w ) \in D_{\wt{x}}(1,1) \times D_{\wt{y}}(1,1): z_d < w_d\} }
+\int_{\{(z, w ) \in D_{\wt{x}}(1,1) \times D_{\wt{y}}(1,1): z_d  \ge w_d\} }
\right) \\
&\quad \times \left(  \frac{x_d}{|w-x|} \wedge 1\right)^p
\left(  \frac{y_d}{|z-y|} \wedge 1\right)^p
\frac{(w_d \wedge z_d)^{\wh \beta_1}(w_d \vee z_d)^{\wh \beta_2}}{|x-w|^{d-\alpha}|y-z|^{d-\alpha}} dzdw \\
&=c_1 \int_{D_{\wt{x}}(1,1)}\left(  \frac{x_d}{|w-x|} \wedge 1\right)^p
\frac{w_d^{\wh \beta_2}}{|x-w|^{d-\alpha}} \left(\int_{D_{\wt{y}}(1, w_d) }
\left(  \frac{y_d}{|z-y|} \wedge 1\right)^p
\frac{z_d^{ \wh \beta_1}dz}{|y-z|^{d-\alpha}} \right) dw \\
&\quad+ c_1 \int_{D_{\wt{y}}(1,1)} 
\left(  \frac{y_d}{|z-y|} \wedge 1\right)^p
\frac{z_d^{\wh \beta_2}}{|y-z|^{d-\alpha}} \left( \int_{D_{\wt{x}}(1, z_d)}
\left(  \frac{x_d}{|w-x|} \wedge 1\right)^p
\frac{ w_d^{ \wh \beta_1}  dw}{|x-w|^{d-\alpha}} \right) dz.
\end{align*}
By symmetry, we only need to bound the last term above.

Since $\wh \beta_1+\alpha > p > \alpha-1$, 
we can apply 
Corollary \ref{c:key} (ii) (with $R=1$, $a=z_d$, 
$q=p$, $\gamma=\wh \beta_1$ and $\beta=\delta=0$) and get
\begin{align*}
& \int_{D_{\wt{y}}(1,1)} 
\left(  \frac{y_d}{|z-y|} \wedge 1\right)^p
\frac{z_d^{\wh \beta_2}}{|y-z|^{d-\alpha}} \left( \int_{D_{\wt{x}}(1, z_d)}
\left(  \frac{x_d}{|w-x|} \wedge 1\right)^p
\frac{ w_d^{ \wh \beta_1}  dw}{|x-w|^{d-\alpha}} \right) dz \nn\\
\le & c_4 x_d^p\int_{D_{\wt{y}}(1,1)} 
\left(  \frac{y_d}{|z-y|} \wedge 1\right)^p
\frac{z_d^{ \wt \beta_2 +\alpha+\beta_1-p}}{|y-z|^{d-\alpha}}  dz. 
\end{align*}
By \eqref{e:fa1} we have that 
$$
 ( \wt  \beta_2  +\alpha+\beta_1-p)+\alpha > p. 
$$
Thus, we can apply Corollary \ref{c:key} (ii) again (with $R=1$, $a=1$, 
$q=p$, $\gamma= \wt \beta_2  +\alpha+\beta_1-p$  and $\beta=\delta=0$)
and  conclude that
$$\int_{D_{\wt{y}}(1,1)} 
\left(  \frac{y_d}{|z-y|} \wedge 1\right)^p
\frac{z_d^{\wh \beta_2}}{|y-z|^{d-\alpha}} \left( \int_{D_{\wt{x}}(1, z_d)}
\left(  \frac{x_d}{|w-x|} \wedge 1\right)^p
\frac{ w_d^{ \wh \beta_1}  dw}{|x-w|^{d-\alpha}} \right) dz \le c_5 x_d^py_d^p.
$$
\qed

\begin{lemma}\label{l:f4}
Suppose \eqref{e:fa1} holds.
There exists $C_{25}>0$ such that 
for all $x, y\in \R^d_+$ with $|\wt x- \wt y|>4$ and $0<x_d, y_d<1/4$,
$$
G(x, y) \le C_{25} x_d^p y_d^p.
$$
\end{lemma}
\pf 
Assume $x=(\wt 0, x_d)$ with $0<x_d<1/4$, and
let $D=D(1,1)$ and $V= D_{\wt y}(1,1)$.
By Lemma \ref{l:prelub},
$$
G(w, y) \le c_1 \left(\frac{y_d}{|w-y|} \wedge 1\right)^p \le c_2 
y_d^p, \quad w \in \R^d \setminus  V.
$$
Thus
by Lemma \ref{l:new-lemma}, 
\begin{align*}
\E_x\left[G(Y_{\tau_{D}}, y); Y_{\tau_{D}}\notin V\right] \le 
c_3 y_d^p \P_x(Y_{\tau_{D}}  \in \R^d_+ )\le c_4  y^p_d x^p_d.
\end{align*}

On the other hand, 
since $2<|z-w|<8$  for  $(w, z) \in D \times V$, we have that $\log(1+\frac{|z-w|}{(w_d\vee z_d)\wedge |z-w|})\le \log(1+\frac{8}{w_d\vee z_d})$, and thus 
$$
J (w, z) \le  
c_5 (w_d \wedge z_d)^{ \beta_1}(w_d \vee z_d)^{ \beta_2}
\left(\log\Big(1+\frac{w_d \vee z_d}{w_d \wedge z_d}\Big)\right)^{ \beta_3}  \left( \log\Big(1+\frac{8}{w_d\vee z_d}\Big)\right)^{\beta_4}, 
   \quad (w, z) \in D \times V.
$$
By using the L\'evy system formula \eqref{e:levys3} (with $f= G(\cdot, y)$) in the first equality, and \eqref{e:Gu} in the third line,
\begin{align*}
&\E_x\left[G(Y_{\tau_{D}}, y); Y_{\tau_{D}}\in  V\right]\\
=&\int_D  G^D (x,w) \int_V  J (w, z)  G(z,y)dzdw
\le
\int_D G (x,w) \int_V   J(w, z) G(z,y)dzdw\\
\le 
& c_8 \int_{D}\left(  \frac{x_d}{|w-x|} \wedge 1\right)^p
\frac{1}{|x-w|^{d-\alpha}} \int_{V}
(w_d \wedge z_d)^{ \beta_1}(w_d \vee z_d)^{ \beta_2}
\left(\log\left(1+\frac{w_d \vee z_d}{w_d \wedge z_d}\right)\right)^{ \beta_3} 
\nn\\
& \qquad \times
\left(\log\Big(1+\frac{8}{w_d \vee z_d}\Big)\right)^{\beta_4} 
\left(  \frac{y_d}{|z-y|} \wedge 1\right)^p
\frac{dz}{|y-z|^{d-\alpha}} dw,
\end{align*}
which is less than or equal to $c_6x^p_dy^p_d$ by Lemma  \ref{l:f3}.
Therefore
\begin{align*}
G(x, y)=\E_x\left[G(Y_{\tau_{D}}, y); Y_{\tau_{D}}\notin V\right]+ \E_x\left[G(Y_{\tau_{D}}, y); Y_{\tau_{D}}\in V\right]\le 
c_{7}x^p_dy^p_d.
\end{align*}
\qed

\subsection{Green function estimates for $p\in [\alpha+\frac{\beta_1+\beta_2}2, \alpha+\beta_1)$}\label{ss:thm1.1(2)}
In this subsection we deal with the case
\begin{align}\label{e:fa1e}
\alpha+\frac{\beta_1+\beta_2}2 \le p <\alpha+\beta_1.
\end{align}
Note that \eqref{e:fa1e} implies $\beta_2 < \beta_1$ and
\begin{align}\label{e:fa3e2}
\alpha+ \beta_2 <  p,
\end{align}
\begin{align}\label{e:fa3e}
2\alpha-2p+\beta_1+\beta_2 \le 0.
\end{align}

Recall that $B^+_a:=B(0, a)\cap \R^d_+$, $a>0$.
The lower bound in the following theorem sharpens the lower bound in Lemma \ref{l:GB_4} under the assumption  \eqref{e:fa1e}. 

\begin{thm}\label{t:f4e}
Suppose \eqref{e:fa1e} holds.
For every $\eps \in (0, 1/4)$, there exists a constant $C_{26}>0$
such that for all $w \in \partial \R^d_+$,  $R>0$ and $x,y \in 
B(w, (1-\eps)R)\cap \R^d_+$, it holds that
\begin{align*}
&G^{B(w, R)\cap \R^d_+}(x, y)\ge
\frac{C_{26}}{|x-y|^{d-\alpha}}
\left( \frac{ x_d \wedge y_d}{|x-y|} \wedge 1 \right)^p
\\
&\qquad \times\begin{cases}
(\frac{ x_d \vee y_d}{|x-y|}\wedge 1)^{2\alpha-p+\beta_1+\beta_2}  
\left(\log (1+\frac{|x-y|}{(x_d \vee y_d)\wedge |x-y|})\right)^{\beta_4}  \quad \text{ if } \alpha+\frac{\beta_1+\beta_2}2 < p <\alpha+\beta_1;&\\
(\frac{ x_d \vee y_d}{|x-y|}\wedge 1)^p 
\left(\log (1+\frac{|x-y|}{(x_d \vee y_d)\wedge |x-y|})\right)^{\beta_4+1  }, 
\quad \text{ if }  p=\alpha+\frac{\beta_1+\beta_2}2.&
 \end{cases}
\end{align*}
\end{thm}

\pf 
By scaling, translation and symmetry, without loss of generality, we assume that 
$w=0$,
$R=1$ and $x_d \le y_d$. Moreover, 
by Theorem \ref{t:GB}, we only need to show that  there exists a constant 
$c_1>0$ such that for
all $x,y \in B^+_{1-\eps}$ with $x_d \le y_d$ satisfying $|x-y|\ge (40/\eps) y_d $, it holds that
\begin{align}
\label{e:flower}
G^{B^+_1}(x, y)\ge  
\frac{c_1  x^p_d}{|x-y|^{d+\alpha+\beta_1+\beta_2}}
 \begin{cases}
y_d^{2\alpha-p+\beta_1+\beta_2}
\left(\log (|x-y|/y_d)\right)^{\beta_4}
 & \text{ if } 2\alpha-2p+\beta_1+\beta_2<0;\\
y_d^p 
\left(\log (|x-y|/y_d)\right)^{ \beta_4+1 } & \text{ if } 2\alpha-2p+\beta_1+\beta_2=0.
 \end{cases}\end{align} 
We assume that $x,y \in B^+_{1-\eps}$ with $x_d \le y_d$ satisfying $|x-y|\ge (40/\eps) y_d $.
By the 
Harnack inequality (Theorem  \ref{t:uhp}), 
we can further assume that $4 x_d  \le  y_d$. Let $M = 40/\eps$ and $r=4|x-y|/M$.

By the L\'evy system formula \eqref{e:levys3} (with $f=G^{B^+_1}(\cdot, y)$) 
 and regular harmonicity of 
$w \mapsto G^{B^+_1}(w, y)$ on $D_{\wt x}(2r,2r)$, 
\begin{align}
&G^{B^+_1}(x, y) \ge
\E_x\left[G^{B^+_1}(Y_{\tau_{D_{\wt x}(2r,2r)}}, y); 
 Y_{\tau_{D_{\wt x}(2,2)}}\in  D_{\wt y}(r,r)\right]\nn\\
=&\int_{D_{\wt x}(2r,2r)} G^{D_{\wt x}(2r,2r)} (x,w) \int_{D_{\wt y}(r,r)} 
J(w, z) G^{B^+_1}(z,y)dzdw\nn\\
 \ge &\int_{D_{\wt x}(r,r)} G^{D_{\wt x}(2r,2r)} (x,w) \int_{D_{\wt y}(r,r)} 
 J(w, z) G^{B^+_1}(z,y)dzdw\nn\\
 \ge &\int_{D_{\wt x}(r,r)} G^{B((\wt x, 0), 2r) \cap \R^d_+} (x,w) \int_{D_{\wt y}(r,r)} 
 J(w, z) G^{B^+_1}(z,y)dzdw. \label{e:nneww}
\end{align}

Since $D_{\wt x}(r,r) \subset B((\wt x, 0), \sqrt 2 r) \cap \R^d_+ $ and  
$ D_{\wt y}(r,r) \subset B^+_{(1-\eps/4)}$, 
we have by Theorem \ref{t:GB},
\begin{equation}\label{e:new1}
G^{B((\wt x, 0), 2r) \cap \R^d_+}(x,w) \ge c_2
\left(  \frac{x_d}{|w-x|} \wedge 1\right)^p\left(  \frac{w_d}{|w-x|} \wedge 1\right)^p
\frac{1}{|x-w|^{d-\alpha}}, \quad w \in D_{\wt x}(r,r), 
\end{equation}
and
\begin{equation}\label{e:new2}
G^{B^+_1}(z,y) \ge c_3
\left(  \frac{y_d}{|z-y|} \wedge 1\right)^p
\left(  \frac{z_d}{|z-y|} \wedge 1\right)^p
\frac{1}{|y-z|^{d-\alpha}}, \quad z \in D_{\wt y}(r,r).
\end{equation}
Moreover, since $(w_d \vee z_d) \le |z-w| \asymp r$  for $(w, z) \in D_{\wt x}(r,r) \times D_{\wt y}(r,r)$, we have
\begin{align}
J(w, z) &\ge c_4
|w-z|^{-d-\alpha}
\Big(\frac{w_d \wedge z_d}{|w-z|}\wedge 1\Big)^{\beta_1}
\Big(\frac{w_d \vee z_d}{|w-z|}\wedge 1\Big)^{\beta_2} 
\left(\log\Big(1+\frac{|w-z|}{(w_d\vee z_d)\wedge |w-z|} \Big)\right)^{\beta_4}   \nn\\
 & \ge c_5 \frac{(w_d \wedge z_d)^{ \beta_1}(w_d \vee z_d)^{ \beta_2}}{r^{d+\alpha+\beta_1+ \beta_2}}  \left(\log\frac{2r}{w_d\vee z_d}\right)^{\beta_4}, 
 \quad  (w, z) \in D_{\wt x}(r,r) \times D_{\wt y}(r,r).
 \label{e:new3}
 \end{align}
Using \eqref{e:new1}--\eqref{e:new3}, we obtain
\begin{align}
 &\int_{D_{\wt x}(r,r)} G^{B((\wt x, 0), 2r) \cap \R^d_+} (x,w) \int_{D_{\wt y}(r,r)} 
 J(w, z) G^{B^+_1}(z,y)dzdw\nn\\
\ge 
& \frac{c_6} {r^{d+\alpha+\beta_1+ \beta_2}}\int_{D_{\wt x}(r,r)}
\left(  \frac{x_d}{|w-x|} \wedge 1\right)^p\left(  \frac{w_d}{|w-x|} \wedge 1\right)^p
\frac{1}{|x-w|^{d-\alpha}} 
\nn \\
& \qquad \times \int_{D_{\wt y}(r,r)}
\left(  \frac{y_d}{|z-y|} \wedge 1\right)^p
\left(  \frac{z_d}{|z-y|} \wedge 1\right)^p
\frac{(w_d \wedge z_d)^{ \beta_1}(w_d \vee z_d)^{ \beta_2}}{|y-z|^{d-\alpha}}  \left(\log\frac{2r}{w_d\vee z_d}\right)^{\beta_4} dzdw \nn\\
\ge  
& \frac{c_7} {r^{d+\alpha+\beta_1+ \beta_2}}
\int_{D_{\wt y}(r,r)\setminus D_{\wt y}(r, 3y_d/2)} 
\left(  \frac{y_d}{|z-y|} \wedge 1\right)^p
\left(  \frac{z_d}{|z-y|} \wedge 1\right)^p
\frac{z_d^{ \beta_2}}{|y-z|^{d-\alpha}}  \left(\log\frac{2r}{z_d}\right)^{\beta_4} 
\nn \\
& \qquad \times \left( 
\int_{D_{\wt{x}}(r, z_d)}
\left(  \frac{x_d}{|w-x|} \wedge 1\right)^p
\left(  \frac{w_d}{|w-x|} \wedge 1\right)^p
\frac{ w_d^{  \beta_1}  dw}{|x-w|^{d-\alpha}} \right) dz
\nn\\
&\ge\frac{c_8x^p_dy^{p}_d} 
{r^{d+\alpha+\beta_1+ \beta_2}}\int_{D_{\wt y}(r,r)\setminus D_{\wt y}(r, 3y_d/2)} \frac{z_d^{p+\beta_2}}{|y-z|^{d+2p-\alpha}}  \left(\log\frac{2r}{z_d}\right)^{\beta_4} 
\nn \\
& \qquad \times
\left(\int_{D_{\wt{x}}(r, z_d)\setminus D_{\wt{x}}(r, 3x_d/2)}
\frac{ w_d^{  p+\beta_1}  dw}{|x-w|^{d+2p-\alpha}} \right) dz.\nn
\end{align}
Now by  applying
Lemma \ref{l:key} (ii) with $R=r$, $a_2=z_d$, $a_3=3x_d/2$, $\gamma=p+\beta_1$, $q=2p$
and $\beta=  \delta=0$ in the inner integral,  
we get that for $z_d \ge 3y_d/2$,
$$
\int_{D_{\wt{x}}(r, z_d)\setminus D_{\wt{x}}(r, 3x_d/2)}
\frac{ w_d^{  p+\beta_1}  dw}{|x-w|^{d+2p-\alpha}}\ge 
c_9(z_d^{\alpha-p+\beta_1}-(3x_d/2)^{\alpha-p+\beta_1})\ge c_{10}z_d^{\alpha-p+\beta_1}.
$$
In the last inequality above, we have used the the assumption  $4 x_d \le y_d$ so that 
for all $z_d \ge 3y_d/2$ it holds $z_d/4   \ge 3x_d/2$.
Thus, we have
\begin{align}
 &\int_{D_{\wt x}(r,r)} G^{B((\wt x, 0), 2r) \cap \R^d_+} (x,w) \int_{D_{\wt y}(r,r)} 
 J(w, z) G^{B^+_1}(z,y)dzdw\nn\\
\ge & \frac{c_{11}x^p_dy^{p}_d} 
{r^{d+\alpha+\beta_1+ \beta_2}}\int_{D_{\wt y}(r,r)\setminus D_{\wt y}(r, 3y_d/2)} \frac{z_d^{\beta_1+\beta_2+\alpha}}{|y-z|^{d+2p-\alpha}}   
\left(\log\frac{2r}{z_d}\right)^{\beta_4}   dz. \label{e:neww}
\end{align}
Finally, applying Lemma \ref{l:key} (ii) with $R=r$, $a_2=r$, $a_3=3y_d/2$, $\gamma=\alpha+\beta_1+\beta_2$, $q=2p$,
$\beta=\beta_4$ and $\delta=0$ 
and using the fact that $y_d<r/4$, we get that the above is 
greater 
than or equal to 
\begin{align}
\frac{c_{12}
 x^p_dy^{p}_d} {r^{d+\alpha+\beta_1+ \beta_2}}
 \left\{ \begin{array}{lll}
 y_d^{2\alpha-2p+\beta_1+\beta_2}  \left(\log\frac{ r }{y_d}\right)^{\beta_4},   & \text{ if } 2\alpha-2p+\beta_1+\beta_2<0;\\
   \left( \log \frac{ r }{y_d}\right)^{\beta_4+1 }, & \text{ if } 2\alpha-2p+\beta_1+\beta_2=0.
\end{array}\right.
 \label{e:nnewww}
\end{align}
Recalling that $r=4|x-y|/M$ and combining  \eqref{e:nneww}, \eqref{e:neww} and \eqref{e:nnewww}, we have proved that \eqref{e:flower} holds.
\qed

We now consider the upper bound of $G(x,y)$.

\begin{lemma}\label{l:f3qq}
Suppose \eqref{e:fa1e} holds.
There exists $C_{27}>0$ such that  for all $x, y\in \R^d_+$ with 
$|\wt x- \wt y|>3$, and $0<4x_d\le y_d<\frac14$ or $0<4y_d\le x_d<\frac14$,
\begin{align}
&\int_{D_{\wt{x}}(1,1)}dw\int_{D_{\wt{y}}(1,1)}dz\left(  \frac{x_d\wedge w_d}{|w-x|} \wedge 1\right)^p\frac{(w_d \wedge z_d)^{ \beta_1}(w_d \vee z_d)^{ \beta_2}}{|x-w|^{d-\alpha}|y-z|^{d-\alpha}} 
\nn \\
&\qquad \qquad\times
\left(\log\left(1+\frac{w_d \vee z_d}{w_d \wedge z_d}\right)\right)^{ \beta_3}  \left(\log\Big(1+\frac{2}{w_d\vee z_d}\Big)\right)^{\beta_4} 
\left(  \frac{y_d\wedge z_d}{|z-y|} \wedge 1\right)^p\nn\\
& \le  C_{27} (x_d\wedge y_d)^p
   \begin{cases}
(x_d \vee y_d)^{2\alpha-p+\beta_1+\beta_2}  \left( \log (1/(x_d \vee y_d))\right)^{\beta_4}  & \text{ if } 2\alpha-2p+\beta_1+\beta_2<0;\\
(x_d \vee y_d)^p \left(\log (1/(x_d \vee y_d)) \right)^{ \beta_4+1}   & \text{ if } 2\alpha-2p+\beta_1+\beta_2=0. 
 \end{cases} \label{e:f30qq}
\end{align}
\end{lemma}
\pf
By symmetry, we only need to consider the case $0\le 4 x_d  \le  y_d\le 1/4$.
Define 
$$\eps_0:=2^{-1}{\bf 1}_{\beta_3>0}
[( \alpha +\beta_1 -p) \wedge (p-\alpha-\beta_2)], \quad \wh \beta_1= \beta_1-\eps_0
\quad\text{ and } \quad \wh \beta_2= \beta_2+\eps_0.$$ 
 Note that $p<\alpha +\wh \beta_1$ and $p>\alpha +\wh \beta_2$ by  \eqref{e:fa3e2}.

By \eqref{e:flog}, 
\begin{align*}
&\int_{D_{\wt{x}}(1,1)}dw\int_{D_{\wt{y}}(1,1)}dz
\left(  \frac{x_d\wedge w_d}{|w-x|} \wedge 1\right)^p\frac{(w_d \wedge z_d)^{ \beta_1}(w_d \vee z_d)^{ \beta_2}}{|x-w|^{d-\alpha}|y-z|^{d-\alpha}} 
\nn\\
&\qquad\qquad \times
\left(\log\left(1+\frac{w_d \vee z_d}{w_d \wedge z_d}\right)\right)^{ \beta_3}  \left(\log\frac{2}{w_d\vee z_d}\right)^{\beta_4} 
\left(  \frac{y_d\wedge z_d}{|z-y|} \wedge 1\right)^p\\
 &\le  c_1 
\left(\int_{\{(z, w ) \in D_{\wt{x}}(1,1) \times D_{\wt{y}}(1,1): z_d < w_d\} }
+\int_{\{(z, w ) \in D_{\wt{x}}(1,1) \times D_{\wt{y}}(1,1): z_d  \ge w_d\} }
\right)
\\
&
\qquad\qquad 
\times 
\left(  \frac{x_d\wedge w_d}{|w-x|} \wedge 1\right)^p
\frac{(w_d \wedge z_d)^{\wh \beta_1}(w_d \vee z_d)^{\wh \beta_2}}{|x-w|^{d-\alpha}|y-z|^{d-\alpha}} \left(\log\frac{2}{w_d\vee z_d}\right)^{\beta_4}   \left(  \frac{y_d\wedge z_d}{|z-y|} \wedge 1\right)^p
dzdw \\
&\le c_1 \int_{D_{\wt{y}}(1,1)} 
\left(  \frac{y_d\wedge z_d}{|z-y|} \wedge 1\right)^p
\frac{z_d^{ \wh \beta_1}}{|y-z|^{d-\alpha}} 
\int_{D_{\wt{x}}(1, 1)\setminus D_{\wt{x}}(1, z_d)}
\left(  \frac{x_d}{|w-x|} \wedge 1\right)^p
\frac{ \left( \log({2} /{w_d})\right)^{\beta_4}  w_d^{  \wh \beta_2}  dw}{|x-w|^{d-\alpha}}
dz \\
&+ c_1 \int_{D_{\wt{y}}(1,1)} 
\left(  \frac{y_d}{|z-y|} \wedge 1\right)^p
\frac{ \left(\log({2}/{z_d})\right)^{\beta_4}  z_d^{ \wh \beta_2}}{|y-z|^{d-\alpha}}  \left( \int_{D_{\wt{x}}(1, z_d)}
\left(  \frac{x_d}{|w-x|} \wedge 1\right)^p
\frac{ w_d^{  \wh \beta_1}  dw}{|x-w|^{d-\alpha}} \right) dz\\ &=:I_1+I_2.
\end{align*}
Since $ \wh \beta_1> p-\alpha > \beta_2 \ge 0$, 
we can apply Corollary \ref{c:key} (ii) to estimate the inner integral in $I_2$ 
to get 
\begin{align}\label{e:hlp1}
I_2   \le c_2 x_d^p\int_{D_{\wt{y}}(1,1)} 
\left(  \frac{y_d}{|z-y|} \wedge 1\right)^p
\frac{z_d^{ \beta_2+\alpha + \beta_1-p}}{|y-z|^{d-\alpha}}  \left(\log\frac{2}{ z_d}\right)^{\beta_4}   dz. 
\end{align}
 By \eqref{e:fa3e},
$$
0 <\beta_2+ \alpha + \beta_1-p  \le  p-\alpha. 
$$
Thus we can apply 
Corollary \ref{c:key} (i) 
to get that 
(and by using $y_d<1/4$ we may replace 2 with 1) 
\begin{align}I_2
 \le c_3 x_d^p \begin{cases}
y_d^{2\alpha-p+\beta_1+\beta_2}\left(\log(1/{y_d})\right)^{\beta_4}   & \text{ if } 2\alpha-2p+\beta_1+\beta_2<0;\\
y_d^p \left(\log (1/y_d)\right)^{ \beta_4+1 } & \text{ if } 2\alpha-2p+\beta_1+\beta_2=0.
 \end{cases}
 \label{e:app1}
\end{align}

We now consider 
\begin{align*}
&I_1\\
&\le\int_{D_{\wt{y}}(1,2x_d)}
 \left(  \frac{z_d}{|z-y|} \wedge 1\right)^p
\frac{z_d^{ \wh \beta_1}}{|y-z|^{d-\alpha}} 
\int_{D_{\wt{x}}(1, 1)\setminus D_{\wt{x}}(1, z_d)}
\left(  \frac{x_d}{|w-x|} \wedge 1\right)^p
\frac{ \left(\log ({2} /{w_d})\right)^{\beta_4}  w_d^{  \wh \beta_2}  dw}{|x-w|^{d-\alpha}}dz\\
&+\int_{D_{\wt{y}}(1,1) \setminus D_{\wt{y}}(1,2x_d)} 
 \left(  \frac{y_d}{|z-y|} \wedge 1\right)^p
 \frac{z_d^{ \wh \beta_1}}{|y-z|^{d-\alpha}} 
\\
&
\qquad\qquad 
\times 
\int_{D_{\wt{x}}(1, 1)\setminus D_{\wt{x}}(1, z_d)}
\left(  \frac{x_d}{|w-x|} \wedge 1\right)^p
\frac{\left(\log ({2}/ {w_d})\right)^{\beta_4}  w_d^{  \wh \beta_2}  dw}{|x-w|^{d-\alpha}}dz\\
&\le\int_{D_{\wt{y}}(1,2x_d)}
\frac{z_d^{ \wh \beta_1 +p}}{|y-z|^{d-\alpha+p}}dz
  \int_{D_{\wt{x}}(1, 1)}
\left(  \frac{x_d}{|w-x|} \wedge 1\right)^p
\frac{ \left(\log ({2}/{w_d})\right)^{\beta_4}   w_d^{  \wh \beta_2}  dw}{|x-w|^{d-\alpha}} \\
&+x_d^p\int_{D_{\wt{y}}(1,1)} \left(  \frac{y_d}{|z-y|} \wedge 1\right)^p
\frac{z_d^{ \wh \beta_1}}{|y-z|^{d-\alpha}} 
 \int_{D_{\wt{x}}(1, 1)\setminus D_{\wt{x}}(1, z_d)}
\frac{ \left(\log ({2}/{w_d})\right)^{\beta_4}  w_d^{  \wh \beta_2}  dw}{|x-w|^{d-\alpha+p}}   dz
\\
&=:I_{11}+x_d^p I_{12}.
\end{align*}
Since  $p \ge \alpha$   and $4x_d \le y_d$, 
we can apply Lemma \ref{l:key} (i) (with $a_1=2x_d, \gamma=p+ \wh \beta_1, 
q=p, \beta=  \delta=0$) to get
$$
\int_{D_{\wt{y}}(1,2x_d)}
\frac{z_d^{ \wh \beta_1 +p}}{|y-z|^{d-\alpha+p}}dz
\le c_4y_d^{\alpha-p-1}x_d^{p+\wh \beta_1+1}.
$$
Since $\alpha+ \wh \beta_2 <  p$,  by Corollary \ref{c:key} (i) it follows that
$$
\int_{D_{\wt{x}}(1, 1)}
\left(  \frac{x_d}{|w-x|} \wedge 1\right)^p
\frac{ \left(\log ({2}/{w_d}\right)^{\beta_4}  w_d^{  \wh \beta_2}  dw}{|x-w|^{d-\alpha}} 
\le c_5x_d^{\alpha+\wh \beta_2}  \left(\log \frac{2}{x_d}\right)^{\beta_4}. 
$$
Thus, we have
\begin{align}
I_{11} &\le 
c_6y_d^{\alpha-p-1}x_d^{p+\wh \beta_1+1}x_d^{\alpha+\wh \beta_2} (\log(2/x_d))^{\beta_4} 
=c_6x_d^px_d^{\alpha+\beta_1+\beta_2+1} (\log(2/x_d))^{\beta_4}  y_d^{\alpha-p-1}\nn\\
&\le c_6 x_d^p y_d^{\alpha+\beta_1+\beta_2+1} (\log(2/y_d))^{\beta_4}    y_d^{\alpha-p-1}
 \le \wt{c_6}  x_d^py_d^{2\alpha-p+\beta_1+\beta_2}  (\log(1/y_d))^{\beta_4}.  
 \label{e:app2}
 \end{align}
Here we used that 
$t\mapsto t^{\alpha+\beta_1+\beta_2+1}\log(2/t)^{\beta_4}$
is almost increasing on $(0,1/4]$. 

Finally, we take care of $I_{12}$. Note that 
for every $z \in D_{\wt{y}}(1,1) \setminus D_{\wt{y}}(1,2x_d)$, we have $z_d > 2x_d$  and 
so, since $\alpha+ \wh \beta_2 <  p$,  by Lemma  
\ref{l:key} (ii) with $R=a_2=1, a_3=z_d, \gamma=\wh \beta_2, q=p, \beta=\beta_4, \delta=0$,
$$ 
\int_{D_{\wt{x}}(1, 1)\setminus D_{\wt{x}}(1, z_d)}
\frac{  \left(\log ({2}/{w_d})\right)^{\beta_4}   w_d^{  \wh \beta_2}  dw}{|x-w|^{d-\alpha+p}}
\le cz_d^{ \alpha+ \wh \beta_2-p}   \left(\log \frac{2}{z_d}\right)^{\beta_4}. 
$$
Thus, 
\begin{align*}
I_{12} \le   c_8 \int_{ D_{\wt{y}}(1,1)} \left(  \frac{y_d}{|z-y|} \wedge 1\right)^p
\frac{z_d^{ \beta_2+\alpha + \beta_1-p}}{|y-z|^{d-\alpha}}     \left(\log \frac{2}{z_d}\right)^{\beta_4}  dz.
\end{align*}
By the same argument as that in in \eqref{e:hlp1} and \eqref{e:app1}, we now have 
\begin{align}
I_{12}  \le c_9
\begin{cases}
y_d^{2\alpha-p+\beta_1+\beta_2}  (\log (1/y_d))^{\beta_4} & \text{ if } 2\alpha-2p+\beta_1+\beta_2<0;\\
y_d^p \left(\log (1/y_d)\right)^{  \beta_4+1 } & \text{ if } 2\alpha-2p+\beta_1+\beta_2=0.
 \end{cases}
  \label{e:app3}
 \end{align}
By combining \eqref{e:app1}--\eqref{e:app3} and symmetry, we have proved the lemma.
 \qed

  \begin{remark}\label{r:6.8}
 {\rm
 In the proof of Lemma \ref{l:f3qq},  if we had used  Tonelli's theorem on $I_1$ and  estimated  it 
as $I_2$ 
  (instead of using the argument to bound $I_{11}$ and $I_{12}$ separately),  
we would not have obtained the sharp upper bound.
}
 \end{remark}
 
\begin{prop}\label{p:f4ef}
Suppose \eqref{e:fa1e} holds.
There exists $C_{28}>0$ such that  for all 
for all $x, y\in \R^d_+$ with $0<x_d, y_d<1/4$ with $|\wt x- \wt y|>4$,
$$
G(x, y) \le C_{28}(  x_d \wedge y_d)^p\begin{cases}
(x_d \vee y_d)^{2\alpha-p+\beta_1+\beta_2}  \left(\log (1/(x_d \vee y_d))\right) ^{\beta_4}  & \text{ if } 2\alpha-2p+\beta_1+\beta_2<0;\\
(x_d \vee y_d)^p \left(\log (1/(x_d \vee y_d))\right)^{\beta_4+1 } & \text{ if } 2\alpha-2p+\beta_1+\beta_2=0.
 \end{cases}
$$
\end{prop}
\pf 
Without loss of generality, we assume $\wt x= \wt 0$. By symmetry, we consider the case $0<x_d \le y_d<1/4$ only. 
By the 
Harnack inequality (Theorem  \ref{t:uhp}), it suffices to deal with the case $0<4x_d \le y_d<1/4$.
Let $D=D(1,1)$ and $V= D_{\wt y}(1,1)$ 
By  the L\'evy system formula \eqref{e:levys3} (with $f=G(\cdot, y)$), 
\eqref{e:Gu}, Lemma \ref{l:f3qq},
and the fact that $2<|z-w|<8 $ below (so that $|z-w|\asymp 2$) 
\begin{align*}
&\E_x\left[G(Y_{\tau_{D}}, y); Y_{\tau_{D}}\in  V\right]\\
=&\int_D G^D (x,w) \int_V J(w, z) G(z,y)dzdw
\le 
\int_D G (x,w) \int_V J(w, z) G(z,y)dzdw\\
\le 
& c_1
\int_{D_{\wt{x}}(1,1)}dw\int_{D_{\wt{y}}(1,1)}dz\left(  \frac{x_d\wedge w_d}{|w-x|} \wedge 1\right)^p\frac{(w_d \wedge z_d)^{ \beta_1}(w_d \vee z_d)^{ \beta_2}}{|x-w|^{d-\alpha}|y-z|^{d-\alpha}} 
\nn \\
&\qquad \qquad\times
\left(\log\Big(1+\frac{w_d \vee z_d}{w_d \wedge z_d}\Big)\right)^{ \beta_3}  \left(\log\Big(1+\frac{2}{w_d\vee z_d}\Big)\right)^{\beta_4} 
\left(  \frac{y_d\wedge z_d}{|z-y|} \wedge 1\right)^p\\
\le &  c_2 (x_d\wedge y_d)^p \begin{cases}
(x_d \vee y_d)^{2\alpha-p+\beta_1+\beta_2}  \left(\log (1/(x_d \vee y_d))\right)^{\beta_4}  & \text{ if } 2\alpha-2p+\beta_1+\beta_2<0;\\
(x_d \vee y_d)^p \log (1/(x_d \vee y_d))^{\beta_4+1  } & \text{ if } 2\alpha-2p+\beta_1+\beta_2=0.
 \end{cases}
\end{align*}
Moreover, by the same argument as that in the proof of Lemma \ref{l:f4}, we also have 
\begin{align*}
\E_x\left[G(Y_{\tau_{D}}, y); Y_{\tau_{D}}\notin V\right] 
\le c_3 y_d^p \P_x(Y_{\tau_{D}}  \in \R^d_+ )\le c_4  y^p_d x^p_d. 
\end{align*}
Therefore
\begin{align*}
G(x, y)&=\E_x\left[G(Y_{\tau_{D}}, y); Y_{\tau_{D}}\notin V\right]+ \E_x\left[G(Y_{\tau_{D}}, y); Y_{\tau_{D}}\in V\right]\\
  & \le  c_5 (x_d\wedge y_d)^p \begin{cases}
(x_d \vee y_d)^{2\alpha-p+\beta_1+\beta_2} \left(\log (1/(x_d \vee y_d))\right)^{\beta_4}  & \text{ if } 2\alpha-2p+\beta_1+\beta_2<0;\\
(x_d \vee y_d)^p \left(\log (1/(x_d \vee y_d))\right)^{ \beta_4+1  }  & \text{ if } 2\alpha-2p+\beta_1+\beta_2=0.
 \end{cases}
\end{align*}
\qed

\subsection{Proof of Theorem \ref{t:Green} and estimates of  potentials}\label{ss:thm1.1}
With the preparations in the previous two subsections, we are now ready to prove
Theorem \ref{t:Green}.
We recall \cite[Theorem 3.14]{KSV} on the 
H\"older continuity of bounded harmonic functions: There exist constants 
$c>0$ and $\gamma\in (0,1)$ such that
 for every $x_0\in \R^d_+$, $r\in (0,1]$ such that $B(x_0,2r)\subset \R^d_+$ and every bounded  $f:\R^d_+\to [0,\infty)$ which is harmonic in $B(x_0,2r)$, it holds that
\begin{align}
\label{eHR}
|f(x)-f(y)|\le c
\|f\|_{\infty}\left(\frac{|x-y|}{r}\right)^{\gamma}\qquad \textrm{for all }x,y\in B(x_0,r).
\end{align}

\bigskip

\noindent \textbf{Proof of Theorem \ref{t:Green}}.
The existence and regular harmonicity of the 
Green function 
were shown in  Proposition \ref{p:existenceGF}. 
We prove now the continuity of $G$.
We fix $x_0, y_0 \in \R^d_+$ and choose a positive $a$ small enough  so that $B(x_0, 4a) \cap B(y_0, 4a)= \emptyset$ and
$B(x_0, 4a) \cup B(y_0, 4a) \subset \R^d_+$.

We recall that by \cite[Proposition 3.11(b)]{KSV}, $\E_y\tau_{B(x_0, 2a)}\le c_1 a^{\alpha}$ for all $y\in B(x_0,a)$. Let $N\ge 1/a$. By using  \eqref{e:levys3} in the second line and Proposition \ref{p:gfcnub} in the fourth, we have for every $y\in B(x_0,a)$, 
\begin{align*}
&\E_y\left[G(Y_{\tau_{B(x_0, 2a)}}, y_0); Y_{\tau_{B(x_0, 2a)}} \in B(y_0, 1/N) \right]\\
& = \E_y \left(\int_0^{\tau_{B(x_0, 2a)}} 
  \int_{B(y_0, 1/N)}
G(w, y_0)J(Y_s,w)dw\, ds\right)\\
& \le 
\left(\sup_{y \in B(x_0, a)}
\E_y\tau_{B(x_0, 2a)}\right) \left(\sup_{
z \in B(x_0, 2a)}
\int_{B(y_0, 1/N)} J(z, w)G(w, y_0)dw \right)\\
 & \le c_1 a^{\alpha} (8a)^{-d-\alpha} \int_{B(y_0, 1/N)} |w-y_0|^{-d+\alpha}dw =c_2 a^{-d} (1/N)^{\alpha}.
\end{align*}
Now choose $N$ large enough so that $c_2 a^{-d} (1/N)^{\alpha}<\epsilon/4$. Then
$$
\sup_{y \in B(x_0, a)}\E_y
\left[G(Y_{\tau_{B(x_0, 2a)}}, y_0); Y_{\tau_{B(x_0, 2a)}} \in B(y_0, 1/N) \right]<\epsilon/4.
$$
Since by Proposition \ref{p:gfcnub}, 
$x\mapsto h(x) 
:= \E_x\left[G(Y_{\tau_{B(x_0, 2a)}}, y_0); Y_{\tau_{B(x_0, 2a)}}\in \R^d_+ \setminus B(y_0, 1/N)\right]$ is a bounded function which is  
 harmonic on $ B(x_0, a)$, 
 it is continuous by \eqref{eHR} so we can choose a $\delta 
\in (0, a)$ such that $|h(x)-h(x_0)|<\eps/2$ for all 
$x\in B(x_0, \delta)$ ,
Therefore, for all  $x\in B(x_0, \delta)$
\begin{align*}
&|G(x, y_0)-G(x_0, y_0)| \\
\le &|h(x)-h(x_0)|+2
\sup_{y \in B(x_0, a)}\E_y
\left[G(Y_{\tau_{B(x_0, 2a)}}, y_0); Y_{\tau_{B(x_0, 2a)}} \in B(y_0, 1/N) \right] <\eps.
\end{align*}

\noindent
(1) 
Combining 
Theorem \ref{t:GB} 
and Lemma \ref{l:f4} 
with \eqref{e:scaling-of-G}, we arrive at Theorem \ref{t:Green}(1). 

\noindent
(2)-(3) Combining Theorem  \ref{t:f4e},  Proposition \ref{p:f4ef} and \eqref{e:scaling-of-G}, 
we arrive at Theorem \ref{t:Green}(2)-(3).
\qed

As an application of Theorem \ref{t:Green},
we get the following estimates on killed potentials of $Y$.

\begin{prop}\label{p:bound-for-integral-new} 
Suppose that 
$p\in ((\alpha-1)_+, \alpha+\beta_1)$.
Then
for any 
$\wt{w} \in \R^{d-1}$, any Borel set $D$ satisfying $D_{\wt{w}}(R/2,R/2) \subset D  \subset D_{\wt{w}}(R,R)$ 
and any 
$x=(\wt{w}, x_d)$ with 
$0<x_d \le R/10$,
\begin{equation}\label{e::bound-for-integral-new} 
\E_x \int_0^{\tau_D}(Y_t^d)^{\gamma }\, dt  = \int_D G^D(x,y)y_d^{ \gamma }\, dy  \asymp  
\begin{cases} 
R^{\alpha+  \gamma   -p}x_d^p, &  \gamma >p-\alpha,\\ 
x_d^p\log(R/x_d), &  \gamma =p-\alpha, \\ 
x_d^{\alpha+ \gamma }, &-p-1<  \gamma <p-\alpha, 
\end{cases}
\end{equation}
where the comparison constant is independent of $\wt{w} \in \R^{d-1}$, $D$, $R$ and $x$.
\end{prop}
\pf
Without loss of generality, we assume  $\wt w=\wt x= \wt 0$.

\noindent
(i) {\it Upper bound:}
Note that, by Lemma \ref{l:prelub}, 
\begin{align*}
&\int_D G^D(x,y)y_d^{\gamma}\, dy\le \int_{D(R,R)} G(x,y)y_d^{\gamma}\, dy\\
&\le c_0\Big( \int_{D(R,x_d/2)}y_d^{\gamma+p}|x-y|^{\alpha-d-p}dy
+\int_{D(R,R)\setminus D(R, x_d/2)}y_d^{\gamma}
\left(\frac{x_d}{|x-y|}\wedge 1\right)^p
|x-y|^{\alpha-d}dy
\Big)\\
&=c_0\Big( \int_{D(R,x_d/2)}f(y; \gamma+p, 0, p, 0, x)dy
+x_d^{\gamma} \int_{D(R,3x_d/2)\setminus D(R, x_d/2)}
g(y; 0, p, 0, x)\, dy \\
&\quad +x_d^p\int_{D(R,R)\setminus D(R, 3x_d/2)}
 f(y; \gamma, 0,p, 0, x) \, dy
\Big)=:c_0(I_1+I_2+I_3).
\end{align*}
Suppose first $-p-1<  \gamma<p-\alpha$. We use Lemma \ref{l:key}(i) on $I_1$(which is allowed since $\gamma+p>-1$) 
and Lemma \ref{l:key}(iii) on $I_2$. Then
$$
I_1 \asymp x_d^{\alpha-p-1}\left(\frac{x_d}{2}\right)^{\gamma+p+1}\asymp x_d^{\alpha+\gamma}\quad \text{and} \quad 
I_2
  \asymp x_d^{\gamma}x_d^{\alpha}=x_d^{\alpha+\gamma}.
$$
 Finally,
\begin{align*}
I_3
\asymp  x_d^p R^{\gamma+\alpha-p}F\left(\frac{3x_d}{2R}; \gamma+\alpha-p-1, 0\right) \asymp   x_d^p R^{\gamma+\alpha-p} \left(\frac{3x_d}{2R}\right)^{\gamma+\alpha-p} \asymp x_d^{\alpha+\gamma}.
\end{align*}
Here the first asymptotic equality follows from Lemma \ref{l:key}(ii) (with $a_2=R$ and $a_3=3x_d/2$) and the second asymptotic equality from the definition of $F(\cdot\ ;  \cdot, 0)$. 

This completes the proof of the upper bound in the case $-p-1<  \gamma<p-\alpha$. The other two cases are similar, but simpler, since one can directly use Corollary \ref{c:key}(i) with Lemma \ref{l:prelub}. We omit the details.

\noindent
(ii) \emph{Lower bound:} 
We first note that by Theorem \ref{t:GB} 
\begin{align*}
\E_x \int_0^{\tau_D}(Y_t^d)^{\gamma}\, dt  &\ge \int_{B^+_{R/2}}  y_d^{\gamma} 
G^{B^+_{R/2}}(x, y)\, dy 
\ge \int_{D(R/5, R/5)}   y_d^{\gamma} 
G^{B^+_{R/2}}(x, y)\, dy\nn\\
&\ge cx_d^p \int_{D(R/5, R/5)\setminus D(R/5, 3x_d/2)}    
\frac{y_d^{p+\gamma} dy} {|x-y|^{d-\alpha+2p}} \\
& =c
x_d^p \int_{D(R/5, R/5)\setminus D(R/5, 3x_d/2)}  
f(y;  \gamma+p, 0, 2p, 0, x)
dy.
\end{align*}
Since $3x_d/2 <3R/20$ (so that $3x_d/(2R/5) \le 3/4$), 
using Lemma \ref{l:key}(ii) (with $a_2=R/5$ and $a_3=3x_d/2$)
and applying
\eqref{e:F1} and \eqref{e:F3}, we immediately get the lower bound. 
\qed

\begin{remark}\label{r:I2}
{\rm
(a)
It follows from the proof of Proposition \ref{p:bound-for-integral-new} and Remark \ref{r:I1} that 
$$ \int_D G^D(x,y)y_d^{ \gamma }\, dy=\infty \quad \text{ if } \gamma  \le -p-1.$$

\noindent
(b)
By Proposition \ref{p:bound-for-integral-new},
for any $ \beta_1\ge 0$, and all $r\in (0, 1]$ and $x\in U(r)$,
\begin{align*}
&r^{\alpha+\beta_1-p}x_d^p \asymp \E_x \int_0^{\tau_{U(r)}} (Y_t^d)^{\beta_1} \, dt \le \E_x \int_0^{\tau_{U(r)}} 
(Y_t^d)^{ \beta_1} 
  |\log Y_t^d|^{ \beta_3  }\, dt\\
&\le c \E_x \int_0^{\tau_{U(r)}} (Y_t^d)^{(p-\alpha+\beta_1)/2} \, dt \asymp r^{\alpha+(p-\alpha+\beta_1)/2-p}x_d^p
=r^{(\alpha+\beta_1-p)/2}x_d^p \le x_d^p\, .
\end{align*}
Thus, Proposition \ref{p:bound-for-integral-new} is a significant generalization of Lemma \ref{l:upper-bound-for-integral}.}
\end{remark}

We end this section with the following corollary, which follows from  
Proposition \ref{p:bound-for-integral-new} and Remark \ref{r:I2} by letting $R \to \infty$.
Recall that $Y_{\zeta-}=\lim_{t\uparrow \zeta}Y_t$ denotes the left limit of the process $Y$ at its lifetime.
\begin{corollary}
Suppose that 
 $p\in ((\alpha-1)_+, \alpha+\beta_1)$.
Then for all $x \in \R^d_+$,
$$
\E_x\int_0^{\zeta} (Y_t^d)^{\gamma}\, dt =   \int_{\R^d_+}G(x,y)y_d^{ \gamma } dy
 \asymp  
\begin{cases} 
\infty & \gamma \ge p-\alpha  \text{ or } \gamma  \le -p-1,\\ 
x_d^{\alpha+\gamma }, & -p-1<   \gamma <p-\alpha.
\end{cases}
$$
In particular,  for all $x \in \R^d_+$,
$\P_x(Y_{\zeta-} \in \R_+^d,  \zeta <\infty) =G\kappa(x) \asymp c>0$
and
$$
\E_x[\zeta]
 \asymp  
\begin{cases} \infty & p \le \alpha,\\ 
x_d^{\alpha}, & p > \alpha.
\end{cases}
$$
\end{corollary}

\section{Boundary Harnack principle}\label{s:BHP}

In this section we give a proof of Theorem \ref{t:BHPnew}. We start with a lemma providing important estimates of the jump kernel $J$ needed in the proof.
Recall that $U=D(\frac{1}2, \frac{1}2).$

Note the exponent $\beta_1-\eps$ in \eqref{e:beta37-new} below is not necessarily  positive, 
but is always strictly larger than $-1$.

\begin{lemma}\label{l:estimates-of-J-for-BHP}
Suppose $p\in ((\alpha-1)_+, \alpha+(\beta_1 \wedge \beta_2))$ and let
\begin{align}
\label{e:k}
k(y)=\frac{(y_d \wedge 1)^{\beta_1}(y_d \vee1)^{\beta_2}}{|y|^{d+\alpha+\beta_1+\beta_2} } (1+|\log(y_d)|)^{\beta_3}
\left(\log\left(1+\frac{|y|}{y_d \vee 1}\right)\right)^{\beta_4}. 
\end{align}
(a) Let 
$z^{(0)}=(\widetilde{0}, 2^{-2})$.  Then for any $z\in B(z^{(0)}, 2^{-3})$ 
and $y \in \R^d_+\setminus D(1, 1)$, it holds that
\begin{align}\label{e:beta30}
J(z,y) \ge  ck(y).
\end{align}

\noindent (b) Let 
$$
 \eps=\beta_1+\alpha-p-\frac{\beta_2+\alpha-p}{M}, \quad \text{where}\quad 
M= 1+ \left( \frac{\beta_2+\alpha-p}{\beta_1+\alpha-p} \vee 1 \right).
$$
Then for any $z\in U$ and $y \in \R^d_+\setminus D(1, 1)$, it holds that
\begin{equation}\label{e:beta37-new}
J(z,y)\le  c z_d^{\beta_1-\eps}k(y).
\end{equation}
\end{lemma}

\pf (a) 
For $z\in B(z^{(0)}, 
2^{-3})$ and $y\in \R^d_+\setminus D(1, 1)$, 
$z_d\asymp z^{(0)}_d =2^{-2}$ 
and $|z-y|\asymp |z^{(0)}-y| \asymp |y| >c$ which immediately implies \eqref{e:beta30}. 

\noindent (b)  
Let $\delta=(1-\frac1M)(\beta_2+\alpha-p)>0$.  
We first note that by the definitions of $M$, $\delta$ and $\epsilon$, 
we have that
\begin{align}\label{e:beta35}
\eps > \beta_1+\alpha-p-\frac{\beta_2+\alpha-p}{ \left( \frac{\beta_2+\alpha-p}{\beta_1+\alpha-p} \vee 1 \right)}
= \beta_1+\alpha-p- (\beta_1+\alpha-p) \wedge (\beta_2+\alpha-p) \ge 0
\end{align}
and 
\begin{align}\label{e:beta34}
\beta_2+\eps
=\beta_2 + \beta_1+\alpha-p-\frac{\beta_2+\alpha-p}{M}
= \beta_1 + (1-\frac1M)(\beta_2+\alpha-p) 
= \beta_1 + \delta 
> \beta_1.
\end{align}

Assume that  $z\in U$ and $y\in \R^d_+\setminus D(1, 1)$. 
Since 
$|z-y|\asymp |y| \ge c(z_d \vee y_d)$, it holds that
\begin{align}\label{e:beta31}
J(z,y) \asymp
\frac{
(z_d\wedge y_d)^{\beta_1}(z_d\vee y_d)^{\beta_2}}
{|y|^{d+\alpha+\beta_1+\beta_2}}
\left(\log\left(1+\frac{z_d\vee y_d}{z_d\wedge y_d}\right)\right)^{\beta_3}
 \left(\log\left(1+\frac{ |y|}{(y_d \vee z_d) \wedge |y|}\right)\right)^{\beta_4}.  
\end{align}
Clearly,
if $y_d \ge 3/4>1/2 \ge z_d$, then 
 $$
 \frac{ |y|}{(y_d \vee z_d) \wedge |y|}  \asymp \frac{ |y|}{(y_d \vee 1) \wedge |y|}  \asymp\frac{ |y|}{y_d \vee 1 }  
$$
and 
\begin{align*}
&\log\left( 1+\frac{y_d}{z_d}\right) \le  
3\log\left(\frac{y_d}{z_d}\right) \le3 \left(|\log y_d|+\log\left(\frac1{z_d}\right)\right)\\
&\le   
6|\log y_d|\log\left(\frac1{z_d}\right) + 3\log\left(\frac1{z_d}\right) \le 6\log\left(\frac1{z_d}\right)(1+ |\log y_d|).
\end{align*}
Thus, for $z\in U$ and $y\in \R^d_+\setminus D(1, 1)$ with $y_d \ge 3/4$,
\begin{align}\label{e:beta322}
J(z,y) \asymp
 \frac{z_d^{\beta_1}y_d^{\beta_2}}{|y|^{d+\alpha+\beta_1+\beta_2} } 
\left(\log \left( \frac{y_d}{z_d} \right) \right)^{\beta_3}
\left(\log\left(1+\frac{|y|}{y_d \vee 1}\right)\right)^{\beta_4}
 \le c
z_d ^{\beta_1}  
\left(\log\left(\frac1{z_d}\right)\right)^{\beta_3}
k(y).
\end{align}
It is easy to see from \eqref{e:beta31} that for 
$(z, y)\in U\times (\R^d_+\setminus D(1, 1))$ with $y_d <3/4$ and $z_d>y_d$, 
\begin{align*}
J(z,y)
\le c\frac{y_d^{\beta_1}z_d^{\beta_2}}{|y|^{d+\alpha+\beta_1+\beta_2} } 
\left(\log\left(\frac1{y_d}\right)\right)^{\beta_3}
\left(\log\left(\frac{|y|}{z_d}\right)\right)^{\beta_4}.
\end{align*}
Since $\delta>0$, 
we have
\begin{align}
&z_d^\delta 
\left(\log\left(\frac{|y|}{z_d}\right)\right)^{\beta_4} 
=|y|^\delta \left(\frac{z_d}{|y|}\right)^\delta 
\left(\log\left(\frac{|y|}{z_d}\right)\right)^{\beta_4} 
\nn\\
\le&
c |y|^\delta \left(\frac{2^{-1}}{|y|}\right)^\delta 
\left(\log\left(\frac{|y|}{2^{-1}}\right)\right)^{\beta_4} 
\le
 c  
\left(\log(2|y|)\right)^{\beta_4}, 
\quad  0<z_d  \le 1/2 <1<|y|. \label{e:nbeta1}
\end{align}
Thus, using \eqref{e:beta34}
\begin{align}\label{e:beta32}
J(z,y)
\le c\frac{y_d^{\beta_1}z_d^{\beta_2}}{|y|^{d+\alpha+\beta_1+\beta_2} } 
\left(\log\left(\frac1{y_d}\right)\right)^{\beta_3}
\left(\log\left(\frac{|y|}{z_d}\right)\right)^{\beta_4}
\le c z_d^{\beta_2-\delta} 
k(y)=c z_d^{\beta_1-\eps} 
k(y).
\end{align}

Since $\eps>0$ by \eqref{e:beta35}, 
we have
$$
z_d^\eps  
\left(\log\left(\frac1{z_d}\right)\right)^{\beta_3} 
\le c y_d^\eps  
\left(\log\left(\frac1{y_d}\right)\right)^{\beta_3}, 
\quad  0<z_d \le y_d<3/2,
$$
so that  by using the same argument  as in \eqref{e:nbeta1},  
\begin{align*}
&z_d^\eps  
\left(\log\left(\frac1{z_d}\right)\right)^{\beta_3} \left(\log\left(\frac{|y|}{y_d}\right)\right)^{\beta_4}
\le c y_d^{\eps-\delta}
\left(\log\left(\frac1{y_d}\right)\right)^{\beta_3} 
 y_d^\delta  \left(\log\left(\frac{|y|}{y_d}\right)\right)^{\beta_4}\\
& \le c y_d^{\eps-\delta}
\left(\log\left(\frac1{y_d}\right)\right)^{\beta_3} 
\left(\log(2|y|)\right)^{\beta_4}, 
\quad  0<z_d \le y_d<3/2<1<|y|.
\end{align*}
Thus using \eqref{e:beta34} in the last inequality below, 
we have that,  for 
$(z, y)\in U\times (\R^d_+\setminus D(1, 1))$ with $y_d <3/4$ and $z_d \le y_d$,  
\begin{align}\label{e:beta36}
&J(z,y) 
\le c\frac{z_d^{\beta_1}y_d^{\beta_2}}{|y|^{d+\alpha+\beta_1+\beta_2} } 
\left(\log\left(\frac1{z_d}\right)\right)^{\beta_3}\left(\log\left(\frac{|y|}{y_d}\right)\right)^{\beta_4}\nn\\
&=
cz_d^{\beta_1-\eps}\frac{y_d^{\beta_2}}{|y|^{d+\alpha+\beta_1+\beta_2} }z_d^\eps  
\left(\log\left(\frac1{z_d}\right)\right)^{\beta_3}\left(\log\left(\frac{|y|}{y_d}\right)\right)^{\beta_4}\nn\\
&\le c  z_d^{\beta_1-\eps}
\frac{y_d^{\beta_2+\eps-\delta}}{|y|^{d+\alpha+\beta_1+\beta_2} } 
\left(\log\left(\frac1{y_d}\right)\right)^{\beta_3}\left(\log(2|y|)\right)^{\beta_4}
\le c z_d^{\beta_1-\eps}
k(y).
\end{align}
Combining
\eqref{e:beta322}, \eqref{e:beta32} and \eqref{e:beta36}, 
and using the inequality 
$$
z_d^{\beta_1-\eps}  \vee (z_d^{\beta_1}  
(\log(1/z_d))^{\beta_3}) 
\le c z_d^{ \beta_1-\epsilon}, \quad z \in U,
$$ 
we get the upper bound \eqref{e:beta37-new} for $J(z,y)$.
\qed

\bigskip

\noindent {\bf Proof of Theorem \ref{t:BHPnew}}.
By scaling, we just need to consider the case $r=1$. Moreover, 
by
Theorem  \ref{t:uhp} (b), 
it suffices to prove \eqref{e:TAMSe1.8new}
for  $x, y\in D_{\wt w}(2^{-8}, 2^{-8})$.

Since $f$ is harmonic in $D_{\wt w}(2, 2)$ and vanishes continuously on 
$B((\wt w, 0), 2)\cap \partial \R^d_+$,
it is regular harmonic in $D_{\wt w}(7/4, 7/4)$ and vanishes continuously on 
$B((\wt w, 0), 7/4)\cap \partial \R^d_+$.
Throughout the remainder of this proof, we assume that 
$x\in D_{\wt w}(2^{-8}, 2^{-8})$. 
Without loss of generality we take $\widetilde{w}=0$.

Define 
$z^{(0)}=(\widetilde{0}, 2^{-2})$. By Theorem  \ref{t:uhp} (b)
and Lemma \ref{l:exit-probability-estimate}, we have
\begin{align}\label{e:TAMSe6.37-new}
f(x)&=\E_x[f(Y_{\tau_{
U}})]\ge \E_x[f(Y_{\tau_{U}}); Y_{\tau_{U}}\in D(1/2, 1)\setminus D(1/2, 3/4)]\nonumber\\
&\ge c_1f(z^{(0)})\P_x(Y_{\tau_{D_{\wt x}(1/4, 1/4)}}\in D_{\wt x}(1/4, 1)\setminus 
D_{\wt x}(1/4, 3/4))\ge c_2 f(z^{(0)})x^p_d.
\end{align}

Let $k$ be the function defined in \eqref{e:k}. 
Using \eqref{e:beta30}, 
the harmonicity  of $f$, the L\'evy system formula and \cite[Proposition 3.11(a)]{KSV}, 
\begin{align}\label{e:POTAe7.27}
&f(
z^{(0)})\ge \E_{z^{(0)}} \left[f(Y_{\tau_{U}}); Y_{\tau_{U}}\notin D(1, 1)\right]\nonumber\\
&\ge \E_
{z^{(0)}}\int^{\tau_{
B(z^{(0)}, 2^{-3})}
}_0\int_{\R^d_+\setminus D(1, 1)}J(Y_t, y)f(y)dydt\nonumber\\
&\ge c_{10}
\E_{z^{(0)}}\tau_{B(z^{(0)}, 2^{-3})}
\int_{\R^d_+\setminus D(1, 1)}
k(y)
f(y)dy\ge c_{11}\int_{\R^d_+\setminus D(1, 1)}
k(y)f(y)dy.
\end{align}

Now we assume that  $z\in U$ and $y\in \R^d_+\setminus D(1, 1)$. Let $\epsilon $ be defined as in Lemma \ref{l:estimates-of-J-for-BHP}.
Since  
$\beta_1-\eps > \beta_1-(\alpha+\beta_1-p)= 
p-\alpha$, 
 by Proposition \ref{p:bound-for-integral-new}
and \eqref{e:beta37-new}, we have
\begin{align}\label{e:4:09}
&\E_x\left[f(Y_{\tau_{U}}); Y_{\tau_{U}}\notin D(1, 1)\right]=\E_x\int^{\tau_{U}}_0\int_{\R^d_+\setminus D(1, 1)}J(Y_t, y)f(y)dydt\nonumber\\
&\le c\E_x\int_0^{\tau_{U}}(Y_t^d)^{ \beta_1-\epsilon } dt\int_{\R^d_+\setminus D(1, 1)} 
k(y)f(y)dy\le cx_d^p\int_{\R^d_+\setminus D(1, 1)} 
k(y)f(y)dy.
\end{align}
Combining this with \eqref{e:POTAe7.27},  we now have
\begin{align}\label{e:POTAe7.29}
&\E_x\left[f(Y_{\tau_{U}}); Y_{\tau_{U}}\notin D(1, 1)\right]\le c x^p_df(w).
\end{align}

On the other hand, 
since $f$ is a non-negative function  in $\R^d_+$ which is harmonic in $D_{\wt{w}}(2, 2)$ with respect to $Y$ and vanishes continuously on 
$B(({\wt{w}},0), 2)\cap \partial \R^d_+$, by Theorem  \ref{t:uhp} (b) and  
Carleson's estimate (Theorem \ref{t:carleson}), it holds that $f(v)\le c_{16}f(z^{(0)})$ for all $v \in D(1, 1)$.
Therefore, by Lemma \ref{l:POTAl7.4}, we have
\begin{align}\label{e:POTAe7.30}
\E_x\left[f(Y_{\tau_{U}}); Y_{\tau_{U}}\in D(1, 1)\right]&\le c_{16}f(z^{(0)})\P_x\left(Y_{\tau_{U}}\in D(1, 1)\right)\le c_{17}f(z^{(0)})x^p_d.
\end{align}
Combining \eqref{e:POTAe7.29}, \eqref{e:POTAe7.30} and \eqref{e:TAMSe6.37-new} we get that
$f(x)\asymp x^p_df(z^{(0)})$ for all $x\in  D(2^{-8}, 2^{-8})$, 
which implies that 
 that for all 
$x, y\in  D(2^{-8}, 2^{-8})$,
$$
\frac{f(x)}{f(y)}\le c_7\frac{x^p_d}{y^p_d},
$$
which is same as the conclusion of the theorem.
\qed

\medskip
\noindent {\bf Proof of Theorem \ref{t:counterexample}}.
The case 
$\alpha+\beta_2<p<\alpha+\beta_1$
 has been dealt with in \cite[Theorem 1.4.]{KSV}, so we only need to deal with the case $p=\alpha+\beta_2$. In the rest of the proof, we assume $p=\alpha+\beta_2$. The proof is the same as that of \cite[Theorem 1.4.]{KSV}, except  that we now can use 
 Proposition \ref{p:bound-for-integral-new} to get for all 
$r>0$ and $x=(\wt 0, x_d)$ with $0<x_d \le r/10$,
\begin{align}
\label{e:newcounter00}
\E_x\int_{0}^{\tau_{U(r)}}  (Y_t^d)^{\beta_2}  dt  \asymp   x_d^{\beta_2+ \alpha} \log(r/x_d)
= x_d^{p} \log(r/x_d). 
\end{align}
Moreover, using \eqref{e:newcounter00}, we also get 
that for every $r>0$ and $x \in U(r)$,
\begin{align}
\label{e:newcounter001}
\E_x\int_{0}^{\tau_{U(r)}}  (Y_t^d)^{\beta_2}  dt  
\le \E_x\int_{0}^{\tau_{D_{\wt x}(5r, 5r)} } (Y_t^d)^{\beta_2}  dt  
\le c_0 x_d^{p} \log(r/x_d). 
\end{align}
The displays \eqref{e:newcounter00} and \eqref{e:newcounter001} will be used to replace the roles played by \cite[Lemmas 5.11 and 5.12]{KSV}.
We provide the full proof for the convenience of the reader.

Suppose that the non-scale invariant BHP holds near the boundary of $\R^d_+$
(see the paragraph before Theorem \ref{t:counterexample}).

Note that by taking 
$g(x)= \P_x(Y_{\tau_{U}}\in D(1/2,1)\setminus D(1/2,3/4))$, 
we see from 
 Lemma \ref{l:exit-probability-estimate} that 
there exists $\wh{R} \in (0, 1)$ such that for any $r \in (0, \wh{R}\, ]$ 
there exists a constant $c_1=c_1(r)>0$ 
such that for any non-negative function $f$ in $\R^d_+$ which is harmonic  in $\R^d_+ \cap B(0, r)$ with respect to $Y$ and vanishes continuously on $ \partial \R^d_+ \cap B(0, r)$, 
\begin{equation}\label{e:bhp_mfx}
\frac{f(x)}{f(y)}\,\le \, c_1 \,\frac{x_d^p}{y_d^p}, \quad  \hbox{for all } x, y\in  
\R^d_+ \cap B(0, r/2).
\end{equation}

Let $r_0=\wh{R}/4$  and choose a point $z_0 \in \partial \R^d_+$  with $|z_0| = 4$.
For  $n \in \N$, $B(z_0,1/n)$ does not intersect $B(0, 2r_0)$. We define
$$
K_n:=\int_{ \R^d_+   \cap B(z_0, 1/n)} |\log(y_d)|^{\beta_3  +\beta_4 }dy\,  ,\quad  \quad f_n(y):=
K_n^{-1}y_d^{-\beta_1} {\bf 1}_{\R^d_+ \cap B(z_0, 1/n)}(y),
$$
and 
$$
g_n(x):=\E_x\left[f_n(Y_{\tau_{U(r_0)}})\right]=\E_x\int_{0}^{\tau_{U(r_0)}}  \int_{\R^d_+ \cap B(z_0, 1/n)} 
J(Y_t,y)f_n(y) dydt, \quad x \in {U(r_0)}.
$$

We claim that there exists $c_2>0$ 
such that 
\begin{align}
\label{e:gnlow}
\liminf_{n \to \infty}g_n(x)  
 \ge c_2 x_d^{\beta_2+\alpha}\log(r_0/x_d)= c_2 x_d^{p}\log(r_0/x_d)
\end{align}
for all  $x=x^{(s)}=(\wt 0, s)\in \R^d_+$
with 
$s \in (0,  r_0/10)$. 

Here is a proof of the claim above.
Since 
$$
 6 >   |z-y|  >2 >y_d \wedge z_d \quad 
\text{ for } (y, z)  \in \left(\R^d_+ \cap B(z_0, 1/n)\right) \times   {U(r_0)},
$$
using 
\textbf{(A3)} 
we have for $(y, z)  \in \left(\R^d_+  \cap B(z_0, 1/n)\right) \times   {U(r_0)}$,
\begin{align*}
J(z,y)&\asymp
(z_d \wedge y_d)^{\beta_1}  (z_d \vee y_d)^{\beta_2}  
\left( \log\left(1+\frac{z_d \vee y_d}{z_d \wedge y_d}\right)\right)^{\beta_3}  \left(\log\left( \frac{1}{z_d \vee y_d}\right)\right)^{\beta_4}\\
&\asymp
 \frac{ z_d^{\beta_1} y_d^{\beta_1} }{ (z_d \vee y_d)^{\beta_1-\beta_2}} 
\left(   \log\left(1+\frac{z_d \vee y_d}{z_d \wedge y_d}\right)\right)^{\beta_3} \left(  \log\left( \frac{1}{z_d \vee y_d}\right)\right)^{\beta_4} .
\end{align*}
Therefore, for $ x \in {U(r_0)}$,
\begin{align}
\label{e:jnew2}
g_n(x)\asymp 
& 
K_n^{-1}
\E_x\int_{0}^{\tau_{U(r_0)}}  (Y_t^d)^{\beta_1} \int_{\R^d_+ \cap B(z_0, 1/n)} 
(Y_t^d \vee y_d)^{-(\beta_1-\beta_2)}
\nn\\
&\quad \times \left(\log\left(1+\frac{Y_t^d \vee y_d}{Y_t^d \wedge y_d}\right)\right)^{\beta_3} \left( \log\left(\frac{1}{Y_t^d \vee y_d}\right)\right)^{\beta_4}    dydt.
\end{align}
Note that, using 
$\sup_{ t\ge 1}t^{-(\beta_1-\beta_2)} (\log(1+t))^{\beta_3}<\infty$, 
 for $z\in U(r_0)$,
\begin{align}
\label{e:jnew12} 
 &  K_n^{-1}
 \int_{\R^d_+ \cap B(z_0, 1/n)}( z_d \vee y_d)^{-(\beta_1-\beta_2)}  
 \left( \log\left(1+\frac{z_d \vee y_d}{z_d \wedge y_d}\right)\right)^{\beta_3} 
\left(  \log\left(\frac{1}{z_d \vee y_d}\right)\right)^{\beta_4}
 dy  \nn\\
\le & 
\frac{K_n^{-1}}
{{z_d}^{\beta_1-\beta_2}}  \int_{\R^d_+ \cap B(z_0, 1/n)\cap \{z_d \le y_d\}} (z_d/y_d)^{\beta_1-\beta_2}\
\left( \log\left(1+\frac{y_d}{z_d}\right)\right)^{\beta_3}  
\left( \log\left(\frac{1}{ y_d}\right)\right)^{\beta_4}  dy\nn\\
&+ \frac{K_n^{-1}}
{{z_d}^{\beta_1-\beta_2}} 
 \int_{\R^d_+ \cap B(z_0, 1/n)\cap \{z_d > y_d\}  
}
\left(\log\left(1+\frac{z_d}{y_d}\right)\right)^{\beta_3}
\left(   \log\left(\frac{1}{z_d }\right)\right)^{\beta_4} dy \nn\\
 \le& c_3 \frac{K_n^{-1}}
{{z_d}^{\beta_1-\beta_2}} 
 \int_{\R^d_+ \cap B(z_0, 1/n)}  
 \left( \log\left(\frac{1}{y_d}\right) \right)^{\beta_3  +\beta_4 }
 dy   \le c_4 z_d^{-(\beta_1-\beta_2)}   
\end{align} 
and
 \begin{align*}
& \lim_{n \to \infty}
 K_n^{-1}
 \int_{\R^d_+  \cap B(z_0, 1/n)} ( z_d \vee y_d)^{-(\beta_1-\beta_2)} \left(  \log\left(1+\frac{z_d \vee y_d}{z_d \wedge y_d}\right)\right)^{\beta_3}   \left( \log\left(\frac{1}{z_d \vee y_d}\right)\right)^{\beta_4}  dy \\
 &=
  z_d^{-(\beta_1-\beta_2)}\, .
 \end{align*}
Moreover, by \eqref{e:newcounter00},
$\E_x\int_{0}^{\tau_{U(r_0)}}  (Y_t^d)^{\beta_2}dt  < \infty$ for all 
$x\in U(r_0)$. 
 Thus we can use  the dominated convergence theorem to get that  for all  
 $x\in U(r_0)$,
 \begin{align}\label{e:notD231}
&\lim_{n \to \infty} K_n^{-1}
\E_x\int_{0}^{\tau_{U(r_0)}}  (Y_t^d)^{\beta_1} \int_{\R^d_+ \cap B(z_0, 1/n)} 
 (Y_t^d \vee z_d)^{-(\beta_1-\beta_2)} 
\nn\\
&\qquad\qquad \times  \left(   \log\left(1+\frac{Y_t^d \vee y_d}{Y_t^d \wedge y_d}\right)\right)^{\beta_3} \left(  \log\left(\frac{1}{Y_t^d \vee y_d}\right)\right)^{\beta_4}  dydt \nn \\
=&\E_x\int_{0}^{\tau_{U(r_0)}}  (Y_t^d)^{\beta_1} (Y_t^d)^{-(\beta_1-\beta_2)}dt
=  \E_x\int_{0}^{\tau_{U(r_0)}}  (Y_t^d)^{\beta_2}dt. 
\end{align} 
Combining
\eqref{e:notD231} with \eqref{e:newcounter00}
 we conclude that \eqref{e:gnlow} holds true.
 
 From \eqref{e:jnew2}, \eqref{e:jnew12} and
 \eqref{e:newcounter001}
 we see 
that for all  $x\in U(r_0)$,
\begin{align}
\label{e:jnew3}
 g_n(x)
 &\le  c_5
 K_n^{-1}
\E_x\int_{0}^{\tau_{U(r_0)}}  (Y_t^d)^{\beta_1} \int_{\R^d_+ \cap B(z_0, 1/n)} 
 (Y_t^d)^{-(\beta_1-\beta_2)}
 \nn \\
&  \qquad \times   \left( \log\left(1+\frac{Y_t^d \vee y_d}{Y_t^d \wedge y_d}\right)\right) ^{\beta_3}   \left(\log\left(\frac{1}{Y_t^d \vee y_d}\right)\right)^{\beta_4}  dydt \nn\\
& 
\le c_6 \E_x\int_{0}^{\tau_{U(r_0)}}  (Y_t^d)^{\beta_2}  dt \le 
 c_7 x_d^{p} \log(r_0/x_d)
 . \end{align}
Thus the $g_n$'s are non-negative functions  in $\R^d_+ $ 
which are harmonic  in $\R^d_+  \cap B(0,2^{-2}  r_0)$ with respect to 
$Y$ and vanish continuously on $ \partial \R^d_+  \cap B(0,2^{-2}  r_0)$. Therefore, by \eqref{e:bhp_mfx},
$$
\frac{g_n(y)}{g_n(w)} \le 
c_1 \frac{y_d^p }{w_d^p }  \quad 
\text{for all } y\in D \cap B(0,2^{-3}  r_0),
$$
where $w=(\wt{0}, 2^{-3} r_0)$ and 
$c_1=c_1(2^{-2}  r_0)$. 
Thus by \eqref{e:jnew3}, for all 
$y\in \R^d_+  \cap B(0,2^{-3}  r_0)$,
$$
\limsup_{n \to \infty}g_n(y) \le 
c_1 \limsup_{n \to \infty}g_n(w) \frac{y_d^p }{w_d^p } \le c_8 y_d^p.
$$
This and \eqref{e:gnlow} imply that  for all  
$x=x^{(s)}=(\wt 0, s)\in \R^d_+$
with $s \in (0, r_0/10)$, $x_d^{p} \log(r_0/x_d) \le c_9 x_d^{p} $, which  gives 
a contradiction. 
\qed

\vspace{.1in}
\textbf{Acknowledgment}:
We thank the anonymous referees 
for very helpful comments which lead to improvements of this paper.
Part of the research for this paper was done while the second-named author was visiting
Jiangsu Normal University, where he was partially supported by  a grant from
the National Natural Science Foundation of China (11931004) and by
the Priority Academic Program Development of Jiangsu Higher Education Institutions

\vskip 0.1truein

\parindent=0em

{\bf Panki Kim}

Department of Mathematical Sciences and Research Institute of Mathematics,

Seoul National University, Seoul 08826, Republic of Korea

E-mail: \texttt{pkim@snu.ac.kr}

\bigskip

{\bf Renming Song}

Department of Mathematics, University of Illinois, Urbana, IL 61801,
USA

E-mail: \texttt{rsong@illinois.edu}

\bigskip

{\bf Zoran Vondra\v{c}ek}

Department of Mathematics, Faculty of Science, University of Zagreb, Zagreb, Croatia,

Email: \texttt{vondra@math.hr}

\end{document}